%% file: paper-cluster-sa.tex
\documentclass[a4paper]{elsarticle_arxiv}

\usepackage[numbers]{natbib}
\usepackage[a4paper, left=30mm, right=30mm, top=25mm, bottom=25mm]{geometry}

\usepackage[utf8]{inputenc}
\usepackage[english]{babel}

\usepackage{amssymb}
\usepackage{amsbsy}
\usepackage{amsmath}
\usepackage{amsthm}
\usepackage{amsfonts}
\usepackage{accents}
\usepackage{xcolor}
\usepackage{esint}
\usepackage{todonotes}
\usepackage{xspace}
\usepackage{hyperref}

\usepackage{algpseudocode}
\usepackage{algorithm}

\input macros.tex

\newdimen\figscale
\figscale=\textwidth

\begin{document}

\begin{frontmatter}
    
\title{Reduced-order computational homogenization for hyperelastic media using gradient based sensitivity analysis of microstructures} 

\author[1]{Vladim\'\i{}r Luke\v{s}}
\ead{vlukes@kme.zcu.cz}

\address[1]{Department of Mechanics \& NTIS -- New Technologies for Information Society Research Centre,
Faculty of Applied Sciences, University of West Bohemia in Pilsen,
{Univerzitn\'\i~8}, {Plze\v{n}}, {30100}, {Czech Republic}}

\author[1]{Eduard Rohan\corref{cor1}}
\ead{rohan@kme.zcu.cz}

\cortext[cor1]{Corresponding author}

\begin{abstract}
    We propose an algorithm for the computational homogenization of locally
    periodic hyperelastic structures undergoing large deformations due to
    external quasi-static loading. The algorithm performs clustering of
    macroscopic deformations into subsets called ``centroids'', and, as
    a new ingredient, approximates the homogenized coefficients using
    sensitivity analysis of micro-configurations with respect to the macroscopic
    deformation. The novel ``model-order reduction'' approach significantly
    reduces the number of microscopic problems that must be solved in nonlinear
    simulations, thereby accelerating the overall computational process. The
    degree of reduction can be controlled by a user-defined error tolerance
    parameter. The algorithm is implemented in the finite element framework
    \SfePy, and its performance effectiveness is demonstrated using
    two-dimensional test examples, when compared with solutions obtained by the
    proper orthogonal decomposition method, and by the full ``FE-square''
    simulations. Extensions beyond the present implementations and the scope of
    tractable problems are discussed.
\end{abstract}


\end{frontmatter}

\input part_intro.tex

\input part_main.tex
\input part_num.tex
\input part_conclusion.tex

\input statements.tex
\section*{Appendix}\label{sec:appendix}
\input appendix_sa.tex

\bibliographystyle{elsarticle-num}
\bibliography{biblio}

\end{document}

%% file: macros.tex
\algnewcommand\algorithmicnot{\textbf{not}}
\algdef{SE}[IF]{IfNot}{EndIf}[1]{\algorithmicif\ \algorithmicnot\ #1\ \algorithmicthen}{\algorithmicend\ \algorithmicif}%
\algdef{SE}[DOWHILE]{Do}{doWhile}{\algorithmicdo}[1]{\algorithmicwhile\ #1}%

\def\veps{\varepsilon}
\def\sveps{^\veps}
\def\vphi{\varphi}
\def\Om{\Omega}
\def\Ome{\Om\sveps}
\def\sigmabf{\boldsymbol{\sigma}}
\def\taubf{\boldsymbol{\tau}}
\def\omegabf{\boldsymbol{\omega}}
\def\Phibf{\boldsymbol{\Phi}}
\def\phibf{\boldsymbol{\phi}}
\def\vphibf{\boldsymbol{\varphi}}
\def\Pibf{\boldsymbol{\Pi}}
\def\Xibf{\boldsymbol{\Xi}}
\def\xibf{\boldsymbol{\xi}}
\def\ab{\boldsymbol{a}}
\def\bb{\boldsymbol{b}}
\def\eb{\boldsymbol{e}}
\def\fb{\boldsymbol{f}}

\def\gb{\boldsymbol{g}}
\def\hb{\boldsymbol{h}}

\def\nb{\boldsymbol{n}}
\def\ub{\boldsymbol{u}}
\def\vb{\boldsymbol{v}}
\def\wb{\boldsymbol{w}}
\def\yb{\boldsymbol{y}}
\def\Ib{\boldsymbol{I}}
\def\Fb{\boldsymbol{F}}

\def\Rb{\boldsymbol{R}}
\def\Ub{\boldsymbol{U}}
\def\Vb{\boldsymbol{V}}
\def\qbm{\boldsymbol{\mathrm{q}}}
\def\Ubm{\boldsymbol{\mathrm{U}}}
\def\Vbm{\boldsymbol{\mathrm{V}}}
\def\Kbm{\boldsymbol{\mathrm{K}}}
\def\qbm{\boldsymbol{\mathrm{q}}}
\def\fbm{\boldsymbol{\mathrm{f}}}
\def\dub{\delta\ub}
\def\du{\delta{u}}
\def\Vbzero{\Vb\kern-0.25em_0}
\def\acirc#1{\mathring{#1}}

\def\intOm{\int\limits_{\Om}}

\def\intOmb{\int\limits_{\partial_\sigma\Om}}

\def\intOmet{\int\limits_{\Ome(t)}}

\def\intOmebt{\int\limits_{\partial_\sigma\Ome(t)}}
\def\intOmek{\int\limits_{\Om\tke}}
\def\intOmebk{\int\limits_{\partial_\sigma\Om\tke}}
\def\dOm{\,\mbox{d}\Om}
\def\dS{\,\mbox{d}S}
\def\dt{\delta t}
\def\dtau{\delta_\tau}
\def\dvel{{\boldsymbol{\mathcal{V}}}}

\def\Lie#1{\mathsterling_#1\,}

\def\Dop{\pmb{\mathbb{D}}}
\def\Aop{\pmb{\mathbb{A}}}
\def\tk{^{(k)}}
\def\tkp{^{(k + 1)}}
\def\tke{^{\veps,(k)}}
\def\tkpe{^{\veps,(k + 1)}}
\def\hx{\hat x}

\def\ihset#1{\mathcal{I}^{#1}}
\def\Acal{\mathcal{A}}
\def\Scal{\mathcal{S}}
\def\Ahom{{\boldsymbol{\Acal}}}
\def\Shom{{\boldsymbol{\Scal}}}
\def\RR{{\mathbb{R}}}

\def\eq#1{(\ref{#1})}

\DeclareMathOperator*{\argmin}{arg\,min}

\def\Tuf#1{{\mathcal{T}}_\veps{\left ({#1}\right )}}
\def\avint{\fint}
\def\intY{\avint_Y}
\def\intYi#1{\avint_{#1}}

\def\Hdb{{\bf{H}}^1}
\def\HpbY{{\bf H}_\#^1(Y)}
\def\HpbYi#1{{\bf H}_\#^1(#1)}
\def\Vbzero{\Vb\kern-0.25em_0}

\def\alin#1#2{a_{Y} \left(#1,#2\right)}
\def\alinY#1#2#3{a_{{#3}} \left(#1,#2\right)}

\def\Fcal{{\mathcal{F}}}
\def\Fcalbf{{\boldsymbol{\mathcal{F}}}}
\def\dlt{\delta}
\def\dltsh{\delta_\tau}

\newcommand\what[1]{\widehat{{#1}}}

\def\Vcal{\mathcal{V}}
\def\Vcalbf{\boldsymbol{\mathcal{V}}}
\def\pd{\partial}
\def\pdiff#1#2{\frac{\partial {#1}}{\partial {#2}}}

\def\ddt#1#2{\frac{{\rm{d}}{#1}}{{\rm{d}}{#2}}}

\def\dVy{\,\mbox{d}Y}

\def\acron{CSA\xspace}
\def\dtss{FE$^2$\xspace}
\def\dtssx{\rm{FE}^2\xspace}
\def\SfePy{{\it SfePy}\xspace}
\def\ij{_{ij}}
\def\kl{_{kl}}
\def\cdivY{\nabla_y\cdot}
\def\algas{:=}
\def\new{\mathrm{NEW}}
\def\snew{^\new}
\newcommand\wrt{{w.r.t.{~}}}
\newcommand\rhs{{r.h.s.{~}}}
\newcommand\ie{{i.e.~}}
\newcommand\eg{{e.g.~}}
\newcommand\cf{{{cf.~}}}
\definecolor{dgreen}{rgb}{0.0,0.4,0.11}
\definecolor{dorange}{rgb}{0.8,0.01,0.15}
\definecolor{dred}{rgb}{0.8,0.,0.}

\definecolor{lgray}{rgb}{0.3686, 0.5255, 0.6235}

\newcommand\itMa[1]{_{[{#1}]}}
\newcommand\itmi[1]{^{\langle{#1}\rangle}}

\newtheorem{remark}{\it Remark\/ \rm}

%% file: part_intro.tex
\section{Introduction}\label{sec:intro}

The \emph{computational homogenization (CH)}, as an alternative to the
phenomenological macroscopic modelling, see \eg \cite{Charalambakis_2010,Zohdi-Book2017},
belongs to the most relevant modelling methods for analyses of nonlinear
problems in continuum mechanics when the model can be characterized at multiple
scales, whereby at least  a ``microscopic'' scale can be described in terms of
continuum models defined in locally periodic structures using
geometrical features which are related to microstructure heterogeneities. In
general, the CH replaces the solution of a problem defined at the scale of
distinguished locally periodic heterogeneities by solutions obtained for a
smeared continuum occupied by an equivalent homogeneous material. This is characterized
by effective material coefficients which are obtained by solving appropriate
boundary value problems in the \emph{representative volume elements (RVE)}, often called the ``reference periodic unit cell'' in the homogenization community. Besides the CH based on ``numerical testing'' the RVE of
the microstructure, by applying the ``macroscopic strain modes'', also the
``asymptotic analysis based homogenization'', \cf
\cite{Cioranescu_1999, Cioranescu_2018}, leads to analogous CH algorithms,
featured by the similar hurdles arising due to the two-scale coupling
between the global ``macroscopic'' and the local ``microscopic'' nonlinear
responses.


For linear problems, especially when perfectly periodic media are considered, the
homogenization significantly reduces the computational cost of solving the
global macroscopic problem for a given finite heterogeneity size, \ie when the
scale parameter $\veps = \veps_0 > 0$ (the one involved in the asymptotic analysis of homogenization) is nonvanishing, with given the microscopic  and the macroscopic characteristic lengths, denoted by $\ell$ and $L$, respectively, such that $\veps_0 =  \ell/L$. 
By virtue of the asymptotic homogenization result, the effective medium properties, like the elasticity tensor involved in the macroscopic problem, are expressed by 
the characteristic responses of the \emph{representative periodic cell (RPC)} are solved
independently of the macroscopic problem. 
To get refined response at the heterogeneity scale, the
macroscopic solution (such as macroscopic strains) can be combined ``point-wise'' with the characteristic responses, to obtain the micro-scale field (the displacements) at the local RVE, defined as  a copy of the RPC. In this way, the macroscopic domain can be covered by a ``tiling'' so that the refined solution is obtained in the whole structure, \cf \cite{Lukes_2022} where this issue was discussed in the large deformation setting. This so-called
reconstruction step of the two-scale modelling is optional and can also be done only selectively
at locations of interest within the macroscopic domain, see \eg
\cite{Rohan-Nguyen-Naili-CaS-2020-dp-validation} for an application in
porous media modelling. In the context of small deformations and linear material
models, computational homogenization is a well-established and widely used
method for modeling, \eg elastic structures \cite{Oliveira_2009}, viscoelastic and porous media \cite{Terada_1998, Rohan_2012,Rohan_2013, Rohan-Nguyen-Naili-CaS-2020-dp-validation}, as well as
thermo-mechanical processes \cite{Terada_2010}. The idea of locally refined response has been developed in \cite{Oden-Zohdi1997} on the basis of the so-called Homogenized Dirichlet Projection Method. Recently, an efficient method for the global-local analysis related to linear higher-order gradient media was proposed in \cite{Wangermez-Allix-2020}.

In nonlinear problems, see \eg \cite{Geers_2017}, the effective material
parameters depend locally on the macroscopic state, and the full-scale
decoupling that is applicable to linear models cannot be employed. The
finite element method is applied for the spatial discretization at both the
macroscopic and the microscopic scales, such that the so approximated
macroscopic problem formulation and the ``microscopic'' problems imposed in the
local RVEs, are coupled due to the nonlinear behavior. In this context, within
the framework of the so-called \emph{FE-square} (\dtss) computational strategy \cite{Feyel_2003, Schroder_2014, Yvonnet_2007}, 
a microscopic problem is usually
solved at each quadrature point of the macroscopic domain, providing homogenized
parameters that are valid locally only at that macroscopic point. Although the
\dtss strategy can provide accurate solutions of the two-scale problem, it
becomes extremely computationally demanding, as it requires repeated solutions
of the coupled micro-macro problems within many steps of the iterative
algorithm. An extreme computational effort is observed mainly in
three-dimensional problems with complex macroscopic geometries, where achieving
sufficient accuracy of the macroscopic solution requires a large number of
finite elements in the spatial discretization or higher-order approximations,
both of which increase the number of microscopic problems that must be solved.

To bypass this computational cost, numerous approximation strategies have been
proposed (\eg \cite{Temizer_2006, Yvonnet_2007, FRITZEN_2018} for hyperelastic
structures and \cite{HERNANDEZ_2014, Gue_2024} for elasto-plastic materials),
most of which are based on the \emph{proper orthogonal decomposition (POD)} which
enables an efficient reduction of the discretized microscopic subproblems dimension
using a precomputed reduced basis. However, as the number of degrees of freedom
of the microscopic discretization grows, the computation of this basis also
becomes increasingly expensive, because it requires solving an eigenvalue
problem for large dense systems. As an alternative to POD-based approaches,
methods relying on the approximation of the homogenized parameters using trained
artificial neural networks have also been developed \cite{Le_2015_neural},
employing convolutional \cite{Fadi_2023, Eidel_2023_neural} or recurrent
\cite{Kovac_2025_neural} network architectures. Furthermore, in \cite{LIU_2016,
Wulfinghoff_2018}, the reduction of computational complexity in nonlinear models
is achieved through the clustering of subdomains at the level of periodic unit
cells.

Although the above mentioned \emph{model order reduction (MOR)} approaches 
to the multiscale numerical modeling for nonlinear
problems seem to provide significant improvements relative to the direct and
very costly FE$^2$ method of the computational homogenization, the following
challenges can be seen:
\begin{list}{}{}
    \item Ch1: To minimize the off-line computations, as used in the POD
        approach: the precomputed solutions of the microscopic problem for various
        \emph{macroscopic deformation states (MDS)} are typically obtained for a range of
        macroscopic deformations containing a lot of MDS which are not used, \ie do
        not appear in the specific macroscopic problem.
    \item Ch2: To enable an efficient solving of a two-scale design optimization
        problem leading to functionally graded materials: the microstructure
        response depends on the MDS and also on the local design, \ie the
        microscopic configuration requires another extra parameterization by the
        design variables. In this context, the POD off-line stage would be even
        more demanding.
    \item Ch3: The clustering approach used, \eg in \cite{Benaimeche_2022,
        Chaouch_2024}, needs to establish an approximation scheme within each
        cluster of similar MDS; for this, typically AI-based (machine learning)
        methods are employed to solve the micro-problems. A common problem with
        this clustering approach arises from the inconsistency between the
        continuous evolution of the deformation, whereas the macroscopic
        properties (typically the tangent elasticity and the stress) are
        changing abruptly when the particular macroscopic point switches from
        one cluster to another between two consecutive iterations. Heuristic
        algorithms (like cluster freezing, \cf \cite{Chaouch_2024}) have been
        proposed in an effort to remedy this effect, which otherwise may stop
        the convergence of the macroscopic iterations. One may ask for
        alternative approaches to construct an approximation scheme over the set
        of all MDS represented by a limited number of clusters, each associated
        with a single microstructure and its microscopic deformation state.
\end{list}

\paragraph{Paper contribution}
In response to these challenges, we propose a new MOR-based method called
\emph{``Clustering and Sensitivity-based Approximation'' (CSA)} intended for the
two-scale computational homogenization. This method has the following features,
which also characterize the contributions of the present work:
\begin{itemize}
    \item The macroscopic properties of the homogenized medium are recovered in
        the entire domain using an approximation scheme based on the clustering
        approach, providing a continuous transition of the MDS from one centroid to
        another. This feature remedies the convergence issues observed in related
        works, see above Challenge 3.
    \item Within each subdomain of the macroscopic body associated with a
        cluster, the homogenized properties are given by the representative
        microstructure of the cluster centroid associated with a MDS.
To establish the homogenized properties 
at a macroscopic position belonging to the
        cluster, a linear mapping is constructed in the centroid using the sensitivity
        analysis of the macroscopic tangent elastic modulus and the effective stress \wrt the MDS variations. For this, the domain
        method of differentiation known from the shape optimization is adapted. The
        so-called ``design velocity field'' is presented by the microstructure
        displacement correctors (in the language of the asymptotic homogenization),
        the characteristic response of the microstructure \wrt the macroscopic deformation mode.
         The linear mapping is applied to relative MDS variations of stretches evaluated using the polar decomposition of the deformation gradient. Note that, in \cite{Chaouch_2024}, the ``linear approximation within the centroid'' is given for effective stresses only using the tangent elastic modulus which, however, is constant within the cluster.
    \item The CSA method does not rely on any ``off-line'' stage of the
        numerical modeling (see Challenge 1 above): the representative
        microstructures are introduced ``on-flight'' as the macroscopic
        deformation progresses. This approach significantly reduces the number
        of microscopic problems that need to be solved compared to the POD-based
        approach.
    \item Although the proposed CSA method is developed and tested for
        hyperelastic materials, the approximation scheme can be generalized
        within the same framework to include: a) functionally graded materials,
        \ie to account for quasi-periodic microstructures evolving with the position
        variation through the macroscopic body -- the issue pertinent for
        two-scale design optimization, see above Challenge 2; b) dependence
        of the microscopic configuration on internal variables, such as
        inelastic deformation, \cf \cite{Benaimeche_2022, Chaouch_2024}. It is
        worth noting that the extended parametrization by such additional
        variables (design, or internal variables) is much less expensive in the
        sense of computational complexity for the ``on-flight'' approach used in
        the CSA method than the analogous extension for the ``off-line''
        approach of pre-computed microstructure responses, such as the POD
        method, requiring to get responses for the Cartesian product of
        intervals for the considered micro-configuration parameters. Therefore,
        the CSA method can be of particular interest for two-scale structural
        optimization based on homogenization.
\end{itemize}

To demonstrate the computational efficiency and illustrate some of the above
claimed features, the proposed CSA method is compared with the full \dtss
simulation and the POD method. The algorithms of both the MOR methods are
implemented within the SfePy framework \cite{Cimrman_Lukes_Rohan_2019} and the
source code is available in the public repository \cite{Lukes_code_2025}.

\paragraph{Paper orgranization}
The paper is organized as follows. Section \ref{sec:hyperelasticity} introduces
the formulation of the large deformation analysis in heterogeneous hyperelastic
structures. Section \ref{sec:homogen} outlines the two-scale homogenization
framework, which yields the local subproblems described in Section
\ref{sec:homogen-micro} and the global problem with homogenized coefficients in
Section \ref{sec:homogen-macro}. The approximation strategy for the homogenized
coefficients is developed in Section \ref{sec:mor-cluster}. The complete CSA
algorithm is presented in Section \ref{sec:mor-cluster-algorithm}. Section
\ref{sec:mor-pod} provides a brief overview of the POD method. Finally, Section
\ref{sec:num_simulations} presents several numerical experiments designed to
assess the approximation error and demonstrate the efficiency of the CSA in
comparison with the \dtss{} and POD methods.

\paragraph{Notation}
Throughout this paper, the following notation is used. Since the problem
involves two spatial scales, we distinguish between the macroscopic coordinates
$x$ and the microscopic coordinates $y$. By $\pd_i = \pd_i^x$ we abbreviate the
partial derivative $\pd/\pd x_i$. We use $\nabla_x = (\pd_i^x)$ and $\nabla_y =
(\pd_i^y)$ when differentiation \wrt coordinate $x$ and $y$ is used,
respectively. By $\veps =  \ell/L$ we denote the scale parameter, relating the
microscopic length $\ell$ to the macroscopic one $L$.
Boldface symbols are used to represent vectors $\ab = (a_i)$ and second-order
tensors $\bb = (b_{ij})$. The fourth-order elasticity tensor is denoted by $\Aop
= (A_{ijkl})$ or $\Dop = (D_{ijkl})$. The homogenized coefficients are written
as $\Ahom = (\mathcal{A}_{ijkl})$ and $\Shom = (\mathcal{S}_{ij})$. The
second-order identity tensor is represented by $\Ib = (\delta_{ij})$ and
$\eb(\ub) = 1/2(\nabla \ub + \nabla \ub^T)$ is the linear strain tensor of a
vector field $\ub$. The dot symbol `$\cdot$' denotes the scalar product of two
vectors, while the colon `$:$' represents the inner product of two second-order
tensors. The Sobolev spaces $H^1(Y)$ of the square-integrable functions up to
the 1st order generalized derivative, defined in domain $Y$. For vector-valued
functions the bold notation is applied $\Hdb$. Subspaces of $Y$-periodic
functions are employed for vector-valued functions are denoted by $\HpbYi{Y}$.

%% file: part_main.tex
\section{Hyperelastic structures -- problem formulation}\label{sec:hyperelasticity}

We consider a compressible hyperelastic material, occupying an open bounded
domain $\Om\sveps \in \mathbb{R}$, whose mechanical behavior is governed by the neo-Hookean
material model. The Cauchy stress tensor is expressed in terms of the bulk
modulus $K\sveps$ and the shear modulus $\mu\sveps$ as
\begin{equation}\label{eq:neohook}
	\sigmabf\sveps = K\sveps(J - 1)\Ib + \mu\sveps J^{-5/3}\, (\bb - \mathrm{tr}(\bb)/3 \Ib),
\end{equation}
where $\Fb$ is the deformation gradient, $J=\det(\Fb)$ is the relative volume
change, $\bb = \Fb\Fb^T$ is the left Cauchy-Green deformation tensor, $\Ib =
(\delta_{ij})$ is the Kronecker symbol, and $\mathrm{tr}(\bb) = b_{ii}$ is the
trace of $\bb$. The superscript ${}\sveps$ indicates spatially fluctuating
quantities depending on the scale parameter $\veps$, otherwise employed in the asymptotic homogenization analysys.
We demonstrate our computational method using the simple two-parametric
neo-Hookean model. However, the derivations presented in the following sections
are general enough to employ more complex material models such as Mooney-Rivlin,
Ogden, or Yeoh.
In the Eulerian frame, the constitutive equation (\ref{eq:neohook}), together
with the equilibrium equation imposed in the deformed configuration,
\begin{equation}\label{eq:equilibrium}
	- \nabla \cdot \sigmabf\sveps = \fb\sveps,\quad \mbox{ in } \Ome\;,
\end{equation}
and boundary conditions
\begin{equation}\label{eq:bc}
    \begin{aligned}
        \ub\sveps = \bar\ub\sveps\ \mbox{ on } \partial_u \Ome,
        \quad \nb \cdot \sigmabf\sveps = \bar\hb\sveps\ \mbox{ on } \partial_\sigma \Ome,\\
        \partial_u \Ome \cup \partial_\sigma \Ome = \partial \Ome,
        \quad \partial_u \Ome \cap \partial_\sigma \Ome = \emptyset,
    \end{aligned}
\end{equation}
represent the non-linear problem, which is further linearized 
in the framework of the updated Lagrangian formulation. The linearized
incremental form of the problem will be subjected to the homogenization
procedure, resulting in a coupled system of macroscopic and microscopic
equations. External volume forces $\fb\sveps$, surface tractions
$\bar\hb\sveps$, whereby $\nb$ is the unit normal, and the prescribed
displacements $\bar\ub\sveps$ are applied in a series of quasi-static loading
steps.

\subsection{Weak formulation}\label{sec:weak}

Following the procedure described in \cite{Lukes_2022}, we define the sets of admissible
displacements
\begin{equation}\label{eq:asets}
    \Vb\sveps \in \lbrace\vb\sveps \vert \vb \in \Hdb(\Ome)\;\sveps = \bar\ub\sveps\mbox{ on } \partial_u\Ome \rbrace,
    \qquad
    \Vb\sveps_0 \in \lbrace\vb\sveps \vert \vb \in \Hdb(\Ome)\;\sveps = 0\mbox{ on } \partial_u\Ome \rbrace,
\end{equation}
respecting the boundary conditions \eq{eq:bc}. In the above definitions,
sufficient regularity for unknown and test displacements is assumed. By
$\Vb\sveps(t)$, $\Vb\sveps_0(t)$ we express the dependency on the spatial configuration
$\Ome(t)$. The weak form of the non-linear problem \eq{eq:neohook}-\eq{eq:bc} for
the unknown displacement field $\ub\sveps \in \Vb\sveps$ is given by
\begin{equation}\label{eq:weak_equilibrium}
    \Phi_t(\ub\sveps, \vb\sveps) = \intOmet \sigmabf(\ub\sveps) : \nabla \vb\sveps \dOm
    - \intOmet \fb\sveps \cdot \vb\sveps \dOm
    - \intOmebt \bar\hb\sveps \cdot \vb\sveps \dS = 0,
    \quad \forall \vb\sveps \in \Vb\sveps_0(t),
\end{equation}
where $\Phi_t(\ub\sveps, \vb\sveps)$ is the residual function evaluated at time $t$. In
order to get the incremental formulation of the above problem a constant time
step $\dt$ is assumed and the residual equation is rewritten at time $t + \dt$
employing the first order Taylor expansion at time $t$,
\begin{equation}\label{eq:weak_residual}
    \Phi_{t+\dt}(\ub\sveps(t + \dt), \vb\sveps) \approx
    \Phi_{t}(\ub\sveps(t), \vb\sveps) + \delta \Phi_{t}(\ub\sveps(t), \vb\sveps)\circ(\dub\sveps, \dt\dvel).
\end{equation}
Time increment $\delta \Phi_{t}(\ub\sveps(t), \vb\sveps)\circ(\dub\sveps,
\dt\dvel)$ appearing on the right-hand side of \eq{eq:weak_residual} is
expressed in terms of displacement increment $\dub\sveps$, defined such that
$\ub\sveps(t + \dt) \approx \ub\sveps(t) + \dub\sveps = \ub\sveps(t) +
\dot\ub\sveps \dt$, and convection field $\dvel$, which is also used to
approximate of the perturbed configuration of the domain: $\Ome(t + \dt) \approx
\Ome(t) + \dt\dvel$, see\cite{Rohan_2017, Lukes_2022} for details.
The differential of the residual function $\delta \Phi_{t}$ is further approximated
by its total time derivative, so that
%
\begin{equation}\label{eq:weak_residual_derivative}
    \begin{aligned}
        \delta \Phi_{t} \approx \dt \dot\Phi_{t} = 
        \dt
        \bigg(
            &\intOmet \left(\nabla\dvel \sigmabf\sveps + \Lie{\dvel} \sigmabf\sveps \right):\nabla\vb\sveps \dOm\\
        - &\intOmet \left(\dot\fb\sveps \cdot \vb + \fb\sveps \cdot \vb\sveps \nabla \cdot \dvel \right)\dOm
        - \intOmebt \left(\dot{\bar\hb}\sveps \cdot \vb\sveps + \bar\hb\sveps \cdot \vb\sveps \nabla \cdot \dvel \right)\dS
        \bigg)
        ,
    \end{aligned}
\end{equation}
where $\Lie{\dvel} \sigmabf\sveps$ denotes the Lie derivative of the Cauchy
stress tensor with respect to the convection field $\dvel$. This term can be
expressed using the tangential stiffness operator $\Dop\sveps$ and the linear
strain rate $\eb(\dot\ub\sveps)$ as
\begin{equation}\label{eq:weak_lie}
    \Lie{\dvel}\sigmabf\sveps = \Dop\sveps \eb(\dot\ub\sveps).
\end{equation}
The tangential fourth-order stiffness tensor $\Dop\sveps$ is the Truesdell rate of the Cauchy
stress. For the compressible neo-Hookean material described by the stress-strain relation
\eq{eq:neohook}, it takes the form (see \cite{Crisfield_1991} and Appendix~\ref{sec:appendix})
\begin{equation}\label{eq:dop}
    \begin{aligned}
    D\sveps_{ijkl} =
      &K\sveps(2J - 1)\delta_{ij}\delta_{kl}
       -K\sveps(J-1)\left(\delta_{ik}\delta_{lj} + \delta_{il}\delta_{jk}\right)\\
    + &\mu\sveps\; J^{-5/3} \left(
      \frac{2}{9} \mathrm{tr}(\bb) \delta_{ij}\delta_{kl}
    - \frac{2}{3}\left(b_{ij}\delta_{kl} + \delta_{ij}b_{kl}\right)
    + \frac{1}{3}\mathrm{tr}(\bb)\left(\delta_{ik}\delta_{lj} + \delta_{il}\delta_{jk}\right)
    \right).
    \end{aligned}
\end{equation}

The prescribed displacements $\bar\ub\sveps$, volume forces $\fb\sveps$, and
surface forces $\bar\hb\sveps$ are applied in a finite number of times steps
$t_k = k\,\dt$, $t = 0, 1, 2, \dots$, $N_t$ assuming a constant time step $\dt$. Using
this time discretization, the time derivative $\dot\ub\sveps$ appearing in
\eq{eq:weak_lie} can be replaced by the forward finite difference,
\begin{equation}\label{eq:du_fd}
    \dot\ub\sveps(t_k) \approx (\ub\tkpe - \ub\tke) / \dt = \dub\tkpe / \dt,
\end{equation}
and the velocity field $\dvel(t_k)$ can be approximated by the unknown
displacement increment $\dub\tkpe$,
\begin{equation}\label{eq:dvel_du}
    \dvel(t_k) \approx \dub\tkpe / \dt.
\end{equation}
Substituting \eq{eq:weak_residual_derivative}-\eq{eq:dvel_du} into equation
\eq{eq:weak_residual} and with respect to the equilibrium
\eq{eq:weak_equilibrium} we get the following incremental formulation of the nonlinear
problem: For a given Cauchy stress $\sigmabf\tke$ and new loads $\fb\tkpe$, $\bar\hb\tkpe$
find displacement increments $\dub\tkpe \in \delta\Vb\sveps$ ($\delta\Vb\sveps$ is a set of admissible
displacement increments) such that
\begin{equation}\label{eq:weak_incremental}
    \begin{aligned}
    \intOmek & \Aop\sveps \nabla \dub\tkpe : \nabla {\vb\sveps} \dOm
    = - \intOmek \sigmabf\tk : \nabla \vb\sveps \dOm\\
      & + \intOmek \fb\tkp \cdot \vb\sveps \dOm
      + \intOmebk \bar\hb\tkp \cdot \vb\sveps \dS,
      \qquad \forall \vb\sveps \in \Vb\sveps_0(t).
    \end{aligned}
\end{equation}
By $\Om\tke$ we refer to the actual spatial configuration $\Ome(t_k)$. 
The symmetric fourth-order tangent tensor $\Aop\sveps$ is defined as
\begin{equation}\label{eq:weak_elastic_operator}
    A_{ijkl}\sveps = D_{ijkl}\sveps + \sigma\sveps_{jl} \delta_{ik}.
\end{equation}
Note that the incremental problem \eq{eq:weak_incremental} still remains
nonlinear due to the dependence of $\Aop^\veps$ on the current deformation state,
and therefore must be solved in an iterative loop, e.g. by a Newton's method.

\section{Two-scale model of locally periodic structures}\label{sec:homogen}

The standard homogenization procedure is applied to the nonlinear incremental
problem~\eq{eq:weak_incremental} to derive the limit model for
$\veps\rightarrow0$, where the non-dimensional parameter $\veps$ defines the
period size of the microstructure. The justification for applying two-scale
homogenization to materials undergoing large deformations, and thus potentially
losing periodicity, is discussed in detail in \cite{Lukes_2022}, Section 4,
where the concepts of slowly varying microstructures and locally periodic media
are introduced.

The heterogeneous material with a periodic structure can be characterized by a
reference periodic cell $Y$, which may consist of several parts with different
material properties. In the numerical examples presented in
Section~\ref{sec:num_simulations}, we consider a two-phase microstructure
composed of a compliant matrix with embedded stiffer inclusions, as illustrated
in the left part of Figure~\ref{fig:num_micro_macro}. 

The displacement increment in \eq{eq:weak_incremental}
can be expressed in the form of the truncated asymptotic expansion
\begin{equation}\label{eq:expansion}
    \Tuf{\dub\tkpe} = \dub^0(x) + \veps \dub^1(x, y) + O(\veps^2), \qquad x\in\Omega,\quad y\in Y,
\end{equation}
where $\Tuf{}$ is the unfolding operator, see e.g. \cite{Cioranescu_2018},
$\dub^0(x)$ is the macroscopic displacement increment depending on $x$, and $\dub^1(x, y)$
is a two-scale function. The analogous expansion is used for the test function $\vb^\veps$.
For the unfolded elastic tensor $\Aop^\veps$ and the terms at the
r.h.s.\ of \eq{eq:weak_incremental} the following expansions hold:
\begin{equation}\label{eq:homogen2}
    \begin{aligned}
        \Tuf{\Aop\veps} & = \tilde \Aop(x,y) + \veps  \Aop^*(x,y),\\
        \Tuf{\sigmabf\tke} & = \tilde\sigmabf(x,y) + \veps \sigmabf^*(x,y) + O(\veps^2),\\
        \Tuf{\fb\tkpe} & = \tilde\fb(x,y)  + O(\veps),\\    
        \Tuf{\bar\hb\tkpe} & = \tilde\hb(x,y)  + O(\veps).\\
    \end{aligned} 
\end{equation}

Substituting the expansion \eq{eq:expansion} into
\eq{eq:weak_incremental}, taking into account $\nabla_y \dub^0(x) = \boldsymbol{0}$,
and passing to the limit $\veps \rightarrow 0$ yields
\begin{equation}\label{eq:homogen1}
    \begin{aligned}
        \intOm \intY & \tilde\Aop\left(\nabla_x\dub^0 + \nabla_y\dub^1\right)
    :\left(\nabla_x\vb^0 + \nabla_y\vb^1\right) \dVy\dOm
   = \\ 
   &\intOm\intY \tilde\sigmabf :\left(\nabla_x\vb^0 + \nabla_y\vb^1\right) \dVy\dOm
    + \intOm\intY \tilde\fb \cdot \vb^0 \dVy\dOm + \intOmb\intY \tilde\hb \cdot \vb^0 \dVy\dS.
    \end{aligned}
\end{equation}
Terms with tensors $\Aop^*$, $\sigmabf^*$ vanish in the limit, see
Eq.~(58) and Section~4.4 in \cite{Lukes_2022}. The integral $\intY$
denotes the average value over the microscopic domain $Y$
\begin{equation}\label{eq:avgint}
    \intY f \dVy = \frac{1}{\vert Y \vert} \int_Y f \dVy.
\end{equation}

Recall that \eq{eq:homogen1} is obtained from the weak form
\eq{eq:weak_incremental} of the linearized sub-problem associated with one load
increment of the nonlinear problem arising from equilibrium \eq{eq:equilibrium}.
For each load increment, incorporating volume forces and the boundary tractions,
the two-scale nonlinear problem is solved iteratively using corrections
$(\dlt\ub^0, \dlt\ub^1)$ satisfying \eq{eq:homogen1} with the right-hand side
presenting the actual out-of-balance. Within one iteration called ``outer''
labelled by $i$, the couple $(\dlt\ub^0, \dlt\ub^1)$ is computed in the
following steps which are described in detail below, in
Sections~\ref{sec:homogen-micro}--\ref{sec:update-micro}.
\begin{enumerate}
    \item \textbf{Step~M0}: computing $\dlt\ub^0$ and updating the macro-configuration $\Om$, see Section~\ref{sec:homogen-macro};
    \item \textbf{Step~m1}: computing $\dlt\ub^1$ and updating the micro-configurations $Y(x)$, for $x \in \Om$, see Section~\ref{sec:update-micro};
    \item \textbf{Step~m2}: updating the effective tangential (incremental)
    elasticity $\Ahom(x)$ and effective Cauchy stress $\Shom(x)$ defined in~\eq{eq:homog_coefs}, see Section~\ref{sec:homogen-micro}.
\end{enumerate}

\begin{figure}
    \centering
    \includegraphics[width=0.85\figscale]{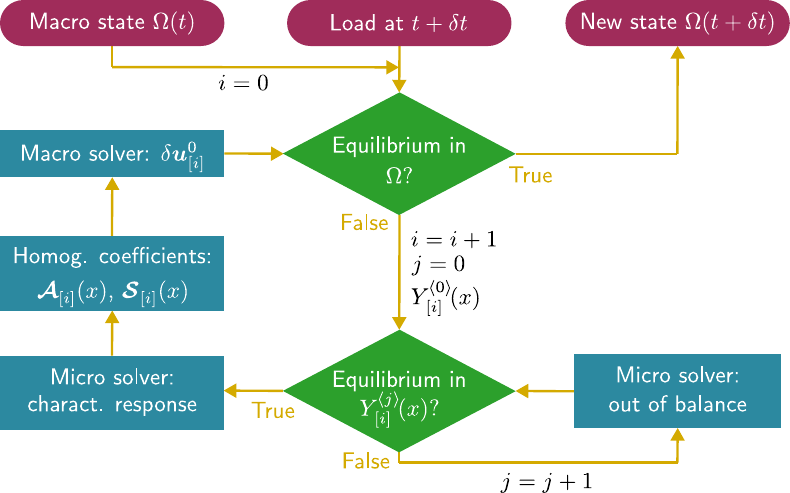}
    \caption{Inner/outer iterations within one ``loading-time'' step.}
        \label{fig:ulf-iterations}
\end{figure}

Thus, the problem is decomposed into macroscopic subproblems for corrections
$\dlt\ub^0\itMa{i}$ in response to the out-of balance (iterations labelled by
$i$), and the local ``microproblems'' which split into 2 steps: \textbf{Step~m1}
involves ``inner'' iterations $j = 0,1,\dots$ associated with computing local
equilibrium for given macroscopic strain -- for a given global iteration $i$, a
sequence of local micro-configurations $Y\itMa{i}\itmi{j}(x)$ is constructed,
starting with $\tilde Y\itMa{i}(x) = Y\itMa{i}\itmi{0}(x)$ which is given by
virtue of the characteristic response computed in $Y\itMa{i-1}(x)$ for the
preceding macroscopic iteration $i-1$; \textbf{Step~m2} consists in solving
characteristic response to applied macroscopic strain increment modes. These
steps are arranged in the computational scheme which is illustrated in
Fig.~\ref{fig:ulf-iterations}.

\begin{remark}\label{rem-mic-mac}
    From \eq{eq:homogen1}, the local problems are obtained for vanishing $\vb^0$,
    and testing by $\vb^1:=\vb \in \HpbY$, hence at almost all $x\in \Om$, 
    \begin{equation}\label{eq:homogen1a}
        \begin{aligned}
        \intY & \tilde\Aop\left(\nabla_x\dub^0 + \nabla_y\dub^1\right)
        :\nabla_y\vb \dVy
    = \intY \tilde\sigmabf : \nabla_y\vb\dVy, \qquad \forall \vb \in \HpbY\;.
        \end{aligned}
    \end{equation}
    Assume the previous iteration $i-1$ yielding $\dlt\ub^0\itMa{i-1}$ due to
    \textbf{Step~M0}, which means also updating the local configuration
    $Y\itMa{i-1}(x) \mapsto \tilde Y\itMa{i}(x)$. Due to \textbf{Step~m1}, the
    local equilibrium in \eq{eq:homogen1a} is achieved: let $\dub^1 =
    \dlt\tilde\wb^\sigma + \dlt\tilde\ub^1$, such that \eq{eq:homogen1a} can be
    written (at any $x \in \Om$) as
    \begin{equation}\label{eq:homogen1b}
        \begin{aligned}
        \intYi{\tilde Y\itMa{i}} & \tilde\Aop\left(\nabla_x\dub^0 + \nabla_y(\dlt\tilde\ub^1 + \dlt\tilde\wb^\sigma)\right)
        :\nabla_y\vb \dVy
    = \intYi{\tilde Y\itMa{i}}\tilde\sigmabf : \nabla_y\vb\dVy, \qquad \forall \vb \in \HpbYi{\tilde Y\itMa{i}}\;.
        \end{aligned}
    \end{equation}
    This problem is decomposed into two subproblems: Find
    $\dlt\tilde\wb^\sigma\in \HpbYi{\tilde Y\itMa{i}}$, such that (we drop the
    macroscopic iteration index $[i]$)
    \begin{equation}\label{eq:homogen1c}
        \begin{aligned}
            \intYi{\tilde Y} & \tilde\Aop\nabla_y \dlt\tilde\wb^\sigma:\nabla_y\vb \dVy
            = \intYi{\tilde Y}\tilde\sigmabf : \nabla_y\vb\dVy, \qquad \forall \vb \in \HpbYi{\tilde Y}\;.
        \end{aligned}
    \end{equation}
    Using the solution $\dlt\tilde\wb^\sigma$, the micro-configuration is
    updated: $\tilde Y := \tilde Y + \dlt\tilde\wb^\sigma$ and the stress
    $\tilde\sigmabf$ is computed using the actual deformation gradients (recall
    we consider a hyperelastic material). Problem \eq{eq:homogen1c} formulated
    in the updated configuration is solved repeatedly in iterations
    $j=0,1,\dots$ until the \rhs evaluated with the actual $\tilde\sigmabf$
    vanishes, \ie until $\nabla_y \tilde\sigmabf = 0$ in $\tilde Y$. In the
    schematic view of the whole algorithm, Fig.~\ref{fig:ulf-iterations}, the
    provisional micro-configurations $\tilde Y = Y\itMa{i}\itmi{j}$, are
    labelled by index $j$.

    Then \eq{eq:homogen1b} considered in the actual local equilibrium
    configuration $Y\itMa{i}:=\tilde Y$ yields the linear problem: Given
    $\dub^0$, find $\dlt\tilde\ub^1 \in \HpbYi{Y\itMa{i}}$, such that
    \begin{equation}\label{eq:homogen1d}
        \begin{aligned}
            \intYi{Y\itMa{i}} & \tilde\Aop\left(\nabla_x\dub^0 + \nabla_y\dlt\tilde\ub^1\right):\nabla_y\vb \dVy
            = 0, \qquad \forall \vb \in \HpbYi{Y\itMa{i}}\;.
        \end{aligned}
    \end{equation}
    This identity can be interpreted in the context of the two-scale limit
    \eq{eq:homogen1}, whereby $\dlt\tilde\ub^1 = \dub^1$, since the local
    configuration is in the local equilibrium, $\dlt\tilde\wb^\sigma = 0$, and
    the \rhs is presented by the macroscopic external forces only. By virtue of
    the linearity, the characteristic response can be introduced, which enables
    obtaining the homogenized elasticity coefficients. Consequently, the
    two-scale computational procedure proceeds with \textbf{Step~M0}. Details
    are given below.
\end{remark}

In subsequent sections, the three main steps embedded in the cycles associated
with the ``outer'' and ``inner'' iterations are presented in a shifter order:
\textbf{Step~m2}, \textbf{Step~M0} and   \textbf{Step~m1}.

\subsection{Local microscopic problem -- Step~m2}\label{sec:homogen-micro}

By virtue of the linearity of the limit equations the fluctuating function
$\dub^1$ can be decomposed in terms of the so called characteristic responses
$\omegabf^{ij}$ and the gradient of the macroscopic deformation $\nabla_x\dub^0$,
\begin{equation}\label{eq:decomposition_responses}
    \dub^1 = \omegabf^{ij} \left(\nabla_x \du^0_i\right)_j.
\end{equation}
Taking $\vb^0=\boldsymbol{0}$ in \eq{eq:homogen1} and using decomposition
\eq{eq:decomposition_responses}, the local problem for the unknown
characteristic response can be defined: Find $\omegabf^{ij} \in \HpbY$, $i,j = 1, 2, (3)$,
satisfying
\begin{equation}\label{eq:micro}
    \alin{\omegabf^{ij} + \Pibf^{ij}}{\vb^1} = 0,\qquad \forall \vb^1 \in \HpbY,
\end{equation}
where $\Pi^{ij}_k = y_j \delta_{ik}$ and $\alin{\ub}{\vb}$ is the bilinear form
\begin{equation}\label{eq:alin}
    \alin{\ub}{\vb} = \intY \tilde\Aop \nabla_y \ub : \nabla_y \vb \dVy.
\end{equation}
The functional space $\HpbY$ denotes the standard Sobolev space of $Y$-periodic
functions. Note that the integral in \eq{eq:alin} is evaluated over the deformed
microscopic configuration $Y = Y(x,t)$.

\subsection{Global macroscopic problem and homogenized coefficients -- Step~M0}\label{sec:homogen-macro}

For the vanishing microscopic test functions $\vb^1 = \boldsymbol{0}$
equation \eq{eq:homogen1} results in
%
\begin{equation}\label{eq:aux}
    \begin{aligned}
        \intOm \intY & \tilde\Aop\left(\nabla_x\dub^0 + \nabla_y\dub^1\right)
        :\nabla_x\vb^0 \dVy \dOm
        = \\ 
        &\intOm\intY \tilde\sigmabf:\nabla_x\vb^0 \dVy \dOm
        + \intOm\intY \tilde\fb \cdot \vb^0 \dVy \dOm + \intOmb\intY \tilde\hb \cdot \vb^0 \dVy \dS,
    \end{aligned}
\end{equation}
%
%
Substituting the decomposition \eq{eq:decomposition_responses} into \eq{eq:aux},
the macroscopic equilibrium problem can be formulated as follows: For given
volume and surface forces included in $L^{\rm new}(\vb^0)$ and homogenized
coefficients $\Ahom$, $\Shom$, find the macroscopic displacement increments
$\dub^0 \in \delta\Vb(\Om)$ such that
\begin{equation}\label{eq:macro}
    \intOm \Ahom \nabla_x\dub^0:\nabla_x\vb^0 \dOm
    = L^{\rm new}(\vb^0) - \intOm \Shom:\nabla_x\vb^0 \dOm, \qquad \forall \vb^0 \in \Vbzero(\Om).
\end{equation}
The homogenized elastic tensor $\Ahom$ and averaged Cauchy stress tensor
$\Shom$ can be identified from \eq{eq:aux} as:
\begin{equation}\label{eq:homog_coefs}
    \begin{aligned}
        \Acal_{ijkl} &= \alin{\omegabf^{kl} + \Pibf^{kl}}{\omegabf^{ij}} = \alin{\omegabf^{kl} + \Pibf^{kl}}{\omegabf^{ij} + \Pibf^{ij}},\\
        \Shom &= \intY \tilde\sigmabf \dVy,
    \end{aligned}
\end{equation}
see \cite{Rohan_2003, Lukes_2022} for details.


\subsection{Update microscopic configurations -- Step~m1}\label{sec:update-micro}
\begin{figure}
    \centering
    \includegraphics[width=0.99\figscale]{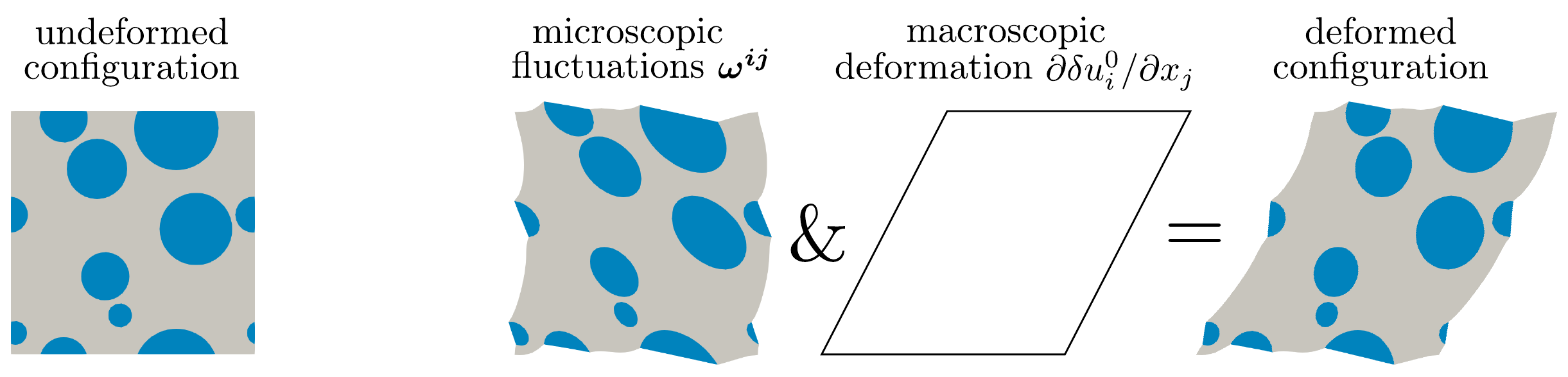}
    \caption{Updating microscopic configurations using the characteristic responses
             and $\omegabf^{ij}$ and the gradient of the macroscopic
             displacement increments $\dub^0$.}
        \label{fig:micro-update}
\end{figure}

The solution of the microscopic subproblems and the evaluation of the
coefficients at a given macroscopic point $x \in \Omega$ must be performed on a
update microscopic configurations $Y(x)$ following the evolution of the
macroscopic deformation $\dub^0$. The updating procedure is illustrated in
Fig.~\ref{fig:micro-update} and can be expressed as
\begin{equation}\label{eq:micro-update}
    \begin{aligned}
        \delta \wb &= \left(\omegabf^{ij} + \Pibf^{ij}\right) (\nabla_x \delta u^0_i)_j\\
        \tilde Y(x) &:= Y(x) + \{\delta \wb\}.
    \end{aligned}
\end{equation}
The stress $\tilde\sigmabf$ updated due to \eq{eq:micro-update} does not automatically satisfy the microscopic equilibrium, so
\begin{equation}\label{eq:micro-equilibrium}
    \nabla_y \cdot \tilde\sigmabf \not = \boldsymbol{0}\;.
\end{equation}
To ensure the local equilibrium (\ie for all local configurations $\tilde Y(x)$, $x \in \Omega$), as explained in Remark~\ref{rem-mic-mac}, see problem \eq{eq:homogen1c},
 corrections
$\delta\wb \in \HpbYi{Y\itmi{j}}$ are computed in iterations $j = 0,1,\dots$ by solving 
\begin{equation}\label{eq:micro-equilibrium2}
    \alinY{\delta\wb^{j+1}}{\vb}{Y\itmi{j}} = - \intYi{Y\itmi{j}} \tilde\sigmabf : \nabla_y \vb \dVy,\qquad \forall \vb \in \HpbYi{Y\itmi{j}},
\end{equation}
where $Y\itmi{j+1} = Y\itmi{j} + \{\delta \wb^{j+1}\}$ and $Y\itmi{0} =\tilde Y$ according to \eq{eq:micro-update}$_2$.

\subsection{Sensitivities of homogenized coefficients}\label{sec:homogen-sa}

The approximation scheme for efficient computing of the deformation dependent
homogenized coefficients proposed in Section~\ref{sec:mor-cluster} employs the
coefficient sensitivities $\dlt\Ahom$ and $\dlt\Shom$ that are introduced
in this section.

Based on the shape sensitivity theory and on the domain parametrization technique
(see \cite{Haslinger_1996} and Appendix~\ref{sec:appendix}),
the following sensitivity of the bilinear form $\alin{\ub}{\vb}$ holds
\begin{equation}\label{eq:sa_alin}
    \begin{aligned}
        \dltsh \alin{\ub}{\vb} = \intY&  \tilde{A}_{irks}\left(
            \delta_{rj}\delta_{sl} \nabla_y \cdot \Vcalbf - \delta_{jr}\pd_s^y \Vcal_l - \delta_{ls}\pd_r^y \Vcal_j
        \right) \pd^y_l u_k \pd^y_j v_i \dVy\\
            & + \intY \dlt\tilde\Aop \nabla_y \ub : \nabla_y \vb \dVy
            - \alin{\ub}{\vb} \intY \nabla_y\cdot \Vcalbf\dVy,
    \end{aligned}
\end{equation}
cf.{} \cite{Rohan_2015} or \cite{Rohan_2024}, where the similar formulae for a
piezoelectric medium are derived. With \eq{eq:sa_alin} in hand, employing the abbreviation
$\Xibf^{ij} = \omegabf^{ij}+ \Pibf^{ij}$, and introducing $\dltsh\Pi^{ij}_k = \Vcal_j \delta_{ik}$,
the sensitivity expression for the elastic tensor can be defined,
\begin{equation}\label{eq:sa_elastic}
        \dlt \Acal_{klij} = \dltsh \alin{\Xibf^{ij}}{\Xibf^{kl}}
            + \alin{\dltsh\Pibf^{ij}}{\Xibf^{kl}}
            + \alin{\Xibf^{kl}}{\dltsh\Pibf^{kl}}.
\end{equation}
For the sensitivity of the averaged Cauchy stress, the following expression holds:
\begin{equation}\label{eq:sa_stress}
    \dlt \Shom =  \intY ( \tilde\sigmabf-\Shom)\nabla_y\cdot \Vcalbf\dVy
    + \intY\dlt\tilde\sigmabf \dVy.
\end{equation}
Note that the first integral, recalling $\Shom$ is constant in $Y$, describes the sensitivity of the stress inhomogeneous on the volumetric change of the RVE (cell $Y$), whereas the second one is the sensitivity of the local stress in $Y$.
Details upon derivation of \eq{eq:sa_elastic} and \eq{eq:sa_stress} are given in
Appendix~\ref{sec:appendix}.

The sensitivity of the homogenized coefficients with respect to the components
of the macroscopic deformation $\left(\nabla_x \du^0_i\right)_j$ is obtained by
substituting $\Pibf^{ij} + \omegabf^{ij}$ for the velocity field $\dvel$ in
\eq{eq:sa_elastic} and \eq{eq:sa_stress}.
The derivatives $\dlt\tilde\sigmabf$, $\dlt\tilde\Aop$ are computed employing
relations~\eq{eq:app2-cs} and~\eq{eq:app2-dtk}.


\section{Approximation of homogenized coefficients}\label{sec:mor-cluster}

Solving the nonlinear macroscopic problem using the finite element method
requires computing the deformation-dependent tangent moduli and stress tensors
at all quadrature points of the discretized macroscopic domain. The tangent
moduli and stresses are associated with the solutions of the local microscopic
subproblems \eq{eq:micro}, and the computation must be repeated at every iteration of
the nonlinear algorithm. This approach, commonly referred to as FE$^2$, is highly
computationally demanding even for relatively simple 2D problems. To avoid the
computational bottleneck of the FE$^2$ method, we propose a novel algorithm, CSA (``Clustering and Sensitivity based
Approximation''), that
reduces the number of solved local subproblems. Our approach leverages the K-means
clustering method to group macroscopic deformations and employs sensitivity
analysis to approximate the homogenized coefficients efficiently.

Let $\{\hx^i\} \subset \Om$, $i = 1, \dots, \hat N$, be the set of all
macroscopic quadrature points where the homogenized coefficients $\Ahom(\hx)$
and $\Shom(\hx)$, which depend on the macroscopic deformation $\Fb(\hx)$, must
be evaluated. The first step of our algorithm involves performing the polar
decomposition of the deformation gradient $\Fb(\hx)$ into the rotation tensor
$\hat\Rb$ and the symmetric stretch tensor $\hat\Ub$,
\begin{equation}\label{eq:polar_decomposition}
    \Fb(\hat x) = \hat\Rb \hat\Ub.
\end{equation}
The deformation state of the microstructure is fully characterized by the
symmetric strain tensor $\eb(\hat x) = \hat\eb = \hat\Ub - \Ib$, as the rotation
tensor $\hat\Rb$ does not influence the local response. The components of
$\hat\eb$ can be interpreted as coordinates of a point in $\RR^3$ for 2D
problems or $\RR^6$ for 3D problems, and the proximity of two deformation states
can be quantified using the $L^2$ norm.

The key aspect of the proposed approach is to find the minimal number
$\acirc N$ of sets $O^j$, $j=1, ..., \acirc N$, such that their union
$\bigcup_{j=1}^{\acirc N} O^j$ covers all current macroscopic strains $\hat\eb$:
\begin{equation}\label{eq:sets_O}
    O^j = \{\hat\eb: \vert \hat\eb - \acirc\eb^j \vert < \rho\},\qquad
    \{ \hat\eb\} \subset \bigcup_{j=1}^{\acirc N} O^j,
\end{equation}
where $\acirc\eb^j$ denotes the mean value (center) of the strains in $O^j$,
$\vert \hat\eb \vert$ is the Euclidean norm in $\RR^3$ or in $\RR^6$
respectively, and $\rho$ is a given parameter. The sets $O^j$ may overlap,
allowing a single strain value to belong to multiple sets, which can improve the
accuracy of the coefficients  approximation. Among the various strategies for
finding the minimal number of sets $O^j$, we adopt the K-means clustering method
\cite{Lloyd_1982}, as suggested in \cite{Benaimeche_2022}. Throughout this
paper, the sets $O^j$ are referred to as centroids.

Instead of solving the local microscopic subproblems for all $\hat N$
macroscopic deformations $\eb(\hx^i)$, we solve them now only for a reduced
number of deformations $\acirc\eb^j$, $j=1,\dots,\acirc{N}$,
to obtain the corresponding homogenized coefficients $\acirc\Ahom^j$ and
$\acirc\Shom^j$. The number of centroids $\acirc N$ may increase dynamically
during the macroscopic iterations whenever newly occurring strains cannot be
classified into the existing centroids.

\begin{figure}
    \centering
    \includegraphics[width=0.65\figscale]{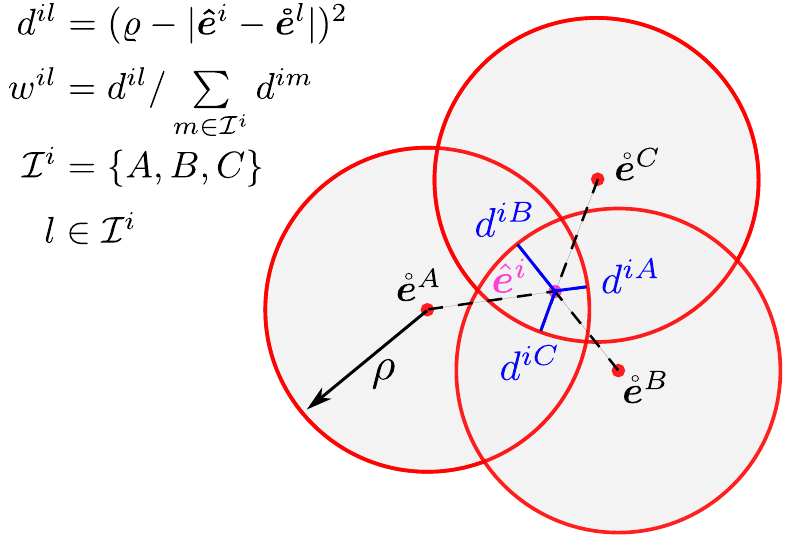}
    \caption{Calculation of weights $w^{il}$ for the case of three overlapping centroids.}
    \label{fig:centroids_weights}
\end{figure}

To approximate the coefficients $\Ahom(\hx)$ and $\Shom(\hx)$ at the quadrature
points $\hx$, we compute not only $\acirc\Ahom$ and $\acirc\Shom$, but also
their sensitivities $\delta \acirc\Ahom$, $\delta \acirc\Shom$ with respect to
deformations. This computation is efficient, as the sensitivities depend only on
the already computed correctors $\omegabf$, see Section~\ref{sec:homogen-sa} and
\cite{Rohan_2003}. Given all coefficients $\acirc\Ahom^j$, $\acirc\Shom^j$ and
their sensitivities $\delta \acirc\Ahom^j$, $\delta \acirc\Shom^j$, coefficients
$\Ahom(\hx^i)$, $\Shom(\hx^i)$ at any quadrature point $\hat{x}^i$  can be
approximated as
\begin{equation}\label{eq:cf_approx}
    \begin{aligned}
        \Acal_{klmn}(\hx^i) &\approx \hat{R}^i_{kp} \hat{R}^i_{lq} \hat{R}^i_{mr} \hat{R}^i_{ns}
            \tilde\Acal^i_{pqrs}, \qquad 
            \tilde\Ahom^i = \sum_{l\in\ihset{i}}  \left(\acirc\Ahom^l + \delta \acirc\Ahom^l \circ \gb^{il} \right) w^{il},\\
        \Scal_{kl}(\hx^i) &\approx \hat{R}^i_{kp} \hat{R}^i_{lq}
            \tilde\Scal^i_{pq}, \qquad 
            \tilde\Shom^i = \sum_{l\in\ihset{i}}  \left(\acirc\Shom^l + \delta \acirc\Shom^l \circ \gb^{il} \right) w^{il},
    \end{aligned}
\end{equation}
where the index set $\ihset{i}$ is defined by $j \in \ihset{i} \Leftrightarrow
\hx^i \in O^j$, so that the sum runs over all centroids containing $\hat\eb^i$.
%
The weights $w^{il}$ are introduced to ensure a smooth approximation when centroids overlap:
\begin{equation}\label{eq:weights}
    w^{il} = d^{il} /
        \sum_{m\in\ihset{i}} d^{im},\qquad d^{il} = (\rho - \vert \hat\eb^i - \acirc\eb^l\vert)^2, \qquad \sum_{l\in\ihset{i}} w^{il} = 1,
\end{equation}
The relative deformation $\gb^{il}$ appearing in \eq{eq:cf_approx} represents
the deformation of $\hat\eb^i$ with respect to the centroid $\acirc\eb^l$ and is
defined using the multiplicative decomposition,
yielding
\begin{equation}\label{eq:rel_g}
    \gb^{il} = (\hat\eb^i + \Ib)(\acirc\eb^l + \Ib)^{-1} - \Ib.
\end{equation}
Figure~\ref{fig:centroids_weights} illustrates the weight calculation for the
case of three overlapping centroids.

The deformed microstructures $Y(\hat\eb^i)$ associated with the macroscopic
deformation states $\hat\eb^i$ are obtained by
\begin{equation}\label{eq:Y_approx}
    \hat\yb^i \approx \left(
        \sum_{l\in\ihset{i}}  \left(\acirc\yb^l + \acirc\Xibf^{rs,l} g^{il}_{rs} \right) w^{il}
    \right)\hat\Rb^i,
\end{equation}
where $\acirc\yb^l$ are the coordinates of the configuration $Y(\acirc\eb^l)$,
and $\hat\yb^i$ define $Y(\hat\eb^i)$.

\begin{figure}
    \centering
    \includegraphics[width=0.49\linewidth]{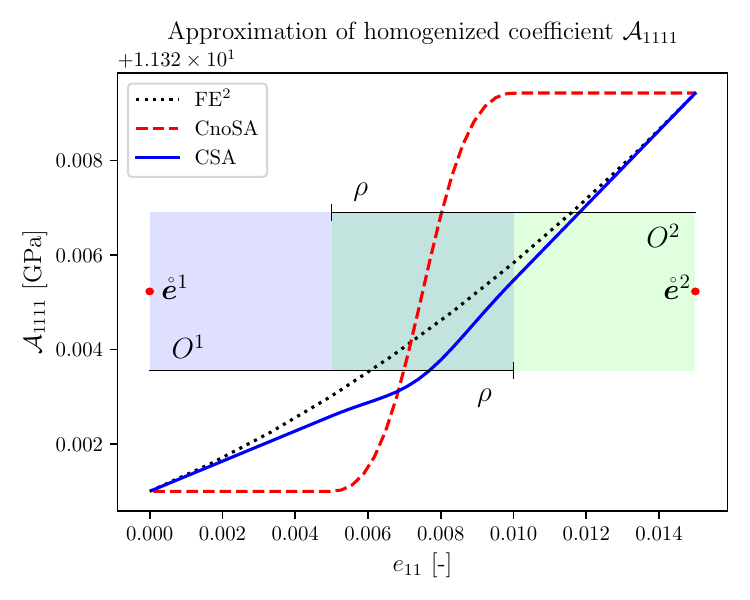}\hfil
        \includegraphics[width=0.49\linewidth]{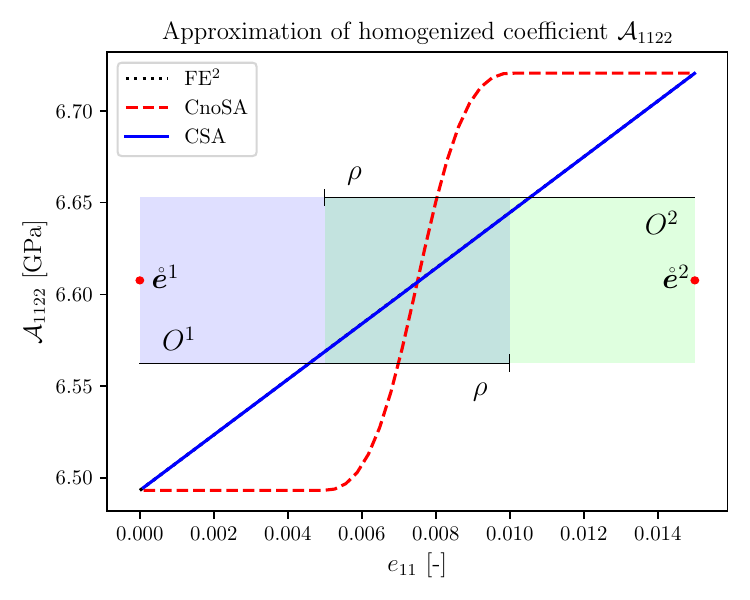}\\
    \includegraphics[width=0.49\linewidth]{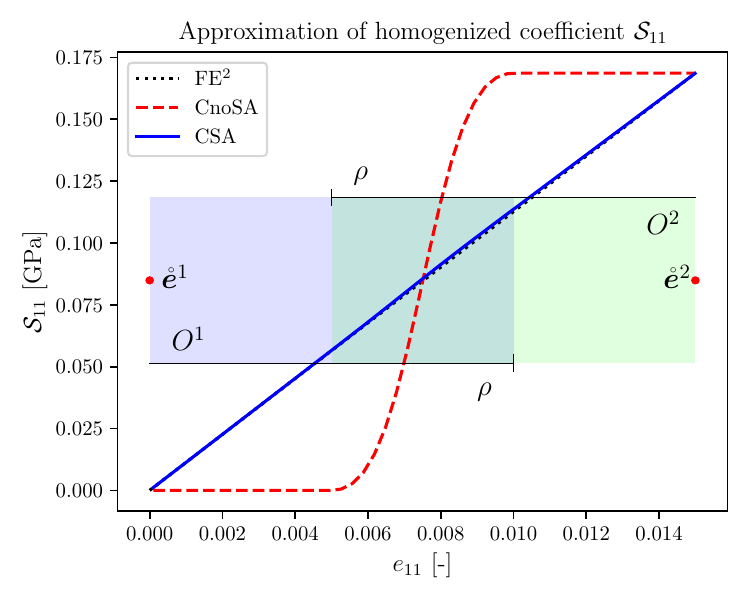}\hfil
        \includegraphics[width=0.49\linewidth]{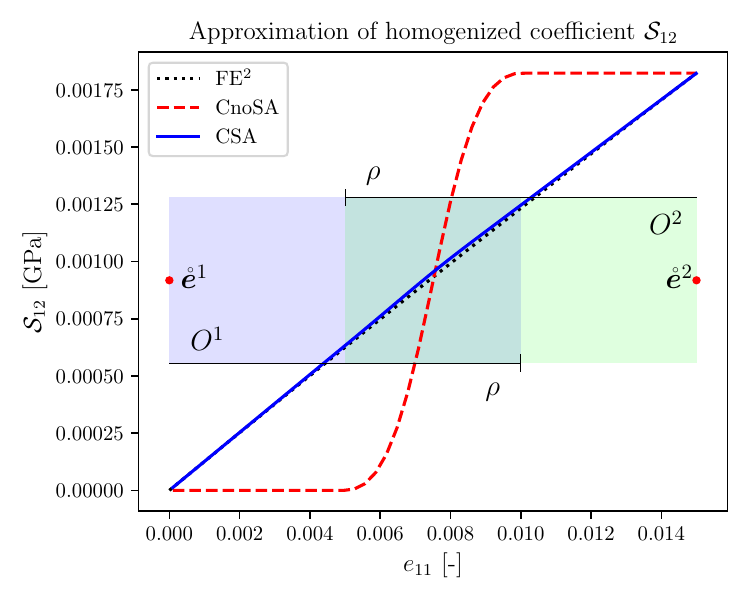}\\
    \caption{Approximations of $\Ahom$, $\Shom$ between two centroids $O^1 = O([0, 0, 0], 0.01)$,
     $O^2 = O([0.015, 0, 0], 0.01)$ with $\rho = 0.01$.
     Comparison of the direct two-scale simulation (\dtss) -- black dotted line,
     approximation using coefficients sensitivities (\acron) -- blue solid line, and 
     approximation without sensitivity information
        (CnoSA, $\delta\acirc\Ahom=\boldsymbol{0}$, $\delta\acirc\Shom=\boldsymbol{0}$) -- red dashed curve.}
    \label{fig:approx-coefs}
\end{figure}

\begin{figure}
    \centering
    \includegraphics[width=0.49\linewidth]{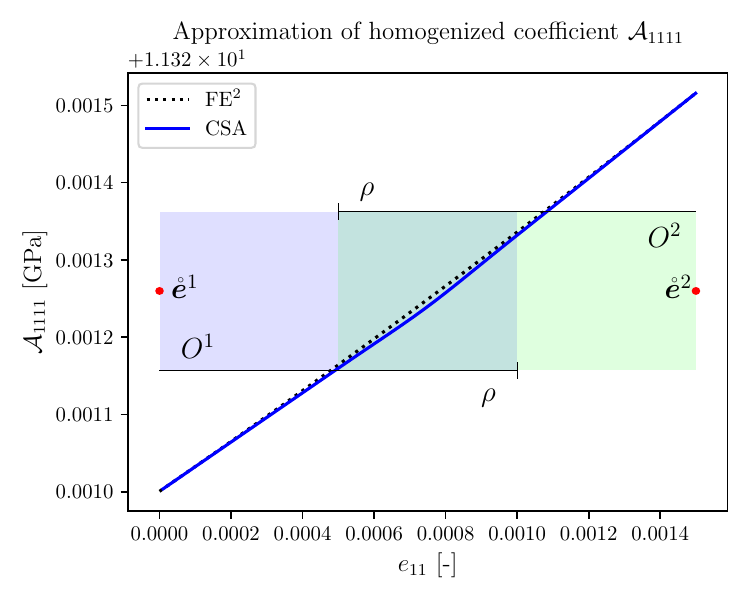}\hfil
        \includegraphics[width=0.49\linewidth]{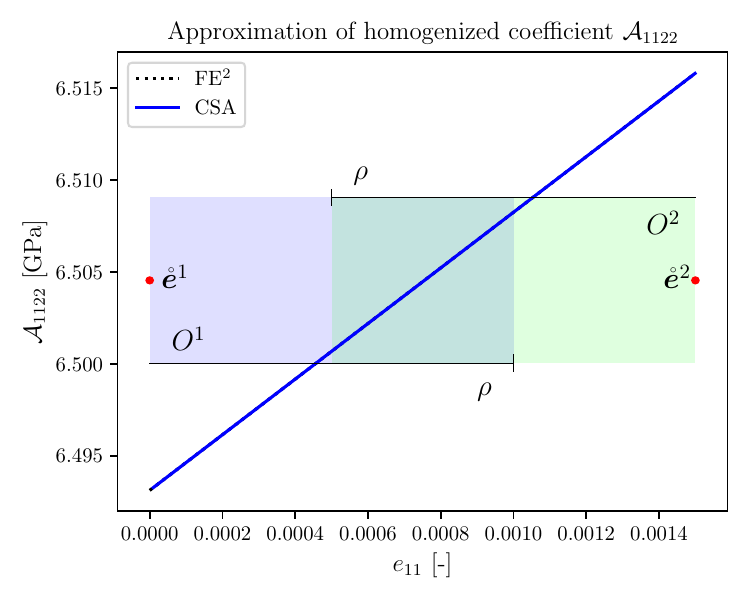}\\
    \includegraphics[width=0.49\linewidth]{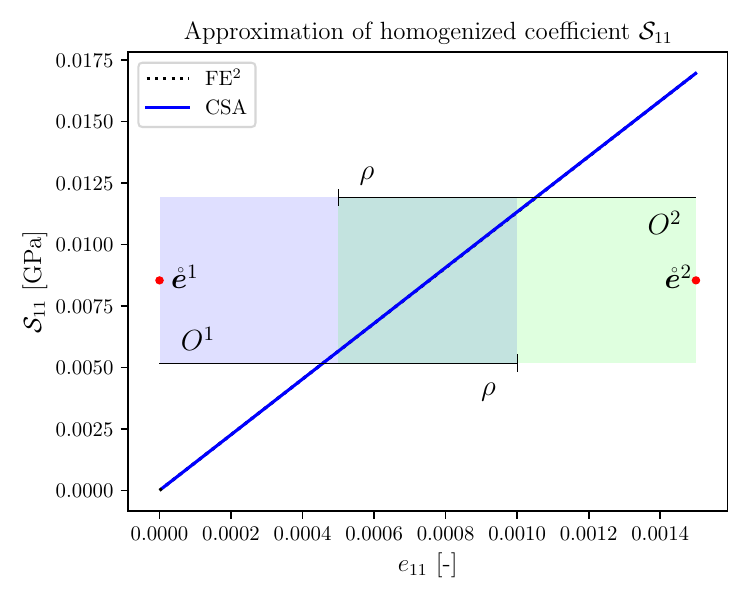}\hfil
        \includegraphics[width=0.49\linewidth]{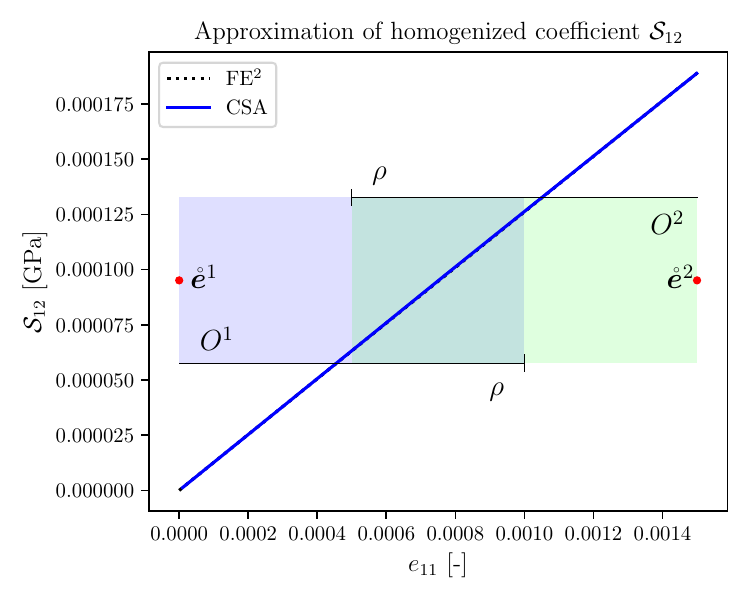}\\
    \caption{Approximations of $\Ahom$, $\Shom$ between two centroids $O^2 = O([0, 0, 0], 0.001)$,
     $O^1 = O([0.0015, 0, 0], 0.001)$ with $\rho = 0.001$.
     Comparison of the direct numerical simulation (\dtss) -- black dotted line and
     approximation using coefficients sensitivities (\acron) -- blue solid line.}
    \label{fig:approx-coefs2}
\end{figure}

To evaluate the difference between the approximated coefficients given by
\eq{eq:cf_approx} and those obtained from the direct numerical solution of
\eq{eq:micro}, a test was performed on a single periodic cell subjected to
macroscopic uniaxial strain. The periodic cell, whose geometry and material
composition are specified in Section~\ref{sec:num_simulations}, was gradually
loaded by the varying macroscopic strain $\eb = [e_{11}, 0, 0]$, for $e_{11} \in
\langle 0, 0.015\rangle$. Two different strategies were employed to approximate the
coefficients between two overlapping centroids $O^1 = O([0, 0, 0], 0.01)$ and
$O^2 = O([0.015, 0, 0], 0.01)$, both of size $\rho = 0.01$. The first strategy
uses the coefficient sensitivities $\delta\acirc\Ahom$,
$\delta\acirc\Shom$ and is represented by the blue solid curve ``\acron'' in
Figure~\ref{fig:approx-coefs}. The second strategy relies solely on the values
$\acirc\Ahom$, $\acirc\Shom$ in $O_1$ and $O_2$ (setting
$\delta\acirc\Ahom=\boldsymbol{0}$, $\delta\acirc\Shom=\boldsymbol{0}$ in
\eq{eq:cf_approx}) and its results are plotted in Figure~\ref{fig:approx-coefs}
by the red dashed curve ``CnoSA''. Both approximations are compared with the direct
two-scale simulation denoted by ``\dtss'' in graphs. The computations in
Fig.~\ref{fig:approx-coefs} were performed for a relatively large centroid size
$\rho$ to demonstrate the properties of the proposed approximation. The
corresponding results for smaller centroids, $O^1 = O([0, 0, 0], 0.001)$ and
$O^2 = O([0.0015, 0, 0], 0.001)$ with $\rho = 0.001$, are presented in
Figure~\ref{fig:approx-coefs2}.

\section{Computational algorithm for reduced two-scale simulation }\label{sec:mor-cluster-algorithm}

At the beginning of the two-scale simulation, the macroscopic configuration is
undeformed, so that $\Fb(\hat{x}^i) = \Ib$ at all quadrature points $\hat{x}^i,
i = 1, \dots, \hat{N}$. In this initial state, the microscopic problem is solved
for a single centroid $O^1 = O(\boldsymbol{0}, \rho)$, and the homogenized
coefficients $\Ahom(\hat{x}^i)$, $\Shom(\hat{x}^i)$ are constant across the
macroscopic domain.
As the macroscopic domain deforms, some strains get out of the existing
centroid, and the K-mean clustering algorithm is invoke to generate the minimal
number of new centroids required to cover all strains $\hat\eb^i$. The
microscopic subproblems are then solved only for those strains $\acirc\eb$
corresponding to the newly created centroids. Finally, based on the values
calculated for centroids, the homogenized coefficients are evaluated for all
quadrature points. These steps are repeated at each iteration of the
macroscopic solution.

\begin{figure}
    \centering
    \includegraphics[width=0.99\figscale]{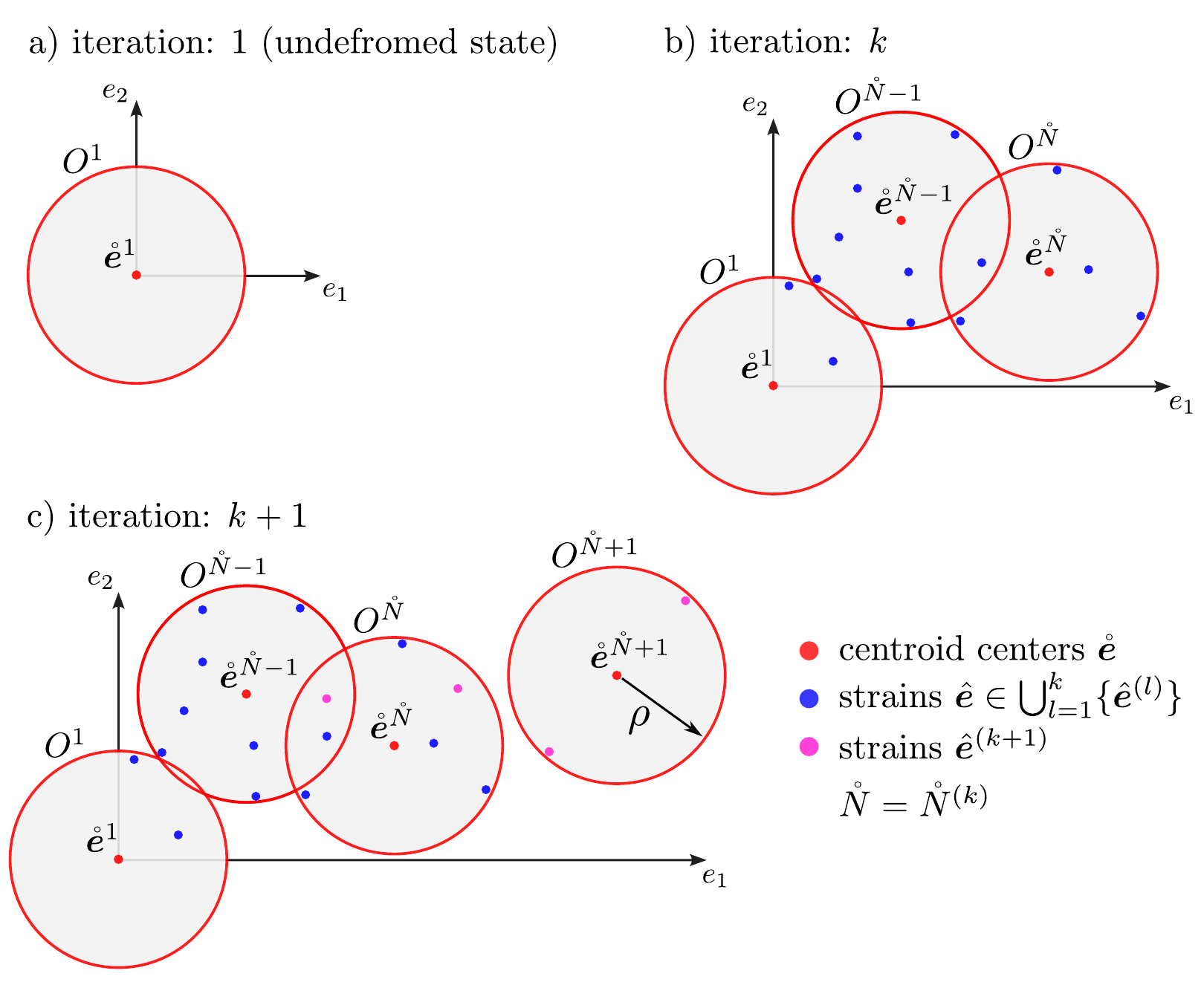}
    \caption{Creation of new centroids:
        a) initial centroid associated with the undeformed macroscopic state
        b) centroids in the $k$-th iteration
        c) a new centroid created in the $(k+1)$-th iteration to cover all strains
        $\hat\eb^{(k+1)}$.}
    \label{fig:centroid_creation}
\end{figure}

The procedure of centroid creation is illustrated in
Fig.~\ref{fig:centroid_creation}, where the deformation space is represented by
two dimensions for clarity, however it is three dimensional for 2D problems.
Fig.~\ref{fig:centroid_creation}a depicts the initial state, while
Fig.~\ref{fig:centroid_creation}b and \ref{fig:centroid_creation}c show two
subsequent iterations.
Assume that at the $k$-th macroscopic iteration
(Fig.~\ref{fig:centroid_creation}b) all strains $\hat\eb^i$ from previous
iterations are covered by centroids $O^j$, $j = 1, \dots, {\acirc N}^{(k)}$,
which are associated with the deformed microscopic configurations $\acirc Y^j$
and coefficients $\acirc\Ahom^j$, $\acirc\Shom^j$, $\delta\acirc\Ahom^j$,
$\delta\acirc\Shom^j$. In the $(k + 1)$-th iteration
(Fig.~\ref{fig:centroid_creation}c), new strains $\hat\eb^{i,(k+1)}$ arise due to the
updated macroscopic state. Some of these strains may fit in the existing
${\acirc N}^{(k)}$ centroids, while others require the creation of $\acirc
N^{new}$ additional centroids to satisfy \eq{eq:sets_O}. For these new centroids
$O^{\acirc{N} + l}$, $l = 1, \dots, \acirc N^{new}$, microscopic subproblems
must be solved on the deformed configurations $\acirc{Y}^l$ to evaluate coefficients
$\acirc\Ahom^l$, $\acirc\Shom^l$, $\delta\acirc\Ahom^l$, $\delta\acirc\Shom^l$.

Due to the nonlinear behavior of the microstructure, the deformed
configurations $\acirc Y^l = Y(\acirc\eb^l)$ cannot be obtained by a ``jump''
deformations $\acirc\eb^l$ but must be reached incrementally to maintain
equilibrium at the microscopic level. To achieve this efficiently,
the closest existing configuration $\acirc Y^m =Y(\acirc\eb^m)$,
\begin{equation}\label{eq:index_m}
    m = \argmin\limits_{i\in \{1,\dots, \acirc N^{(k)}\}} \vert \acirc\eb^i - \acirc\eb^l \vert,
\end{equation}
is selected as a reference, and the update from $\acirc Y^m$ to
$\acirc Y^l$ is done via \eq{eq:micro-update}, where $\nabla_x \ub^0$ is replaced
by the relative deformation $\gb^{lm}$ defined in \eq{eq:rel_g}.

\begin{figure}
    \centering
    \includegraphics[width=0.65\figscale]{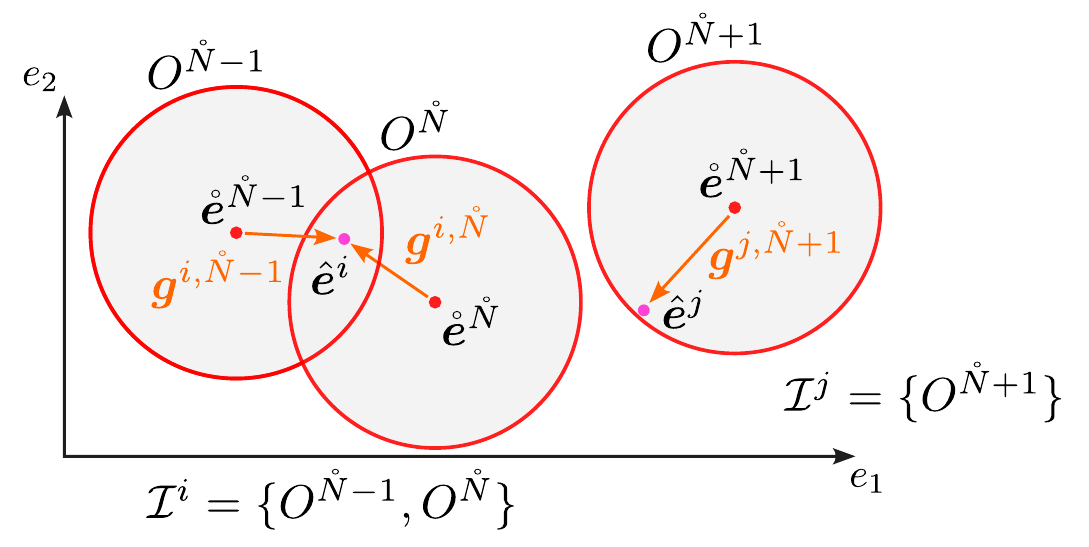}
    \caption{Illustration of homogenized coefficients approximation at points $\hat\eb^i$ and $\hat\eb^j$,
    index sets $\ihset{i}$, $\ihset{j}$, and relative strains $\gb^{il}$, $\gb^{jm}$, where $l\in\ihset{i}$, $m\in\ihset{j}$.}
    \label{fig:centroid_approximation}
\end{figure}

The computational procedure for the reduced two-scale simulation is summarize in
Algorithms~\ref{alg:csa-macro}, \ref{alg:csa-macro-approx},
and~\ref{alg:csa-micro}. Replacing {\sc CSA\_coefs\_approx} function with
{\sc FE2\_coefs} function (defined in  Algorithms~\ref{alg:fe2-coefs}) in
Algorithm~\ref{alg:csa-macro} yields the full FE$^{2}$ computational scheme.

\begin{algorithm}
    \caption{Reduced two-scale nonlinear simulation -- {\sc Macro}}\label{alg:csa-macro}
    \begin{algorithmic}
        \Function{Macro}{}
            \State initial configuration:
            \State \qquad $\Om^{(0)}$, $Y^{(0)}$,
                $\ub^{0, (0)} \algas \boldsymbol{0}$,
                $\acirc N \algas 0$
            \Statex
            \For{time steps $k=0,\dots, N_t-1$}
                \State given: $\Om^{(k)}$, $\ub^{0,(k)}$, load $\fbm_{ext}^{(k+1)}$
                \State $\ub^{0,(k+1)} \algas \ub^{0,(k)}$,
                    $i \algas 1$
                \Statex
                \Do
                    \State calculate homogenized coefficients:
                    \State\qquad $\Shom$, $\Ahom$
                            = \Call{CSA\_coefs\_approx}{$\nabla_x \ub^{0, (k + 1)}$}
                    \Statex
                    \State solve \eq{eq:macro} for displacement increments $\delta\ub^{0,(k + 1,i)}$
                    \State update displacements: $\ub^{0,(k+1)} \algas \ub^{0,(k + 1)} + \delta\ub^{0,(k + 1,i)}$
                    \State evaluate residue: $R \algas \vert \fb_{ext}^{(k+1)} - \fb_{int}^{(k,i)} \vert$
                    \State $i \algas i + 1$
                \doWhile{$R > \epsilon_M$}
                \Statex
                \State step increment: $\delta\ub^{0,(k+1)} = \ub^{0,(k+1)} - \ub^{0,(k)}$
                \State update macroscopic state: $\Om\tkp \algas \Om\tk + \{\delta\ub^{0,(k+1)}\}$
            \EndFor
        \EndFunction
    \end{algorithmic}
\end{algorithm}

\begin{algorithm}
    \caption{Reduced two-scale simulation -- {\sc CSA\_coefs\_approx}}\label{alg:csa-macro-approx}
    \begin{algorithmic}
        \Function{CSA\_coefs\_approx}{$\nabla_x \ub$}
            \State $\mathrm{E} \algas \emptyset$
            \For{$i=1, \dots, \hat{N}$}
                \State get $\hat\Ub^i, \hat\Rb^i$ from polar decomposition: $\Fb(\hat{x}^i) = \Ib + \nabla_x \ub(\hat{x}^i) = \hat\Rb^i\hat\Ub^i$
                \State get symmetric strain tensor: $\hat\eb^i = \hat\Ub^i - \Ib$

                \If{$\hat\eb^i \notin \bigcup_{j=1}^{\acirc N} O^j$}
                    \State $\mathrm{E} \algas \mathrm{E} \cup \{\hat\eb^i\}$
                \EndIf
            \EndFor
            \Statex
            \State using {\it K-mean} method, create $\acirc N^{new}$ new centroids
                such that
            \State\qquad $\mathrm{E} \subset \bigcup_{j=1}^{\acirc N^{new}} O^{\acirc N + j}$
            \Statex
            \For{$l = (\acirc N + 1), \dots, (\acirc N + \acirc N^{new})$}
                \If{$\acirc N > 0$}:
                    \State find the closest centroid $O^m$ to $O^l$
                    \State evaluate relative deformation:
                        $\gb^{lm}$
                    \State solve microscopic subproblem:
                    \State \qquad $\acirc\Shom^l$, $\acirc\Ahom^l$,
                        $\delta\acirc\Shom^l$, $\delta\acirc\Ahom^l$,
                        $\acirc Y^l$, $\acirc\omegabf^{rs, l}$ $=$
                        \Call{Micro}{$\gb^{lm}$, $\acirc Y^m$, $\acirc\omegabf^{rs,m}$}        
                \Else
                    \State solve microscopic subproblem:
                    \State \qquad $\acirc\Shom^l$, $\acirc\Ahom^l$,
                        $\delta\acirc\Shom^l$, $\delta\acirc\Ahom^l$,
                        $\acirc Y^l$, $\acirc\omegabf^{rs, l}$ $=$ 
                        \Call{Micro}{$\boldsymbol{0}$, $\acirc Y^{(0)}$, $\boldsymbol{0}$}    
                \EndIf
            \State $\acirc N \algas \acirc N + \acirc N^{new}$
            \EndFor
            \Statex
            \For{$i=1, \dots, \hat{N}$}
                \State approximate coefficients $\Shom(\hat{x}^i)$, $\Ahom(\hat{x}^i)$ using 
                    $\acirc\Shom$, $\acirc\Ahom$, $\delta\acirc\Shom$, $\acirc\delta\Ahom$, $\hat\Rb^i$ 

            \EndFor
            \Statex
            \State \Return{$\Shom$, $\Ahom$}
        \EndFunction
    \end{algorithmic}
\end{algorithm}

\begin{algorithm}
    \caption{Reduced two-scale nonlinear simulation -- {\sc Micro}}\label{alg:csa-micro}
    \begin{algorithmic}
        \Function{Micro}{$\gb$, $Y$, $\omegabf^{rs}$}
            \State evaluate: $\delta \wb = (\Pibf^{rs} + \omegabf^{rs}) g_{rs}$
            \State update configuration: $Y\snew \algas Y + \{\delta\wb\}$
            \Statex
            \Do
                \State solve nonlinear problem \eq{eq:micro-equilibrium} for corrections $\delta\wb$ 
                \State update configuration: $Y\snew \algas Y\snew + \{\delta\wb\}$
                \State evaluate residue: $R  \algas \vert \delta\wb\vert$
            \doWhile{$R > \epsilon_\#$}
            \Statex
            \State solve \eq{eq:micro} for new characteristic responses $\omegabf^{rs,\new}$
            \State evaluate: $\Shom$, $\Ahom$, $\delta\Shom$, $\delta\Ahom$
            \Statex
            \State \Return{$\Shom$, $\Ahom$, $\delta\Shom$, $\delta\Ahom$, $Y\snew$, $\omegabf^{rs,\new}$}
        \EndFunction
    \end{algorithmic}
\end{algorithm}

\begin{algorithm}
    \caption{Two-scale nonlinear simulation -- {\sc FE2\_coefs}}\label{alg:fe2-coefs}
    \begin{algorithmic}
        \Function{FE2\_coefs}{$\nabla_x \delta\ub$}           
            \For{$i=1, \dots, \hat{N}$}
                \State $\Shom(\hat{x}^i)$, $\Ahom(\hat{x}^i)$, $Y(\hat{x}^i)$, $\omegabf^{rs}(\hat{x}^i)$ = 
                    \Call{Micro}{$\nabla_x \delta\ub(\hat{x}^i)$, $Y(\hat{x}^i)$, $\omegabf^{rs}(\hat{x}^i)$}
            \EndFor
            \Statex
            \State \Return{$\Shom$, $\Ahom$}
        \EndFunction
    \end{algorithmic}
\end{algorithm}


\section{Model reduction using proper orthogonal decomposition}\label{sec:mor-pod}

To assess the efficiency of the proposed \acron approach against a
well-established technique, we implemented the R3M algorithm described in
\cite{Yvonnet_2007}, which is based on the proper orthogonal decomposition (POD)
method \cite{Chatterjee_2000} applied to the microscopic subproblems. The R3M
method constructs a reduced basis from precomputed solutions of the microscopic
problem corresponding to various macroscopic deformation states. This basis is
then used to project the discretized microscopic problem onto a reduced-order
model of significantly lower dimensionality, resulting in a substantial
reduction of computational cost.

Discretization of the microscopic systems \eq{eq:micro} and~\eq{eq:micro-equilibrium2} by means of the
finite element method leads to a system of discrete equations,
\begin{equation}\label{eq:pod-discrete}
    \Kbm_\# \qbm = \fbm_\#,
\end{equation}
which is solved iteratively due to dependence of the tangent matrix $\Kbm_\#$ on
the total deformation. Vector $\qbm$ of dimension $N_\#$
represents the characteristic responses $\omegabf^{ij}$ or microscopic corrections $\delta\wb$ after finite element
discretization, where $N_\#$ denotes the number of degrees of freedom of the
microscopic subproblem.
To construct the reduced basis, the microscopic problem \eq{eq:pod-discrete} is solved for a
series of deformed configurations governed by the macroscopic strain
%
\begin{equation}
    (\eb^{pqR})_{rs} (k) = k / (N_t - 1) \epsilon^{pqR} \delta_{pr} \delta_{qs},
        \quad \mbox{for } k = 0, \dots, N_t - 1,
        \quad p, q = 1, 2, (3), \quad R = \mbox{min}, \mbox{max},
\end{equation}
%
where $\epsilon^{pq,min}$, $\epsilon^{pq,max}$ are the estimated minimum and
maximum values of the strain components $p, q$, and $N_t$ is a given number of
load steps. For each load case $\alpha = pqR$, $N_t$ solutions
$\qbm^{\alpha}(k)$ of \eq{eq:pod-discrete} are calculated
and collected in matrix $\Ubm^{\alpha}$,
%
\begin{equation}
    \Ubm^{\alpha} = [\qbm^{\alpha}(0), \qbm^{\alpha}(1), \dots, \qbm^{\alpha}(N_t - 1)].
\end{equation}
%
Matrices $\Ubm^{\alpha}$ are further collected in 
$\Vbm = [\Ubm^{11,min}, \Ubm^{11,max}, \dots, \Ubm^{33,max}]$ 
and the eigenvalue problem
\begin{equation}\label{eq:pod-evp}
\Vbm \Vbm^T \phibf_m = \lambda_m \phibf_m
\end{equation}
is solved for the eigenvectors $\phibf_m$ and eigenvalues $\lambda_m$. The
eigenvectors corresponding to the first $M(\delta)$ largest eigenvalues form the
reduced orthogonal basis $\Phibf = [\phibf_1, \dots, \phibf_M]$ where $M$ is
chosen according to a given tolerance error parameter $\delta$
\begin{equation}\label{eq:pod-err}
    \frac{\left(\sum_{l=M+1}^{N_\#} \lambda_l \right)^{1/2}}{\left(\sum_{l=1}^{N_\#} \lambda_l \right)^{1/2}}
        < \delta,
    \end{equation}
see \cite{Yvonnet_2007} for details.

Using the pre-calculated basis $\Phibf$, the following approximation
for correctors $\qbm$ is used:
\begin{equation}\label{eq:pod-approx}
    \qbm = \sum\limits_{m=1}^{M} \phibf_m \xi_m = \Phibf \xibf.
\end{equation}
Substituting \eq{eq:pod-approx} into \eq{eq:pod-discrete}, the reduced discrete system of
dimensions $M\times M$ for the unknown multipliers $\xibf$ is obtained,
\begin{equation}\label{eq:pod-reduced}
    \Phibf^T \Kbm_\# \Phibf \xibf = \Phibf^T\fbm_\#.
\end{equation}
Since $M p\ll N_{\#}$, the solution of  \eq{eq:pod-reduced} exhibits a
substantially lower computational complexity compared to the full-order problem
\eq{eq:pod-discrete}. This gain in the ``on-line'' phase is partially offset by
the computational cost associated with the solution of a large-scale eigenvalue
problem for the dense correlation matrix $\Vbm \Vbm^T$ of dimension $N_{\#}
\times N_{\#}$. The efficiency of the reduced basis construction could be
further improved through various approaches, e.g. by adopting more effective
sampling strategies for the macroscopic deformation space, such as Latin
hypercube sampling. A detailed discussion of these improvements is beyond the
scope of the present work.

%% file: part_num.tex
\section{Numerical simulations}\label{sec:num_simulations}

The coupled two-scale numerical simulations have been carried out using the \SfePy package -- Simple Finite Elements in Python, \cite{Cimrman_2014}. As its
name suggests, the package implements the finite element (FE) method in Python
and provides useful tools for the efficient realization of coupled or uncoupled
multiscale simulations, see \cite{Cimrman_Lukes_Rohan_2019}. The source code used to perform the
computations presented in this paper is published in the Mendeley Data
repository \cite{Lukes_code_2025}.

\begin{figure}
    \centering
    \hfil\includegraphics[width=0.3\figscale]{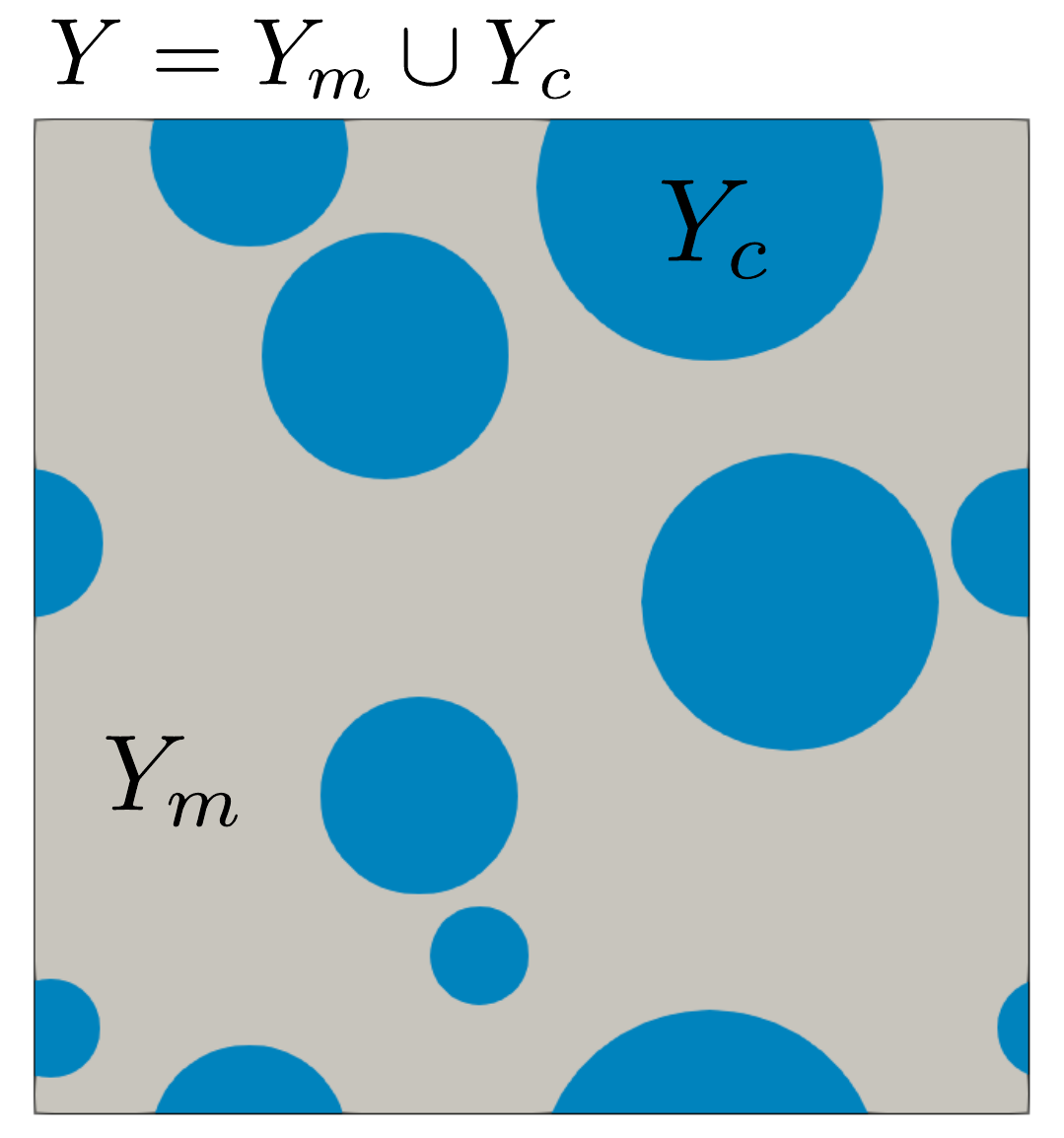}\hfil\hfil
    \includegraphics[width=0.6\figscale]{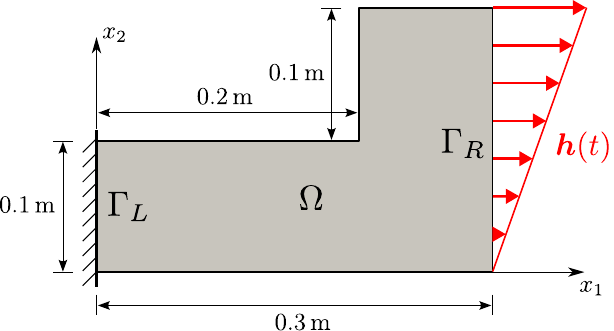}\hfil
    \caption{Left: decomposition of the microscopic reference cell $Y = Y_m \cup Y_c$;
             Right: undeformed macroscopic domain $\Omega$ and applied boundary conditions.}
    \label{fig:num_micro_macro}
\end{figure}

All subsequent simulations concern a two-dimensional problem of a heterogeneous
periodic structure composed of two {hyperelastic} materials with distinct stiffness: a soft ``matrix'' and
stiffer ``inclusions''. The material parameters are summarized in Table~\ref{tab:materials}.
\begin{table}
    \begin{tabular}{|l|c|c|}
        \hline
        Part & shear modulus $\mu$ [GPa] & bulk modulus $K$ [GPa]\\
        \hline\hline
        matrix (epoxy) $Y_m$ & 1.35 & 5.7 \\
        \hline
        inclusions (glass) $Y_c$ & 28.46 & 43.21 \\
        \hline
    \end{tabular}

    \caption{Material properties of the compressible neo-Hookean hyperelastic structure.}\label{tab:materials}
\end{table}
The periodic structure is represented by a reference unit cell $Y = Y_m\cup
Y_c$, where $Y_m$ denotes the soft matrix and $Y_c$ the stiff inclusions. The
undeformed configuration of the unit cell is shown in Fig.~\ref{fig:num_micro_macro}
left.
At the macroscopic level, an L-shaped sample is fixed on its
left {hand} side,
{$\ub^0 = \boldsymbol{0}$ on $\Gamma_L$}, and
subjected to a surface traction on the opposite side, 
{$\tilde\hb = \bar\hb(x, t)$ on $\Gamma_R$}.
The traction is applied incrementally in $N_t = 10$ equal loading steps,
see Figure~\ref{fig:num_micro_macro} right,
%
\begin{equation}\label{eq:num_loading_fun}
    \begin{aligned}
        \bar\hb(x, t_k) =
        \begin{bmatrix}
            10^{8}\,r_k\,x_2 / 0.2\\ 0
        \end{bmatrix}\,\mbox{N/m},
        \quad \mbox{for}\quad r_k = \frac{k}{N_t - 1},
        \quad t_k = r_k \cdot 1\mbox{s},\quad
        \mbox{and}\quad k = 0, \dots, N_t -1.
    \end{aligned}
\end{equation}

The displacement fields at both scales are approximated by the Lagrangian
piece-wise linear elements and the nonlinear systems defined
by \eq{eq:macro} and \eq{eq:micro} are solved using the Newton-Raphson
iterative scheme. Since the macroscopic FE mesh consists of 36 quadrilateral
elements, see Fig.~\ref{fig:num_res_displ} (left), and a four-point quadrature
rule is used for numerical integration in all simulations, $36 \times 4$
deformation-dependent homogenized coefficients must be evaluated in each
macroscopic iteration independently of the employed method (\dtss, \acron, or POD).
In Sections~\ref{sec:num_simul} and~\ref{sec:num_error}, the reference cell $Y$
is represented by a FE mesh (denoted as micro-mesh \#1 in
Tab.~\ref{tab:Y-cells}) consisting of 4166 triangular elements and 2160 nodes.
In Section~\ref{sec:num_comp_time}, several microscopic FE meshes with different
levels of refinement are employed to study the effect of mesh density on
computational time.

Section \ref{sec:num_simul} presents the results of the reduced two-scale
simulation and shows some numerical aspects of centroid formation. In Section
\ref{sec:num_comp_time}, the computational efficiency of the proposed algorithm
is evaluated and compared with the proper orthogonal decomposition approach
briefly described in Sec.~\ref{sec:mor-pod}. The approximation error is analyzed
in Section~\ref{sec:num_error}. The application of the CSA algorithm for a three
dimensional case is demonstrated in Section~\ref{sec:num-3D}.

\subsection{Reduced micro-macro simulation}\label{sec:num_simul}

\begin{figure}
    \centering
    \includegraphics[width=0.5\linewidth]{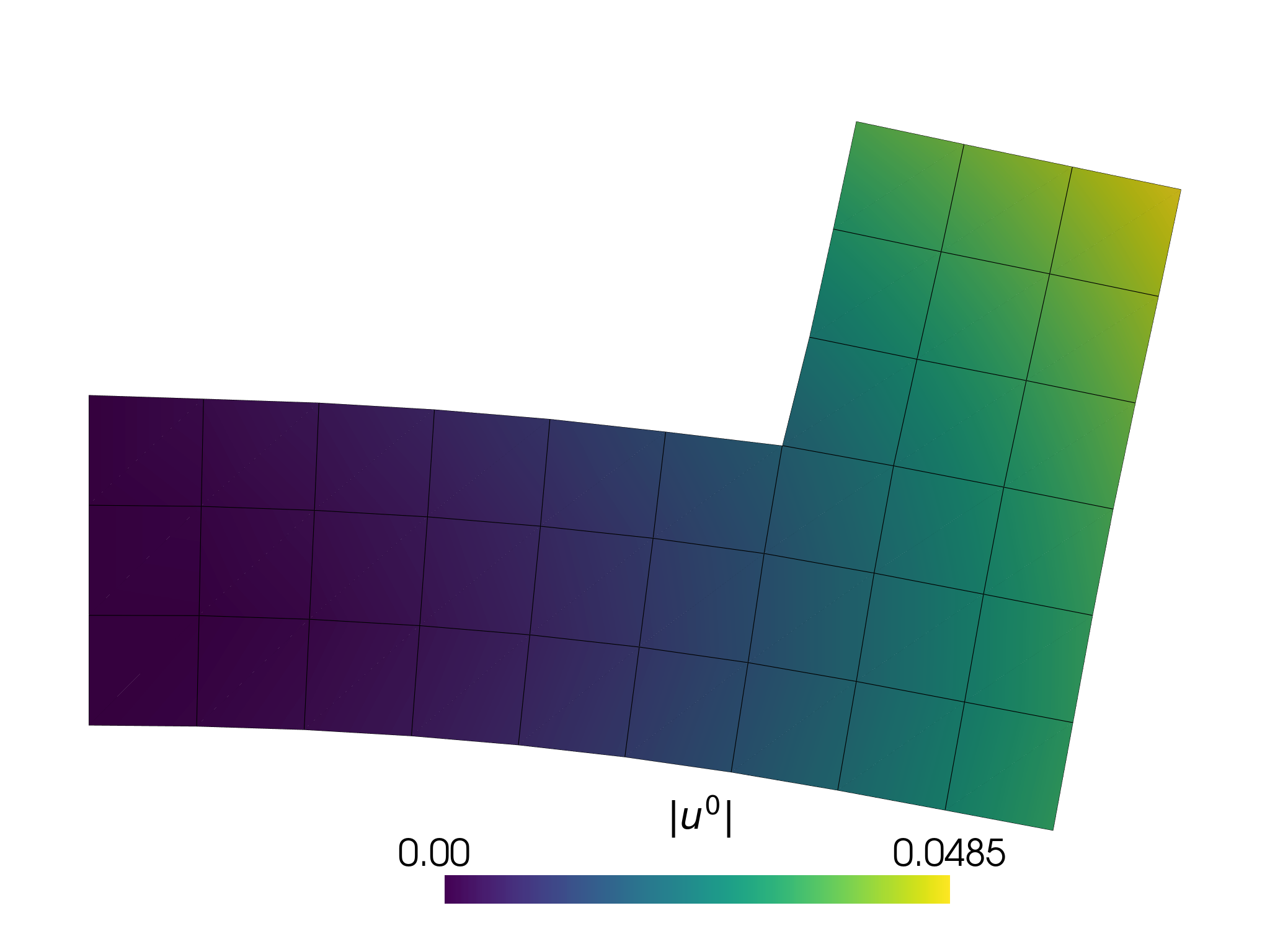}\hfil
    \includegraphics[width=0.49\linewidth]{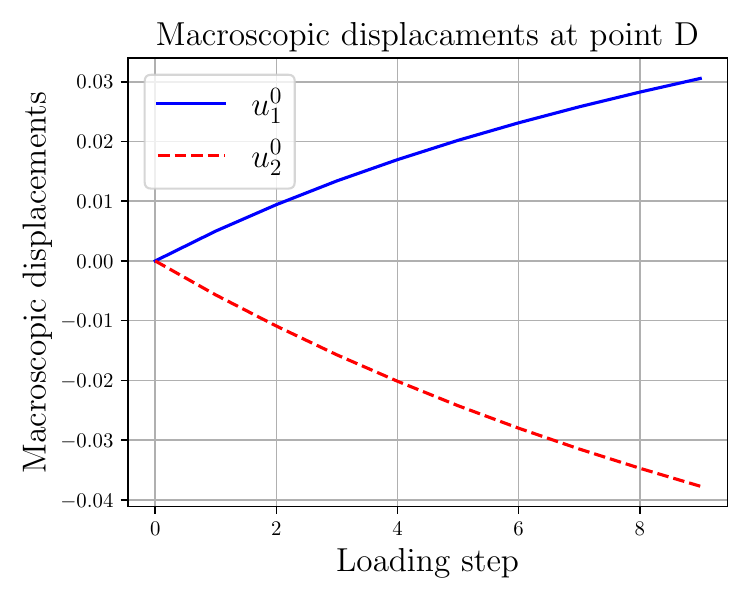}\hfil
    \caption{Left: resulting deformation in the last time step; Right: macroscopic displacements at point D.}
    \label{fig:num_res_displ}
\end{figure}

\begin{figure}
    \centering
    \includegraphics[width=0.49\linewidth]{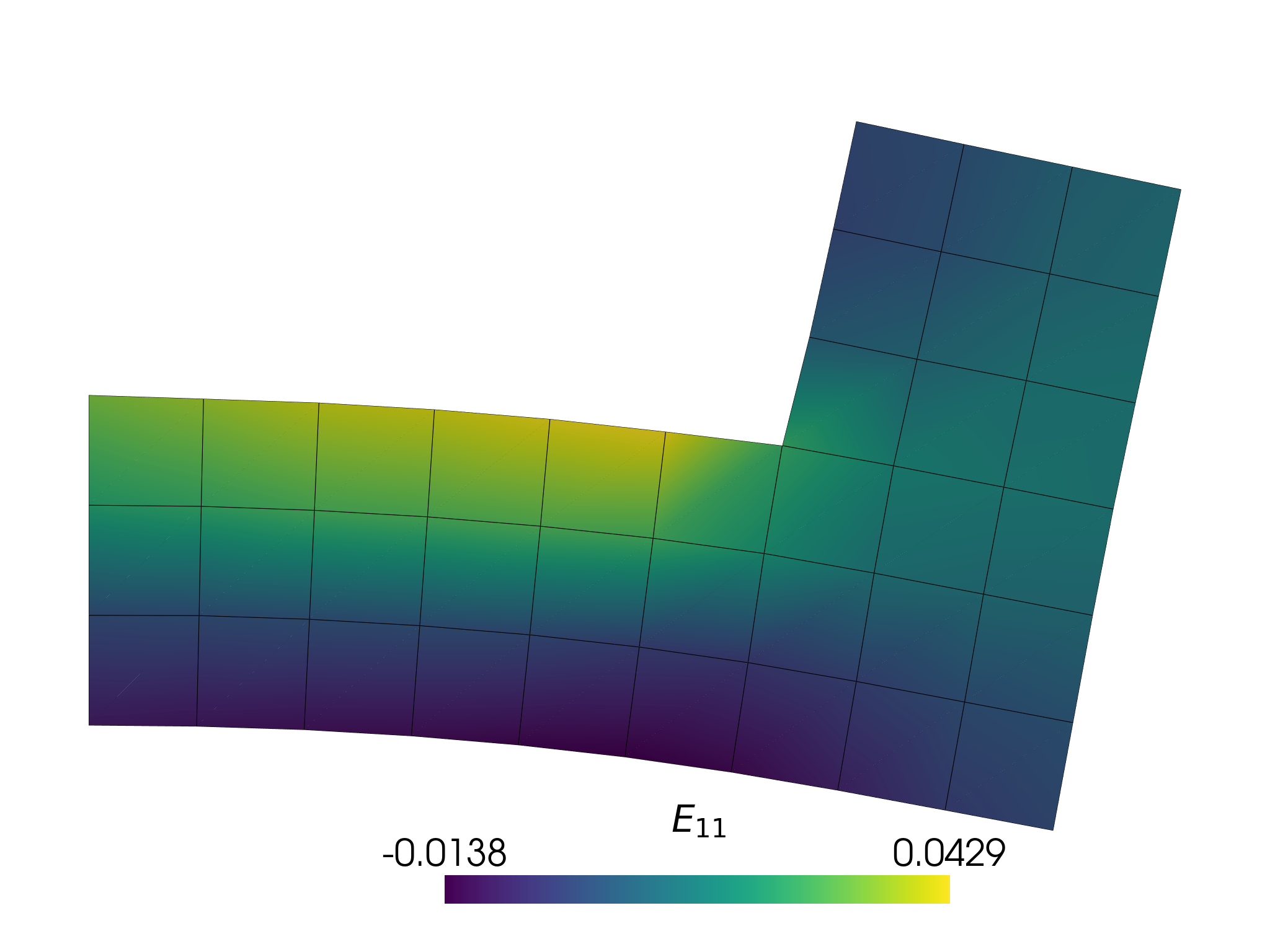}\hfil
        \includegraphics[width=0.49\linewidth]{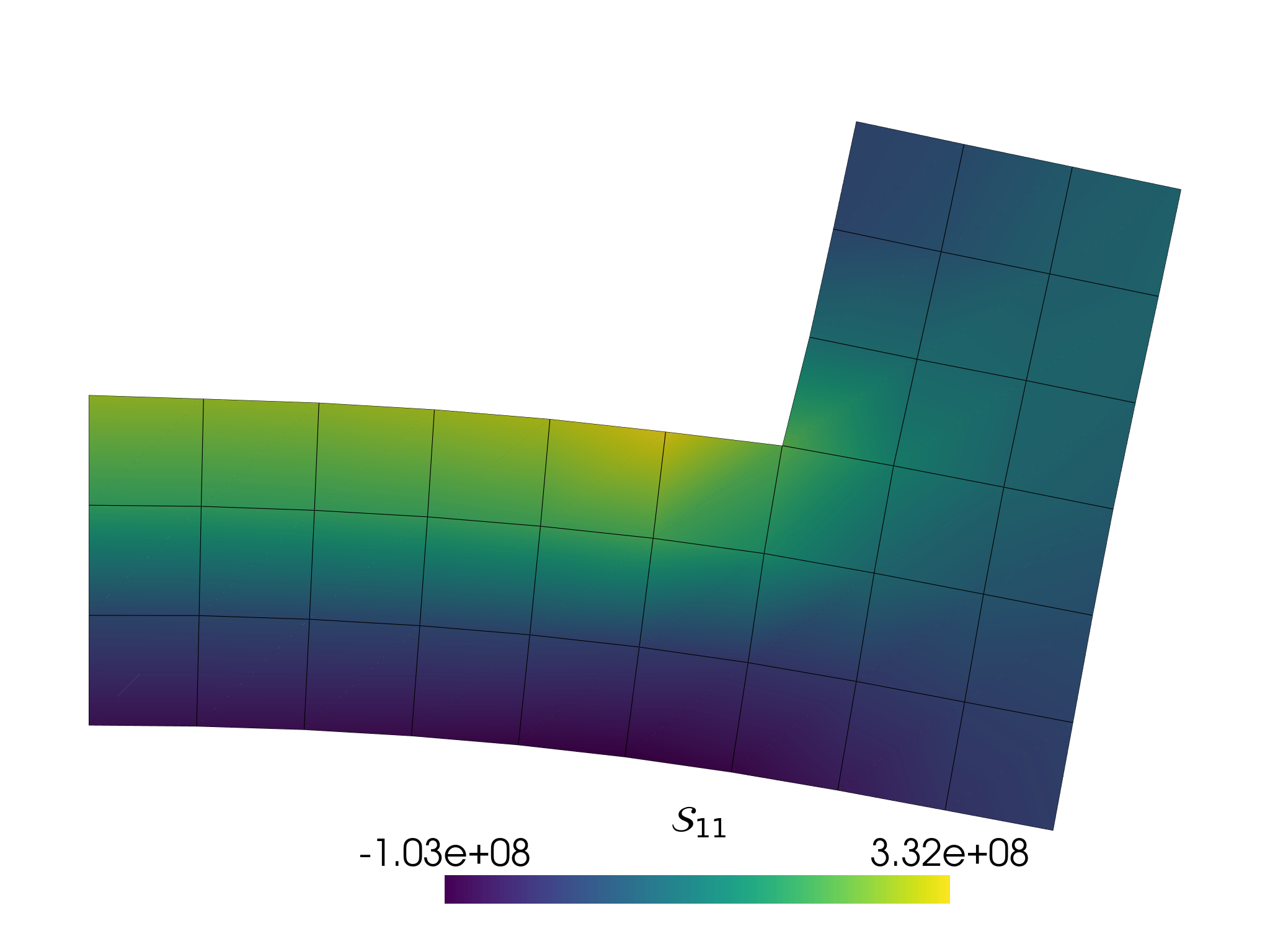}\\
    \includegraphics[width=0.49\linewidth]{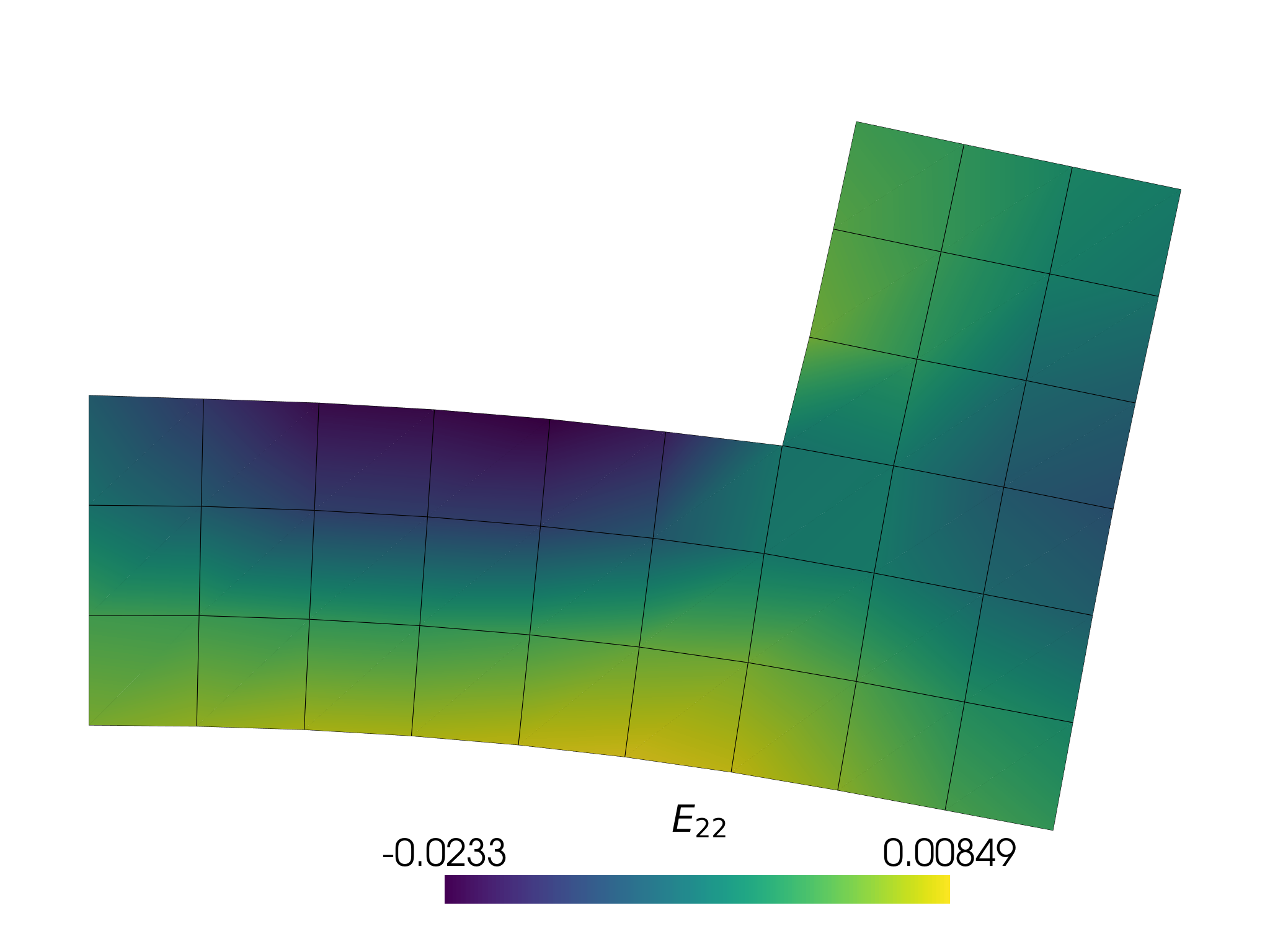}\hfil
        \includegraphics[width=0.49\linewidth]{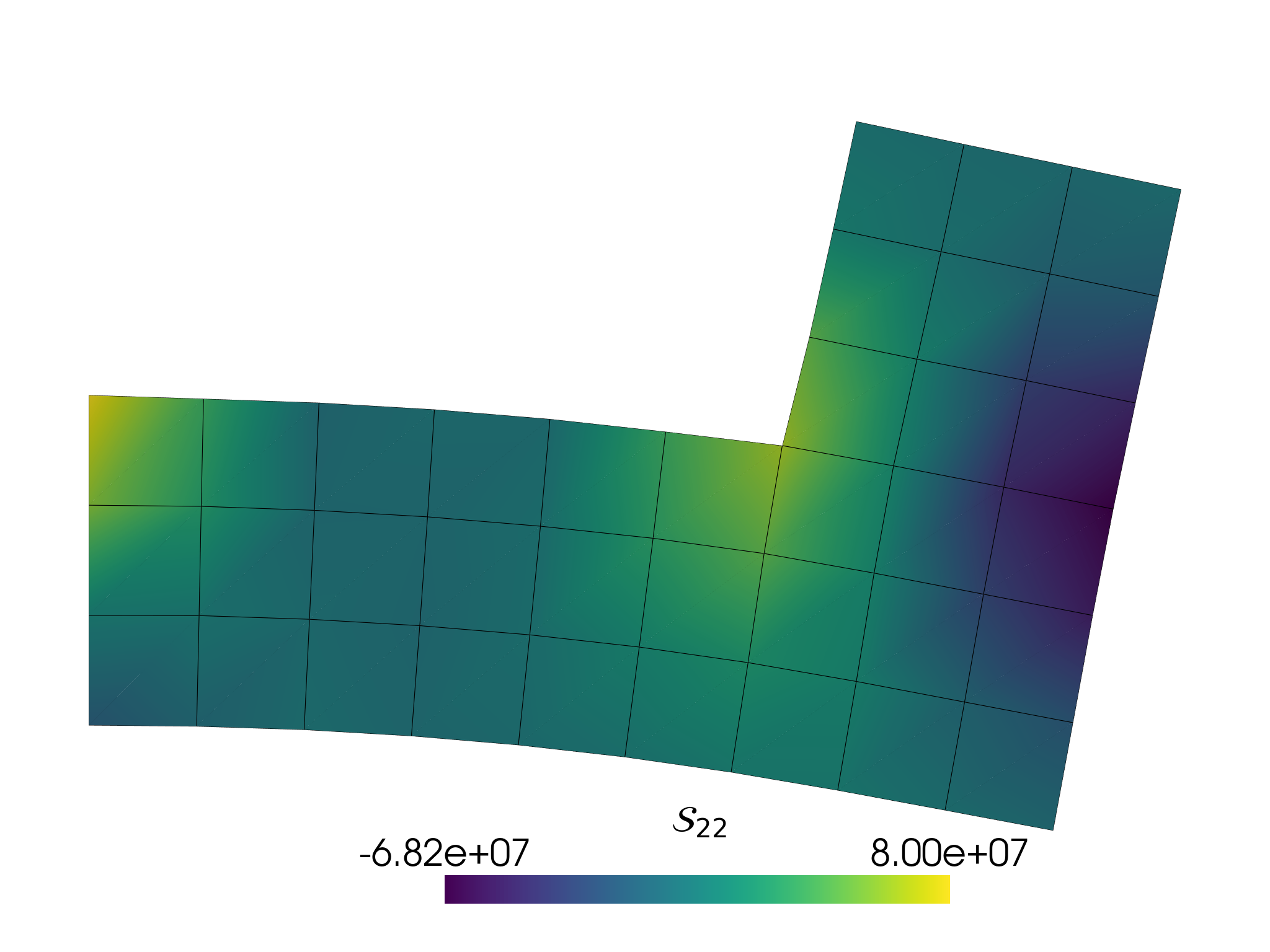}\\
    \caption{Left: components of the Green strain tensor;
             Right: components of the averaged Cauchy stress $\Shom$.}
    \label{fig:num_strain_stress}
\end{figure}

Running the coupled micro-macro simulation for $N_t$ loading steps yields the
deformed macroscopic sample depicted in Fig.~\ref{fig:num_res_displ} (left). The
time evolution of the displacements at point D is plotted in
Fig.~\ref{fig:num_res_displ} (right), where the nonlinear response to the linearly
applied load arises from the nonlinear kinematics and hyperelastic material
behavior. Figure~\ref{fig:num_strain_stress} shows the components of the Green
strain (left) and the averaged Cauchy stress (right) at the end of the
simulation.

\begin{figure}
    \centering
        \includegraphics[width=0.98\linewidth]{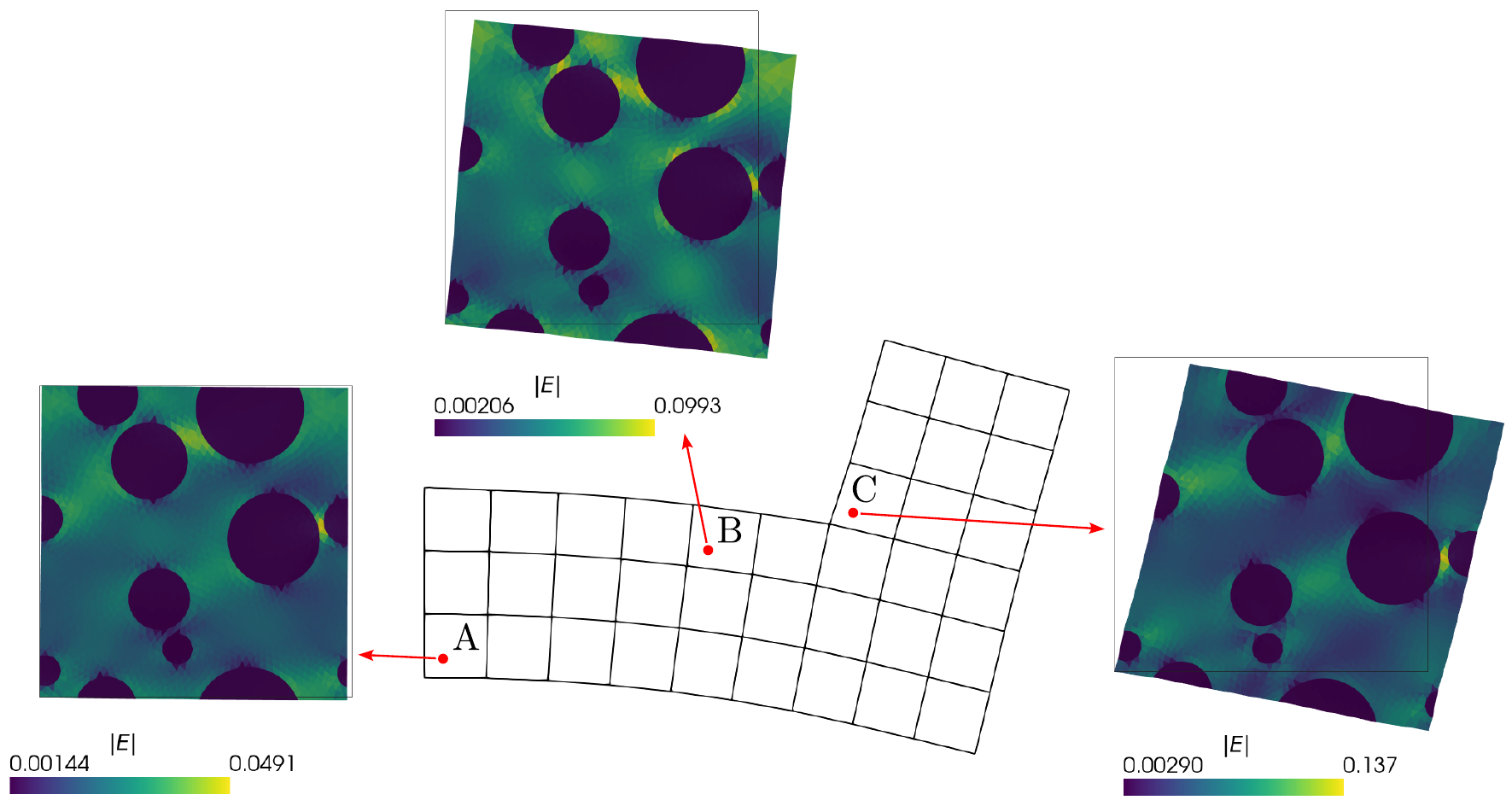}
    \caption{Reconstructed microstructures in the last loading step at quadrature points A, B, and C.}
    \label{fig:num_deform_micro_recovered}
\end{figure}

\begin{figure}
    \centering
    \includegraphics[width=0.49\linewidth]{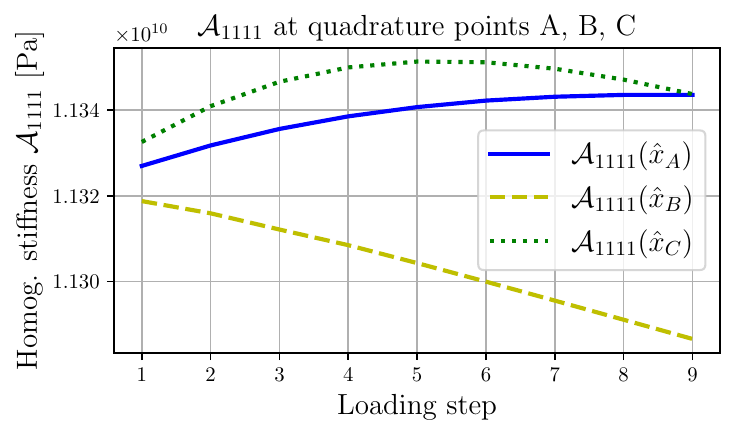}\hfil
        \includegraphics[width=0.49\linewidth]{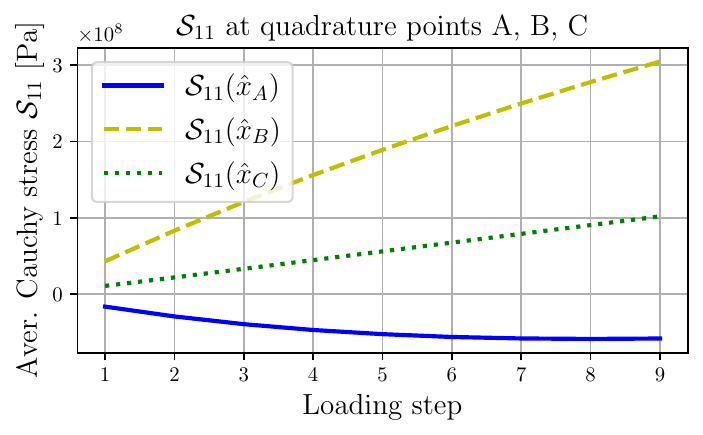}\\
    \includegraphics[width=0.49\linewidth]{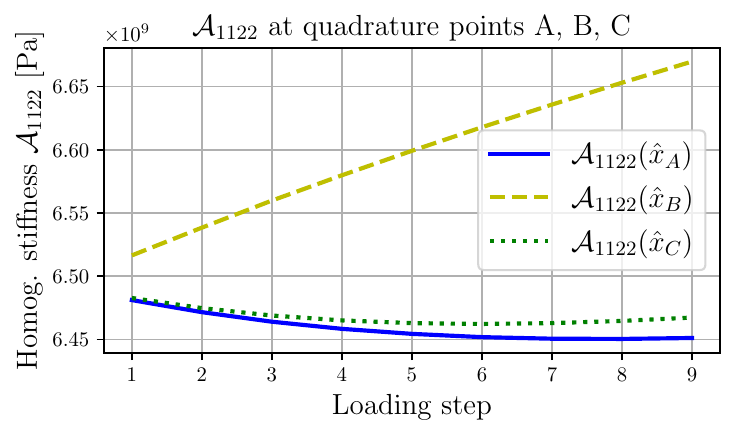}\hfil
        \includegraphics[width=0.49\linewidth]{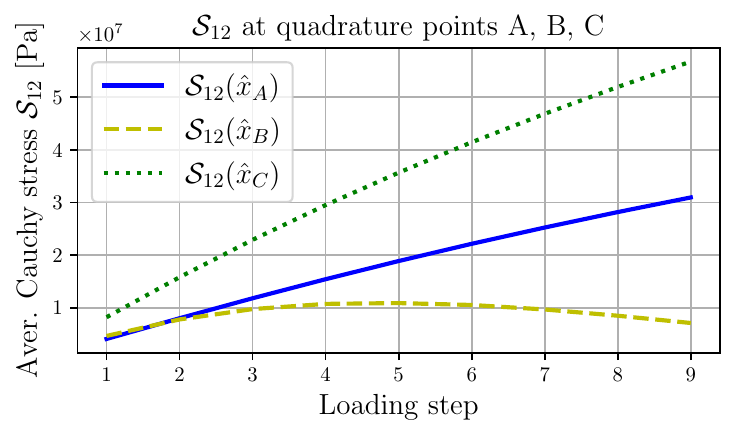}\\
    \includegraphics[width=0.49\linewidth]{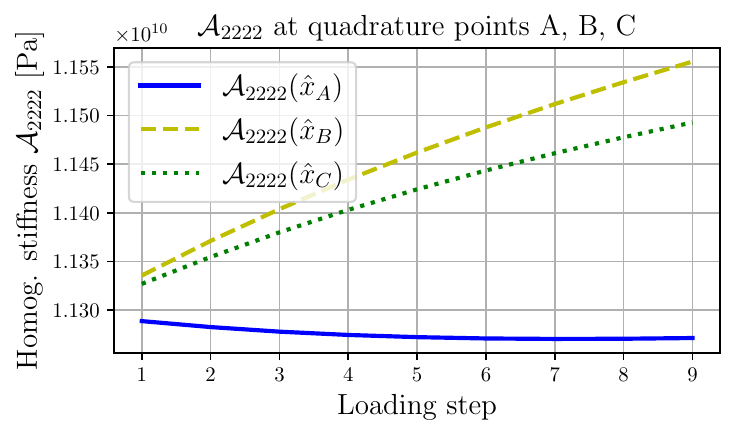}\hfil
        \includegraphics[width=0.49\linewidth]{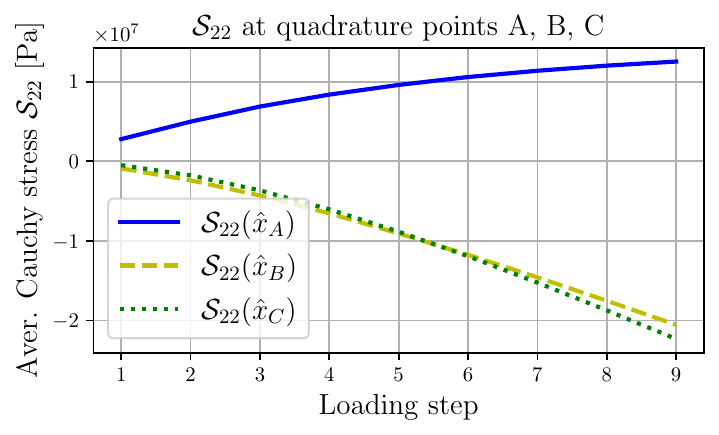}\\
    \caption{Evolution of the homogenized coefficients
              $\Ahom$ (left) and $\Shom$ (right) at macroscopic quadrature points A, B, C.}
    \label{fig:evolution_ABC}
\end{figure}

The reconstructed microstructures at three different points, A, B, and C,
corresponding to the quadrature points of the macroscopic sample, are depicted in
Fig.~\ref{fig:num_deform_micro_recovered}. The evolution of the homogenized
coefficients $\Ahom$, $\Shom$ at these points during the
loading steps is plotted in Figure~\ref{fig:evolution_ABC}.


\begin{figure}
    \centering
    \includegraphics[width=0.49\linewidth]{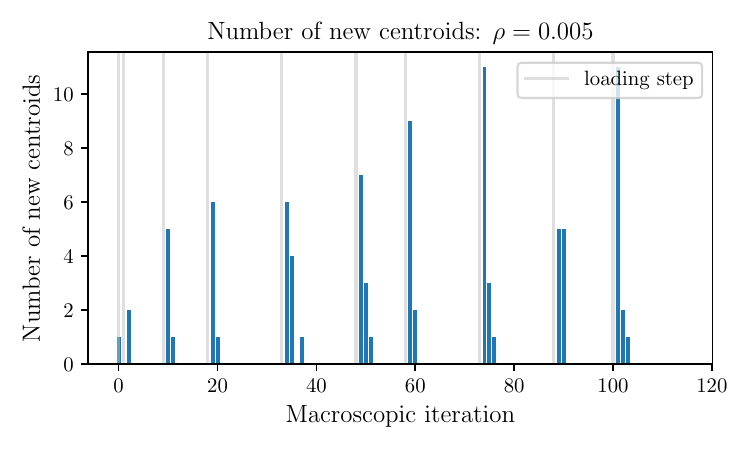}\hfil
        \includegraphics[width=0.49\linewidth]{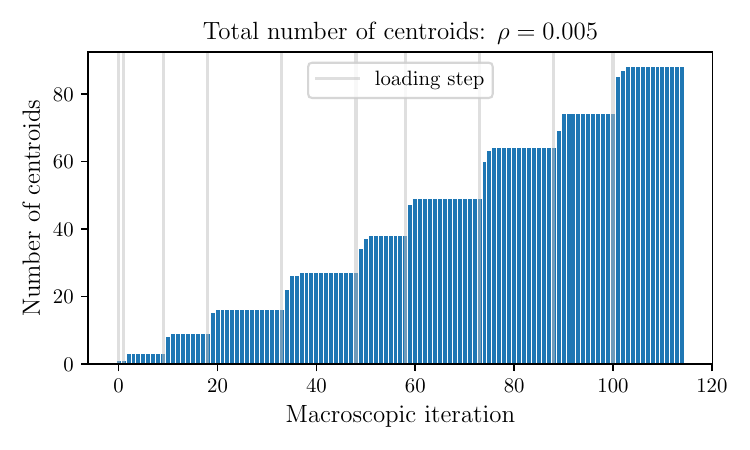}\\
    \includegraphics[width=0.49\linewidth]{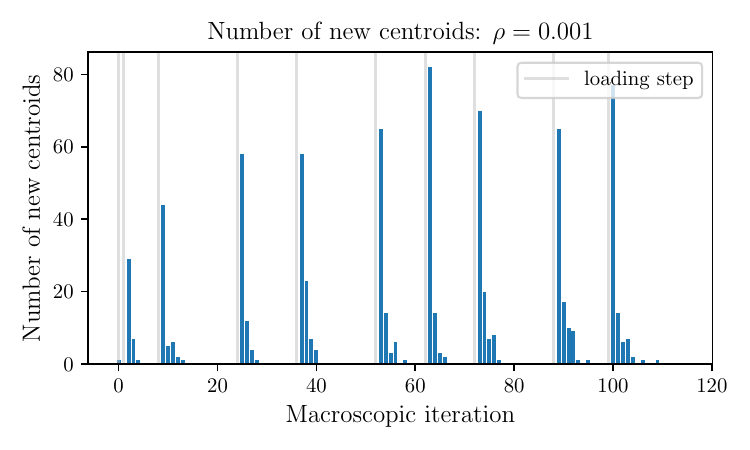}\hfil
        \includegraphics[width=0.49\linewidth]{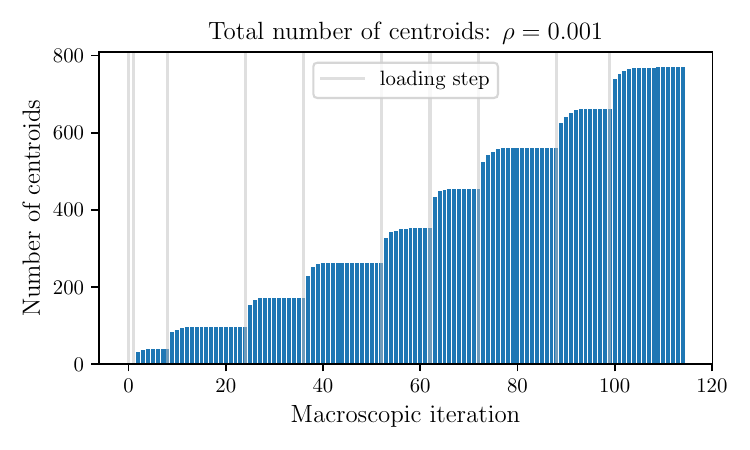}
    \caption{The number of new centroids in each macroscopic iteration step for $\rho = 0.005$ (top) and $\rho=0.001$ (bottom).}
    \label{fig:num_centroids_new}
\end{figure}

As described in Sec.~\ref{sec:mor-cluster}, centroids are generated dynamically
as the deformation progresses due to external loading. The number of new
centroids at each iteration step is shown in Fig.~\ref{fig:num_centroids_new}
(left), while the total (cumulative) number of centroids created throughout the
simulation is presented in Fig.~\ref{fig:num_centroids_new} (right). The results
are shown for two centroid sizes $\rho=0.005$ (top) and $\rho=0.001$ (bottom).

\begin{figure}
    \centering
    \includegraphics[width=0.49\linewidth]{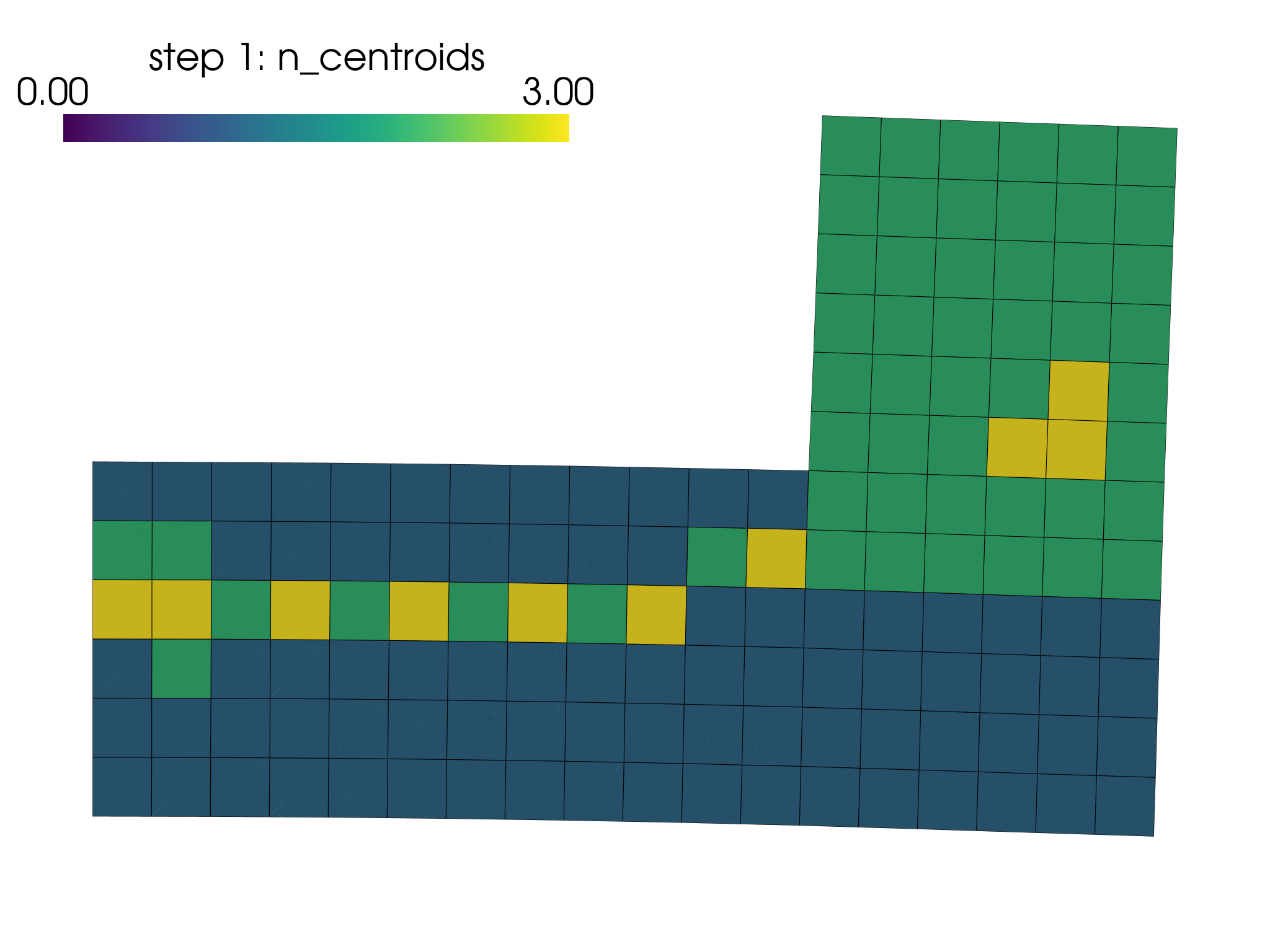}\hfil
        \includegraphics[width=0.49\linewidth]{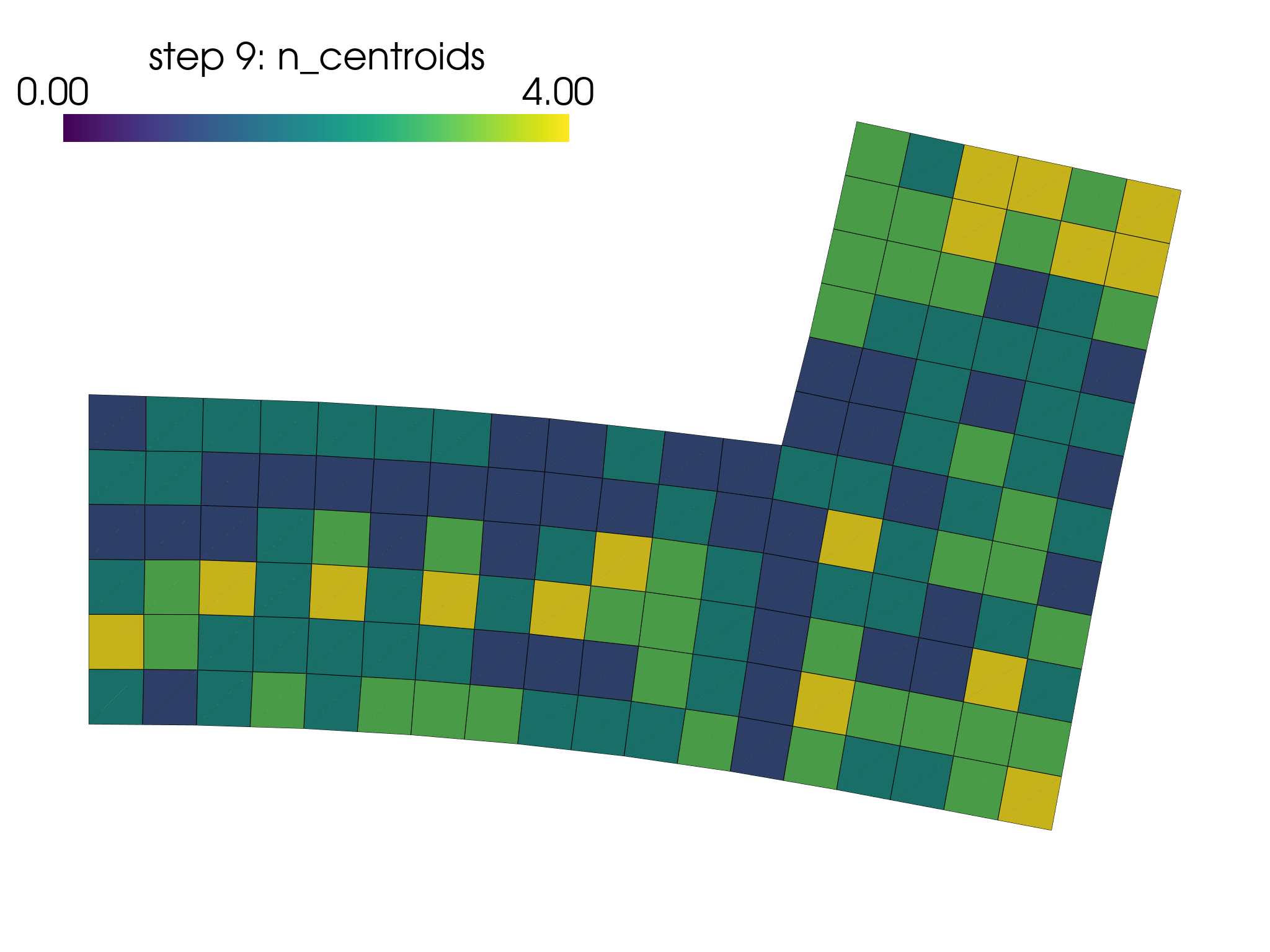}\\
    \includegraphics[width=0.49\linewidth]{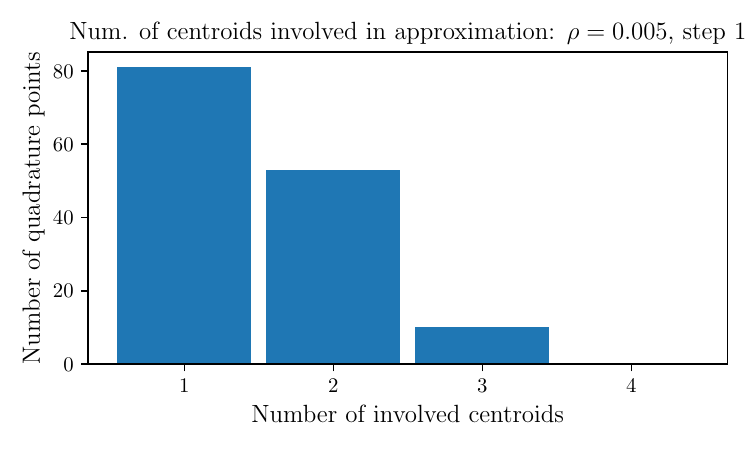}\hfil
        \includegraphics[width=0.49\linewidth]{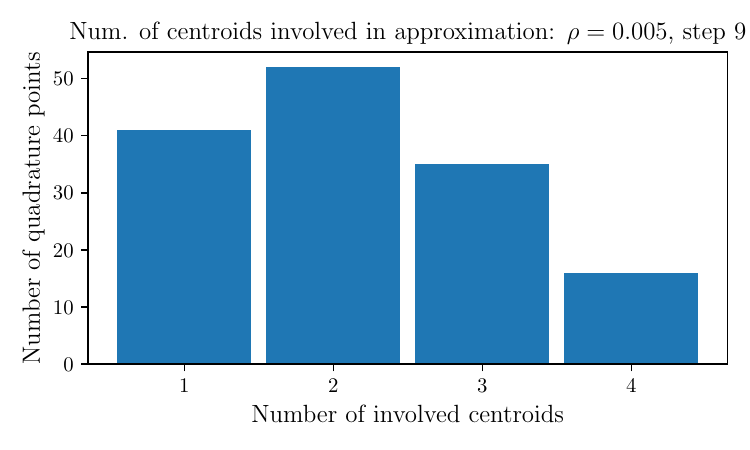}\\
    \caption{The number of centroids approximating the coefficients for $\rho =
        0.005$ in the final iterations of time steps 1 (left) and 9 (right).}
    \label{fig:num_centroids}
\end{figure}

The homogenized coefficients at the quadrature points of the
macroscopic domain are approximated using the values and sensitivities evaluated
at one or more centroids. If the deformation associated with a quadrature point
$\hx^i$ belongs to a single centroid $O^l$, the index set $\ihset{i}$ reduces to
${l}$, the corresponding weight $w^{il}$ equals 1, and \eq{eq:cf_approx}
yields an extrapolation from $\acirc\Ahom^l$, $\acirc\Shom^l$ to $\Ahom(\hx^i)$,
$\Shom(\hx^i)$. When the deformation lies within multiple centroids, the
homogenized coefficients are obtained by interpolation according to
\eq{eq:cf_approx}, using the weights defined in \eq{eq:weights}.
%
The number of centroids approximating the coefficients for $\rho = 0.005$ in the
final iterations of time steps 1 and 9 is shown in the upper part of
Figure~\ref{fig:num_centroids}. 
{To display the which cluster is employed to approximate the homogenized coefficients, the finite elements are
subdivided into four sub-parts each comprising one integration point, according
to the use of the four-point integration scheme.}
%
The lower part of Fig.~\ref{fig:num_centroids} shows the number of macroscopic
quadrature points for which the coefficients are approximated using values from
one, two, three, or four centroids. At load step~1, the coefficients at
approximately 55\% of the quadrature points are obtained by extrapolation from a
single centroid, while those at the remaining points are determined by
interpolation between two or three centroids.

\subsection{Approximation error}\label{sec:num_error}

To test the accuracy of the approximated homogenized coefficients, the relative
and cumulative errors are defined as
\begin{equation}\label{eq:num_err}
    \begin{aligned}
        Err^{rel}_{\mathcal{S}_{ij}^{k}}(\hat{x})
          & = \frac{\vert \mathcal{S}_{ij}^{k, \acron}(\hat{x}) - \mathcal{S}_{ij}^{k, \dtssx}(\hat{x})\vert}{\vert \Shom^{k, \dtssx}(\hat{x}) \vert},\\
        Err^{cum}_{\mathcal{S}_{ij}}(\hat{x})
          & = \sum\limits_k^{N_{iter}} Err^{rel}_{\mathcal{S}_{ij}^{k}}(\hat{x})
    \end{aligned}
\end{equation}
where $\Shom^{\acron}(\hat{x})$ is the approximated coefficient obtained from
\eq{eq:cf_approx}, and $\Shom^{\dtssx}(\hat{x})$ is the direct numerical solution
of the microscopic problem, both evaluated for the same macroscopic deformation
$\Fb(\hat{x})$. The superscript ${}^{k}$ denotes the $k$-th macroscopic
iteration, and $\hat{x}$ is the position of the quadrature point. Analogous
definitions apply to the homogenized elastic coefficient $\Ahom$.

\begin{figure}
    \centering
    \includegraphics[width=0.49\linewidth]{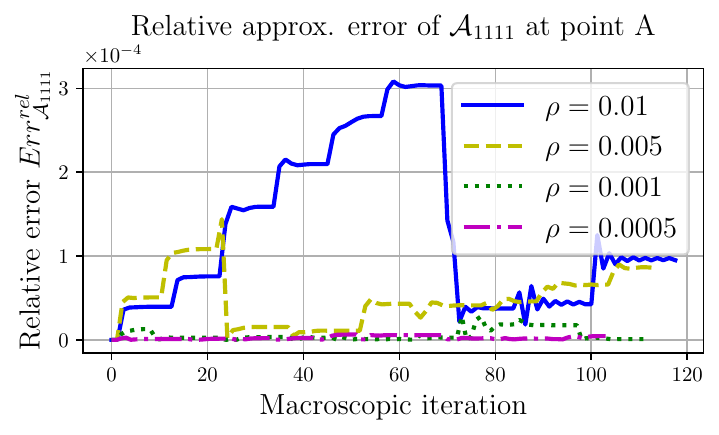}\hfil
        \includegraphics[width=0.49\linewidth]{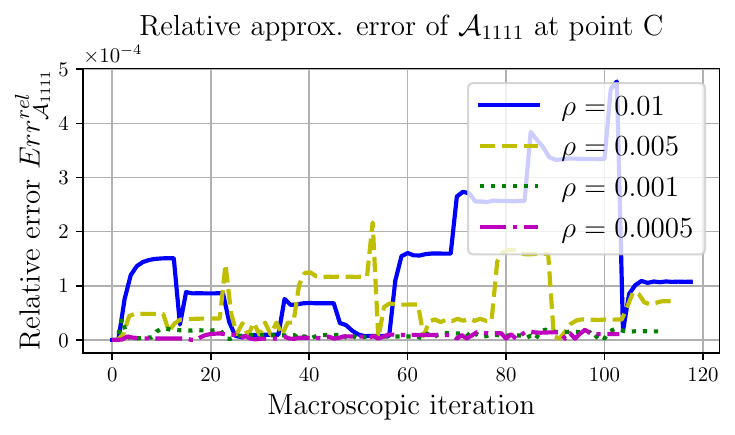}\\
    \includegraphics[width=0.49\linewidth]{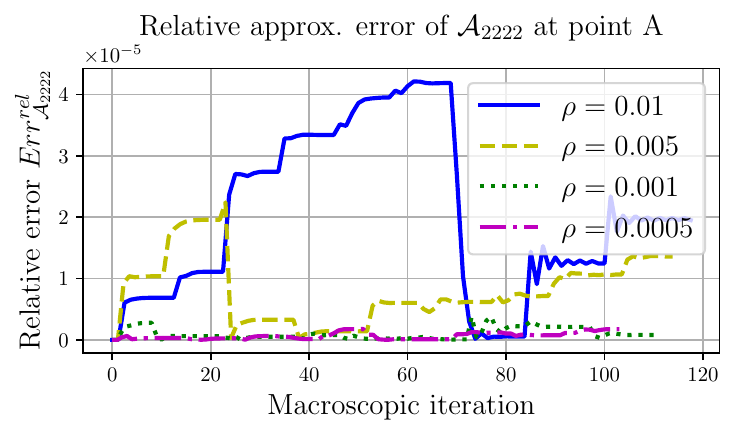}\hfil
        \includegraphics[width=0.49\linewidth]{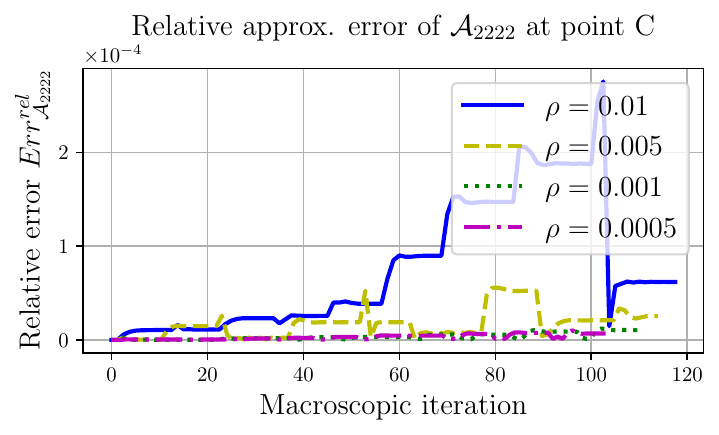}\\
    \includegraphics[width=0.49\linewidth]{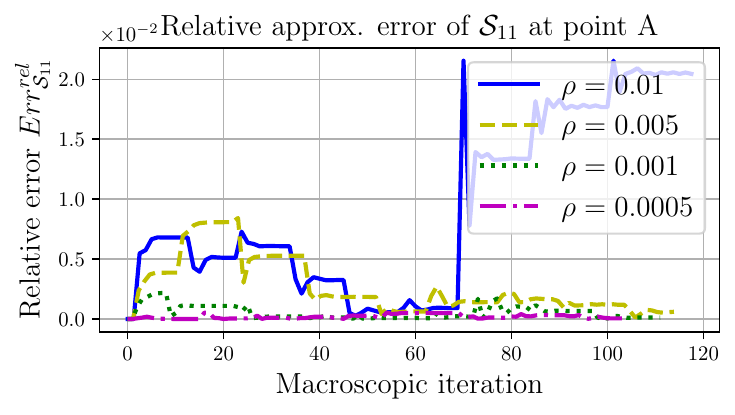}\hfil
        \includegraphics[width=0.49\linewidth]{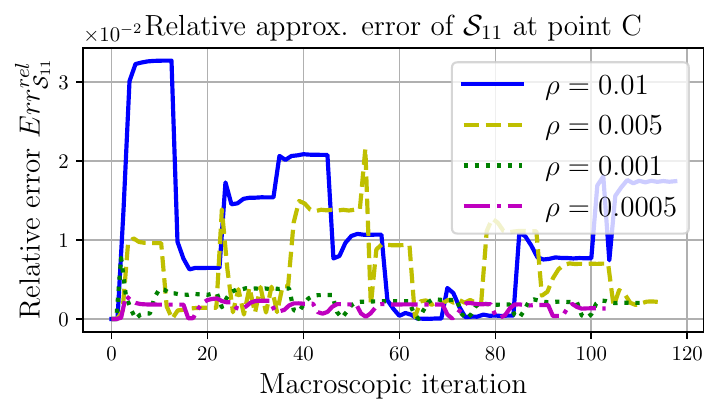}\\
    \includegraphics[width=0.49\linewidth]{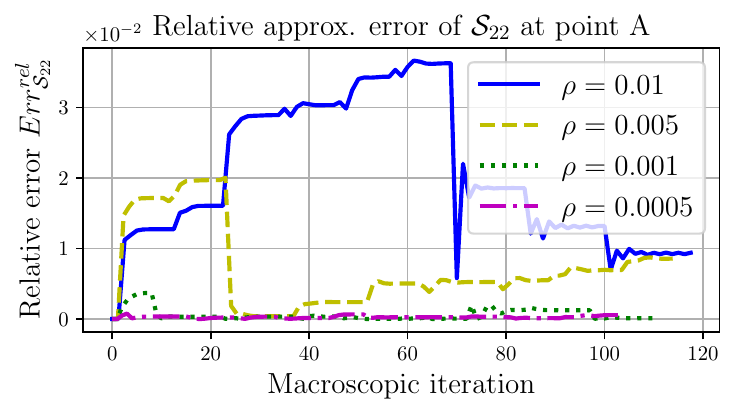}\hfil
        \includegraphics[width=0.49\linewidth]{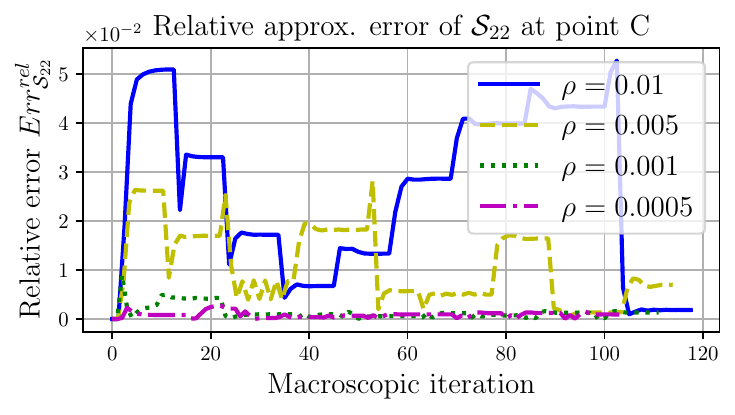}
    \caption{Relative errors of the homogenized coefficients at two quadrature
             points A (left) and C (right) calculated for the centroid size $\rho =
             0.01, 0.005, 0.001, 0.0005$.}
    \label{fig:num-rerr}
\end{figure}

\begin{figure}
    \centering
    \includegraphics[width=0.49\linewidth]{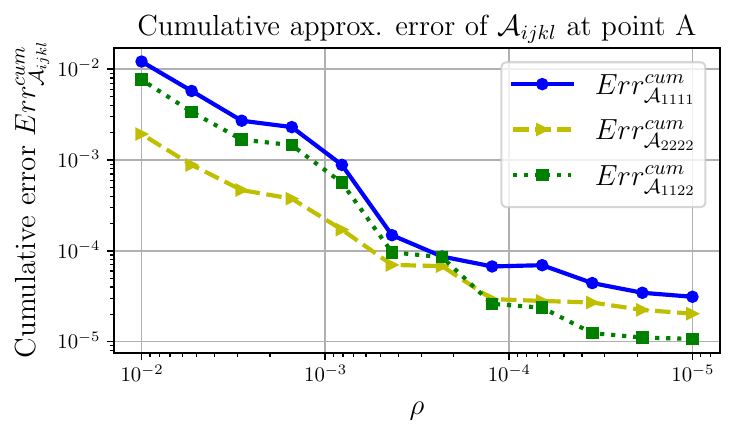}\hfil
        \includegraphics[width=0.49\linewidth]{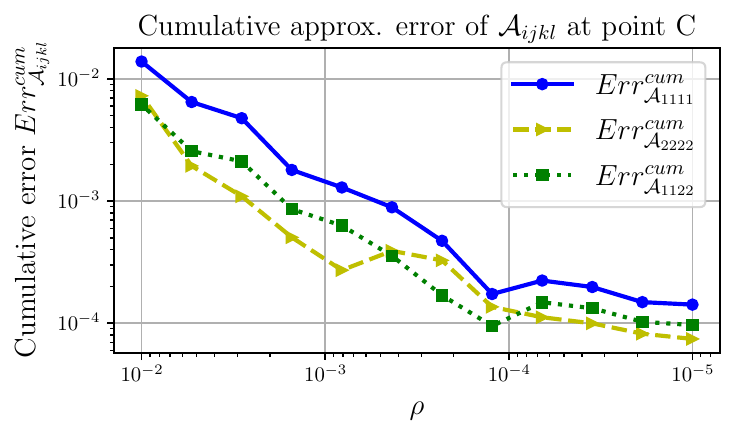}\\
    \includegraphics[width=0.49\linewidth]{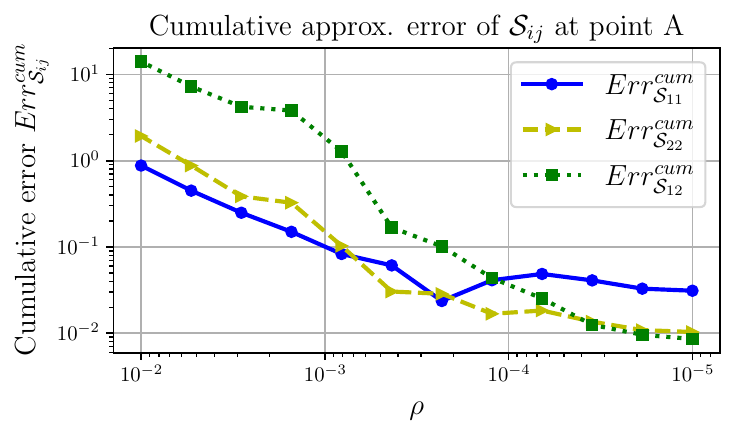}\hfil
        \includegraphics[width=0.49\linewidth]{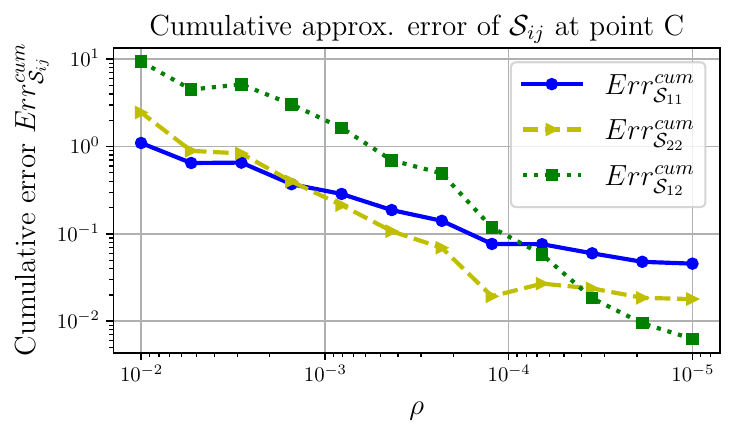}\\
    \caption{Cumulative error of coefficients $\Ahom$, $\Shom$ as a function of $\rho$.}
    \label{fig:num-cerr}
\end{figure}

Figure~\ref{fig:num-rerr} shows the evolution of the relative errors of
$\mathcal{A}_{1111}$, $\mathcal{A}_{2222}$, $\mathcal{S}_{11}$, and
$\mathcal{S}_{22}$ at two selected quadrature points A and C (see
Fig.~\ref{fig:num_deform_micro_recovered}), for the centroid size $\rho = 0.01,
0.005, 0.001, 0.0005$. 
The results indicate that the approximation error decreases as the centroid size
$\rho$ becomes smaller. For $\rho \approx 10^{-3}$, the relative error of
$\Shom$ is on the order of $10^{-3}$ and error of $\Ahom$ on the order of
$10^{-5}$. Insufficient accuracy for $\rho = 0.01$
leads to convergence issues at the macroscopic level, manifested as oscillations of the
approximated coefficients encountered between consecutive iterations.
These oscillations require the simulation to
be terminated after reaching the maximum number of macroscopic iterations.
During such oscillatory behavior, no new centroids are usually created, and the
approximation of the coefficients is based solely on the existing centroids. The
convergence of the approximated coefficients toward the directly computed values
is more clearly illustrated in Figure~\ref{fig:num-cerr}, which shows the
cumulative error as a function of $\rho$.

\subsection{Time complexity of reduced calculations}\label{sec:num_comp_time}

The practical applicability of the proposed approximation method is determined
by its accuracy, which was addressed in the previous section, as well as the
computational efficiency it offers. In the proposed algorithm, the simulation
time is influenced by the parameter $\rho$, which controls the approximation
error, and by the finite element mesh resolution at both the macroscopic and
microscopic scales. Since the number of macroscopic quadrature points is held
constant for all simulations in this paper, our focus is restricted to the
number of FE nodes and elements of the microscopic reference cell. To
evaluate the computational complexity of the method, a series of FE
microscopic meshes with identical geometry but varying numbers of elements
and nodes were generated, see 
{Table~\ref{tab:Y-cells}.}

\begin{table}
    \begin{tabular}{|l|c|c|c|}
        \hline
        Micro-mesh ID & Num. elements & Num. nodes & \\
        \hline\hline
        \#1 & 4166 & 2160 & \includegraphics[trim=0 0 0 -5,width=0.25\linewidth]{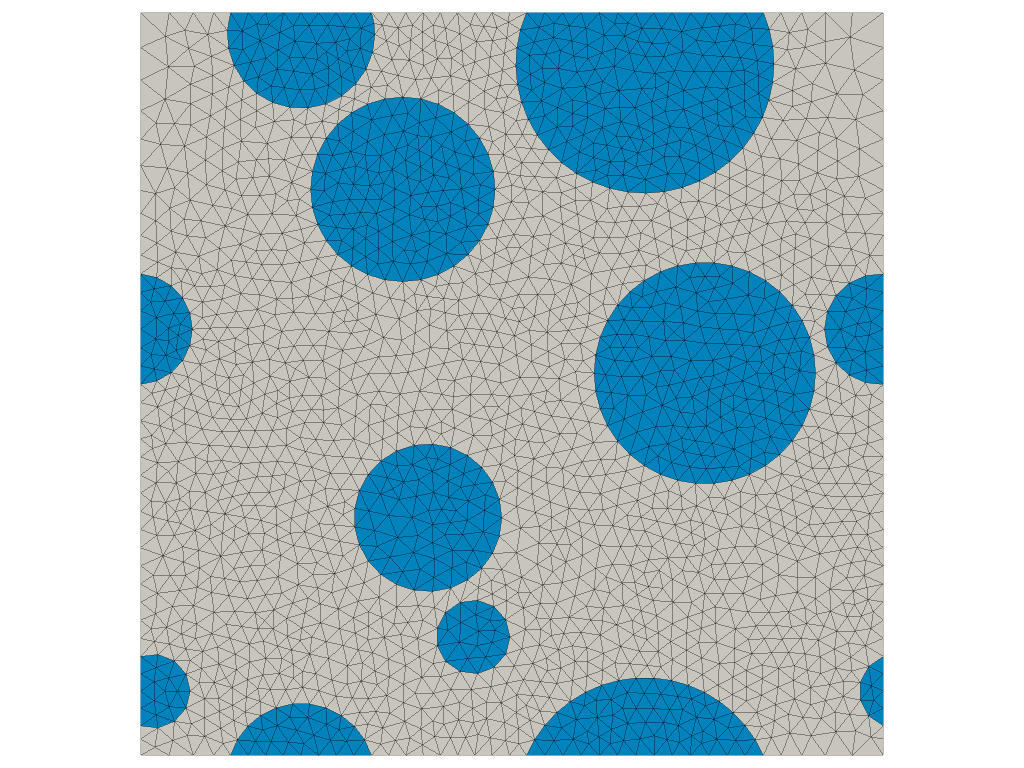}\\
        \hline
        \#2 & 15024 & 7961 & \\
        \hline
        \#3 & 24430 & 12399 & \\ 
        \hline
        \#4 & 42752 & 21620 & \\
        \hline
        \#5 & 57182 & 28872 & \\ 
        \hline
    \end{tabular}
    \caption{Periodic reference cells with varying mesh density.}\label{tab:Y-cells}
\end{table}

\begin{figure}
    \centering
    \includegraphics[width=0.49\linewidth]{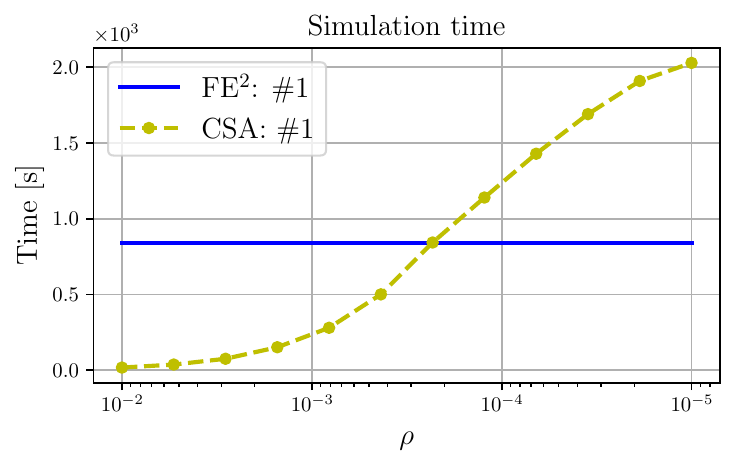}\hfil
        \includegraphics[width=0.49\linewidth]{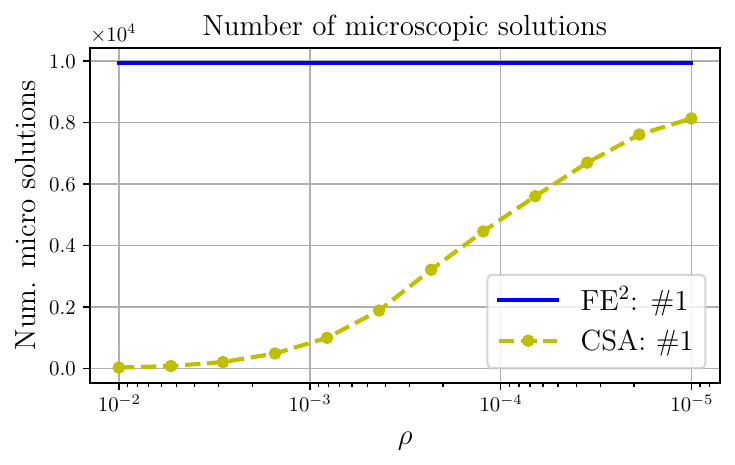}\\
    \caption{Solution time and number of solved microproblems for micro-mesh \#1 as functions of $\rho$.}
    \label{fig:num-complexity1}
\end{figure}

\begin{figure}
    \centering
    \includegraphics[width=0.49\linewidth]{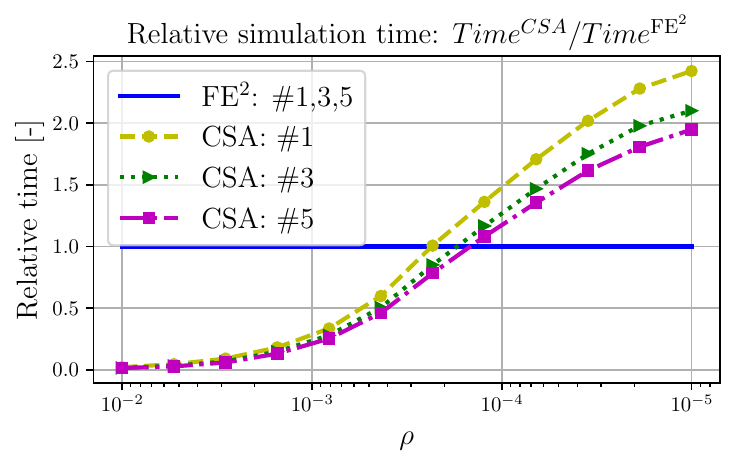}\hfil
     \includegraphics[width=0.49\linewidth]{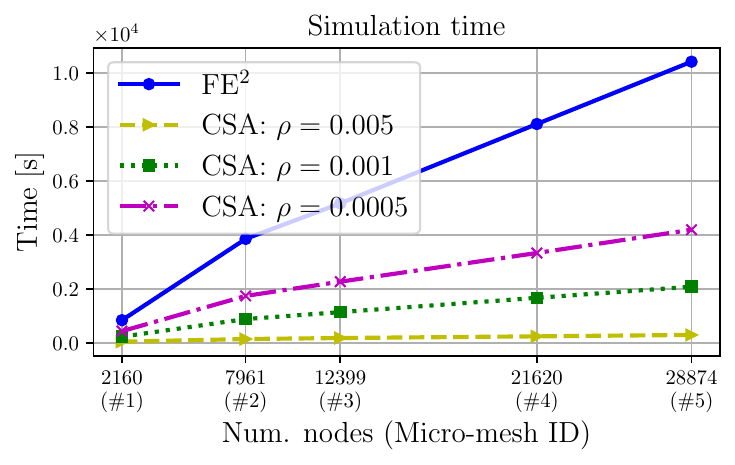}\\
    \caption{Left: relative simulation time as a function of $\rho$ for various microscopic FE mesh size;
             Right: simulation time as a function of microscopic FE mesh size.}
    \label{fig:num-complexity2}
\end{figure}

Figure~\ref{fig:num-complexity1} presents the solution time and the number of
solved microproblems for micro-mesh \#1 as functions of the centroid size
$\rho$, compared to the direct numerical simulation. As seen in
Fig~\ref{fig:num-complexity1} (left), the solution time exceeds that of the
direct simulation as $\rho \rightarrow 10^{-4}$ due to the additional
computational overhead associated with centroid generation, coefficient
approximation, and the evaluation of coefficient sensitivities. With the increasing
number of nodes in the microscopic FE mesh, the cost of solving the microscopic
problems becomes dominant relative to the overhead, resulting in more
significant computational savings provided by the proposed method. This trend is
evident in Fig.~\ref{fig:num-complexity2} (left), where the relative solution
time is defined as the ratio between the computational time of the reduced
model and the direct simulation. The absolute simulation time as a function of
the microscopic mesh size are presented in Fig.~\ref{fig:num-complexity2}
(right).

When the local subproblems \eq{eq:micro} are solved using a direct sparse linear
solver, 
{the total computational complexity for the coefficient evaluation is of
order $o(N_{mac} \cdot N_\# \cdot b^2)$ \cite{Golub_2013} for both the \dtss and \acron simulations, where
$N_\#$ denotes the number of microscopic degrees of freedom, $N_{mac}$ is the
number of local subproblems, and $b$ is the bandwidth of the sparse matrix
arising from the FE discretization of the microscopic problem.}
For the \dtss, $N_{mac}$ equals the number of macroscopic quadrature points
multiplied by the number of macroscopic iterations, whereas for \acron,
$N_{mac}$ represents the number of centroids generated during the simulation.
Since $N_{mac}$ remains constant for a given algorithm and $b$ does not change
with varying mesh density, the overall time complexity exhibits linear scaling
behavior, as confirmed in Fig.~\ref{fig:num-complexity2} (right).

\begin{figure}
    \centering
    \includegraphics[width=0.49\linewidth]{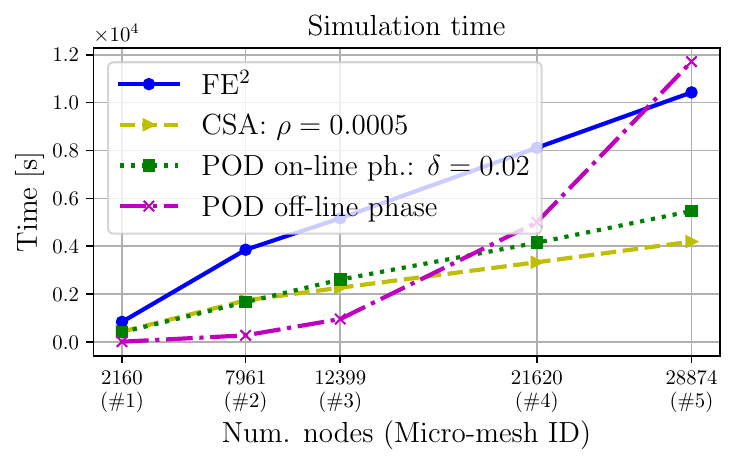}
    \includegraphics[width=0.49\linewidth]{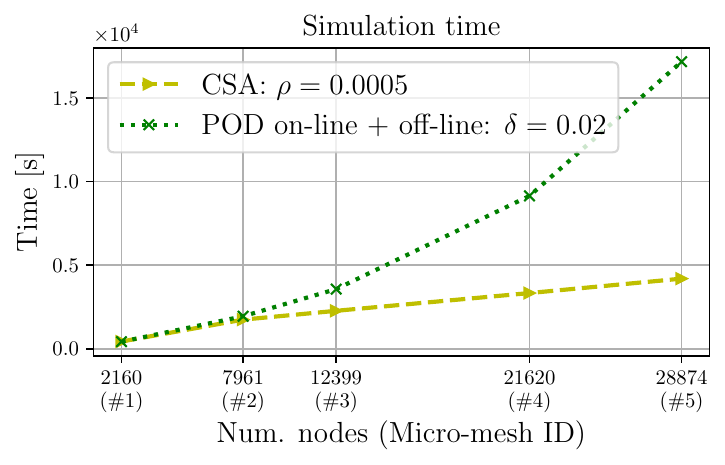}
    \caption{Comparison of the simulation times of the proposed reduction algorithm (for $\rho=0.0005$),
             direct simulation, and the POD method ($\delta=0.02$).}
    \label{fig:num-complexity3}
\end{figure}

\begin{figure}
    \centering
    \includegraphics[width=0.52\linewidth]{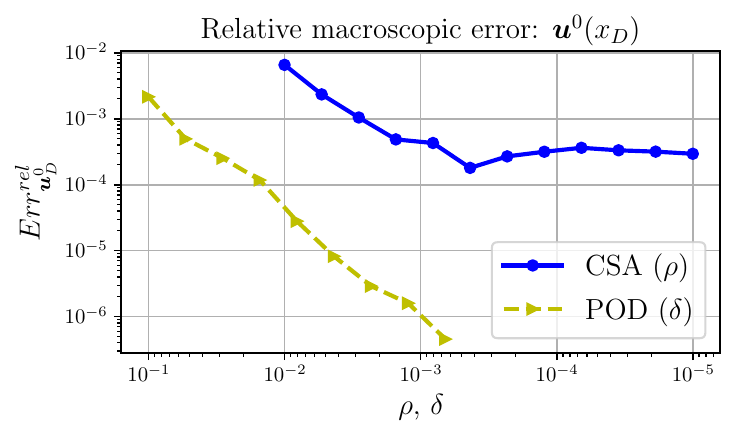}
        %
    \caption{Relative macroscopic error at point D for the CSA and POD algorithm.}\label{fig:num-complexity4}
\end{figure}

The efficiency of the CSA algorithm is compared against the direct simulation and
the proper orthogonal decomposition method, see~\cite{Yvonnet_2007}, in
Fig.~\ref{fig:num-complexity3}. In the figure, ``POD'' refers to the simulation time
using a pre-computed orthogonal basis and ``POD -- basis calc.'' denotes the
time required for computing that basis. The POD simulations are carried out with the error
parameter $\delta = 0.02$, and the \acron simulations use the centroid size
$\rho=0.0005$. These settings yield a comparable order of the macroscopic error
measured by $Err^{rel}_{\ub^0_D}$, as illustrated in Fig.~\ref{fig:num-complexity4}, where
\begin{equation}\label{eq:num_err2}
    Err^{rel}_{\ub^0_D}(x_D)
         = \frac{\vert \ub^{0, \acron}(x_D) - \ub^{0, \dtssx}(x_D)\vert}{\vert \ub^{0, \dtssx}(x_D) \vert}.
\end{equation}

The complexity of the POD ``on-line'' phase is primarily given by the solution
of the reduced (dense) system \eq{eq:pod-reduced} of dimension $M \times M$, which is
independent of $N_\#$, and by the matrix multiplications required for order
reduction. These operations show linear time complexity scaling, as demonstrated
in Fig.~\ref{fig:num-complexity3} by the green dotted line. In contrast, the
``off-line'' phase involves solving the eigenvalue problem \eq{eq:pod-evp} for a
dense symmetric matrix of dimension $N_\# \times N_\#$ with complexity of order
$o(N_\#^3)$, as shown by the cyan dash-dotted line in
Fig.~\ref{fig:num-complexity3}. This cubic growth can considerably limit
the applicability of the POD for very large microscopic FE meshes.
%

All numerical simulations were implemented in the general-purpose FE package
\SfePy \cite{Cimrman_Lukes_Rohan_2019} without any significant optimization.
Using dedicated, highly optimized implementations of individual algorithms would
likely result in shorter simulation times.


\subsection{3D example}\label{sec:num-3D}

\begin{figure}
    \centering
    \includegraphics[width=0.33\linewidth]{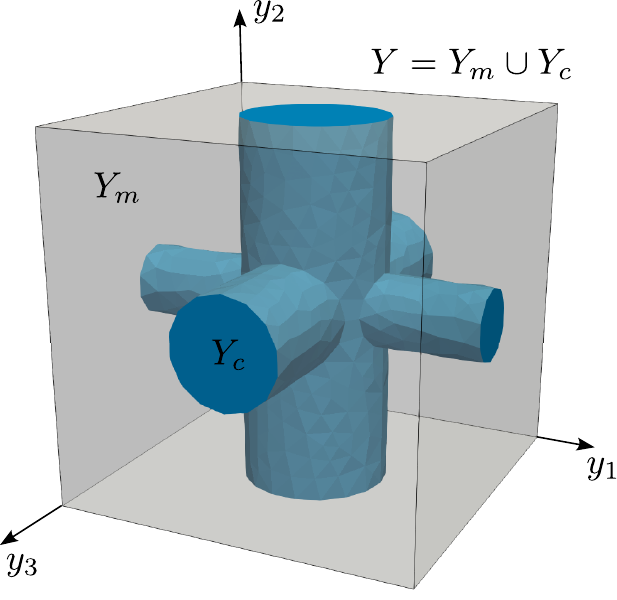}\hfil
    \includegraphics[width=0.63\linewidth]{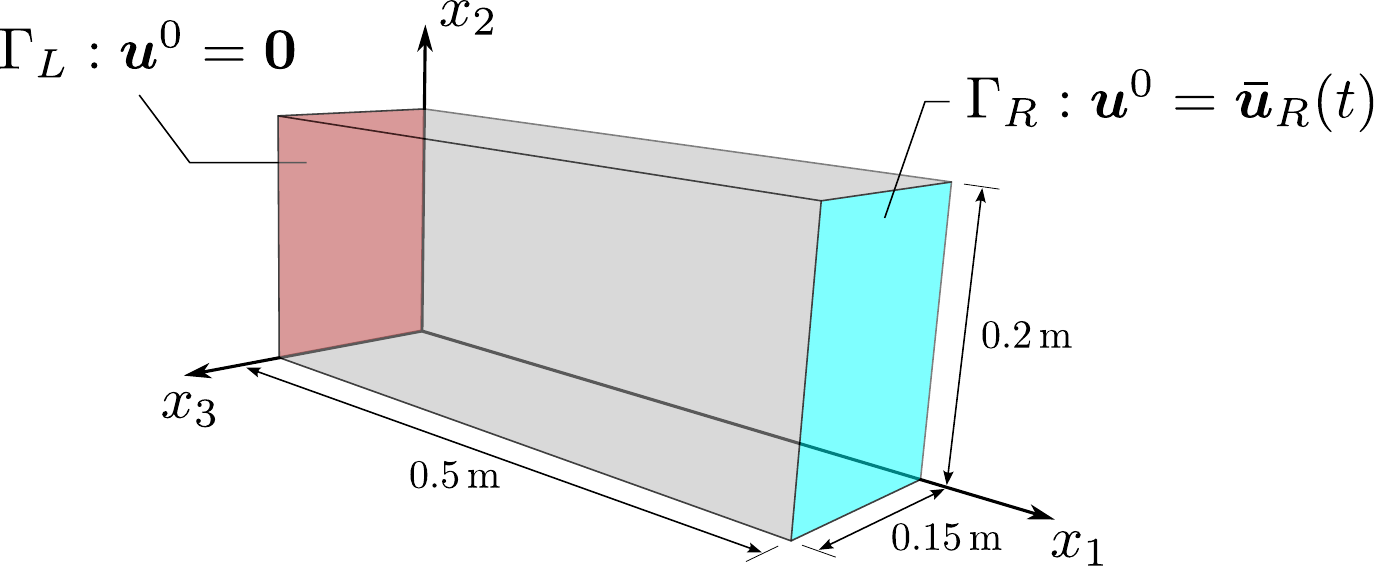}\\
    \caption{Geometries of the reference unit cell (left) and the macroscopic sample (right).}
    \label{fig:porohyperela-d2s-example-geom3D}
\end{figure}

The implementation of the CSA algorithm in \SfePy also enables solution of
three-dimensional problems, as demonstrated in this section. At the macroscopic
level, a block-shaped specimen with dimensions $0.5 \times 0.2 \times 0.15$\,m is
considered, see Fig.~\ref{fig:porohyperela-d2s-example-geom3D} (right). One {of} its
faces is fixed, $\ub^0 = \boldsymbol{0}$ on boundary $\Gamma_L$, while the
opposite face is subjected to prescribed displacements, $\ub^0 = \bar\ub_{R}(t)$ on $\Gamma_R$, where
\begin{equation}\label{eq:num_displ_fun3}
    \begin{aligned}
    \bar\ub_{R}(t_k) = \left(
        \begin{bmatrix}
            1 & 0 & 0\\
            0 & \cos(\varphi_k) & \sin(\varphi_k)\\
            0 & -\sin(\varphi_k) & \cos(\varphi_k)
        \end{bmatrix}
        - \Ib\right) \boldsymbol{x},\\
        \qquad \mbox{for}\quad \varphi_k = \frac{35}{180}\pi\,r_k, \quad r_k = \frac{k}{N_t - 1}, 
        \quad \mbox{and}\quad t_k = r_k \cdot 1\,\mbox{s}, \quad k = 0, \dots, N_t -1.
    \end{aligned}
\end{equation}
The finite element mesh of the macroscopic sample consists of 200 hexahedral Q1 elements,
for which a eight-point numerical integration rule is applied.

The reference unit cell representing the periodic microstructure, shown in
Fig.~\ref{fig:porohyperela-d2s-example-geom3D} (left), is composed of two
continuous subdomains with distinct material properties identical to those used
in the two-dimensional simulation, see Table~\ref{tab:materials}.

\begin{figure}
    \centering
    \includegraphics[width=0.99\figscale]{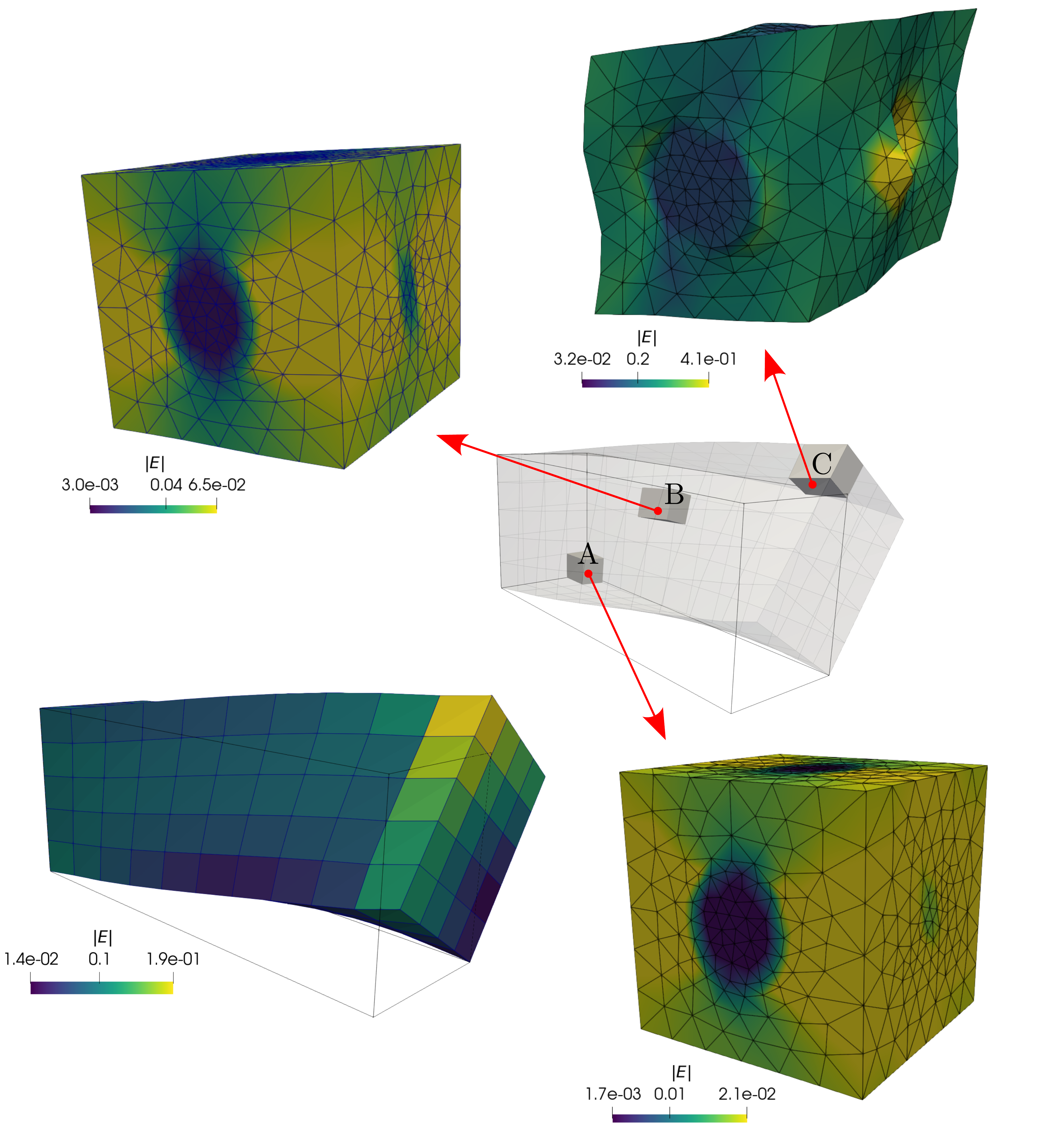}
    \caption{Deformed macroscopic sample and deformed microstructures at selected points.}
    \label{fig:macro-micro-3D-def}
\end{figure}

\begin{figure}
    \centering
    \includegraphics[width=0.49\linewidth]{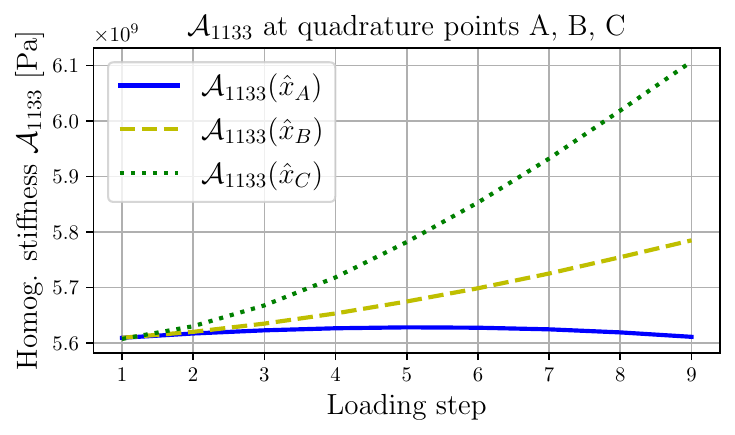}\hfill
    \includegraphics[width=0.49\linewidth]{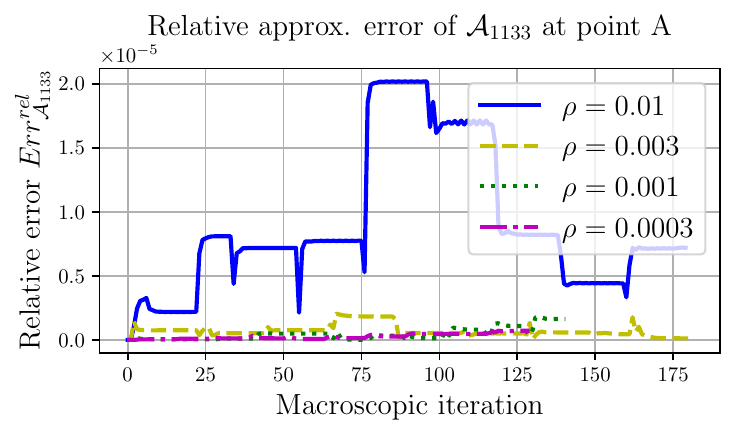}\\
    \hfill
    \includegraphics[width=0.49\linewidth]{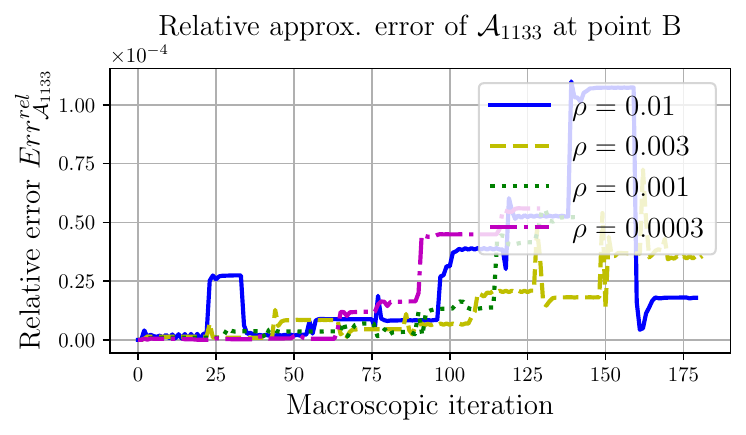}\\
    \caption{Evolution of the homogenized coefficient $\mathcal{A}_{1133}$
    (left) and its relative approximation errors at macroscopic points A and B (right).}
    \label{fig:num-rerr-3D-A}

\end{figure}

\begin{figure}
    \centering
    \includegraphics[width=0.49\linewidth]{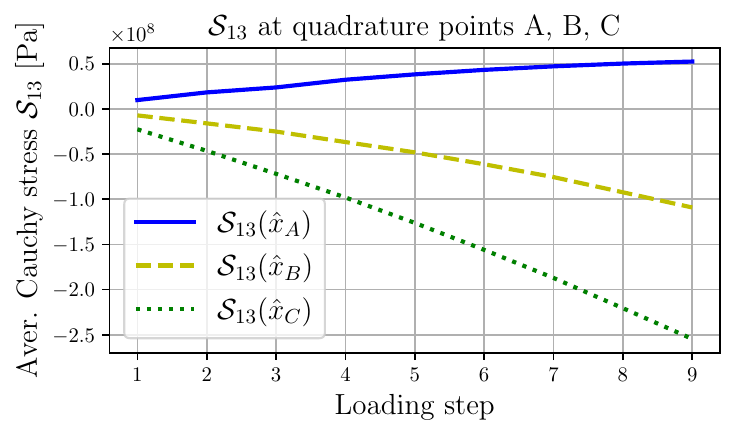}\hfill
    \includegraphics[width=0.49\linewidth]{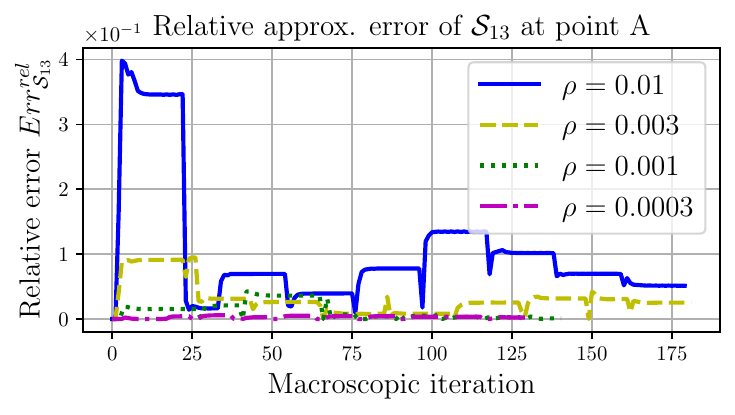}\\
    \hfill
    \includegraphics[width=0.49\linewidth]{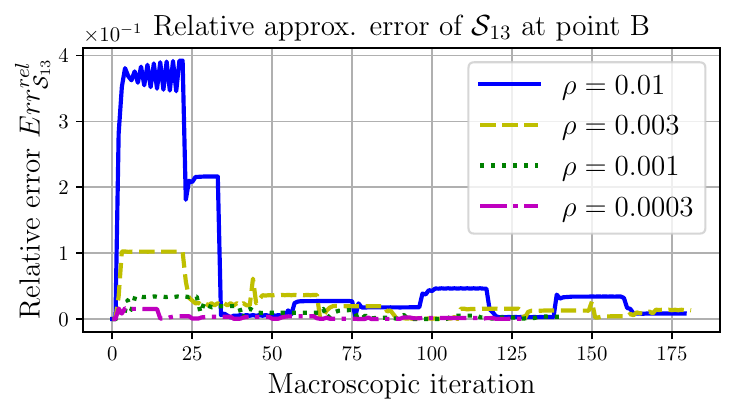}\\
    \caption{Evolution of the homogenized coefficient $\mathcal{S}_{13}$
    (left) and its relative approximation errors at macroscopic points A and B (right).}
    \label{fig:num-rerr-3D-S}
\end{figure}

The deformed macroscopic configuration and the reconstructed microscopic
structures at selected points are displayed in Fig.~\ref{fig:macro-micro-3D-def}
for the final loading step. Figures~\ref{fig:num-rerr-3D-A}
and~\ref{fig:num-rerr-3D-S} show the approximation errors of $\mathcal{A}_{1133}$
and $\mathcal{S}_{13}$ at quadrature points A and B for different values of
$\rho$. As the centroid size $\rho$ decreases, the approximation errors
become smaller and the number of macroscopic iterations required to achieve
equilibrium decreases, as evidenced by the shorter curves for $\rho =
0.001$ and $\rho = 0.0003$.

\begin{figure}
    \centering
    \includegraphics[width=0.49\linewidth]{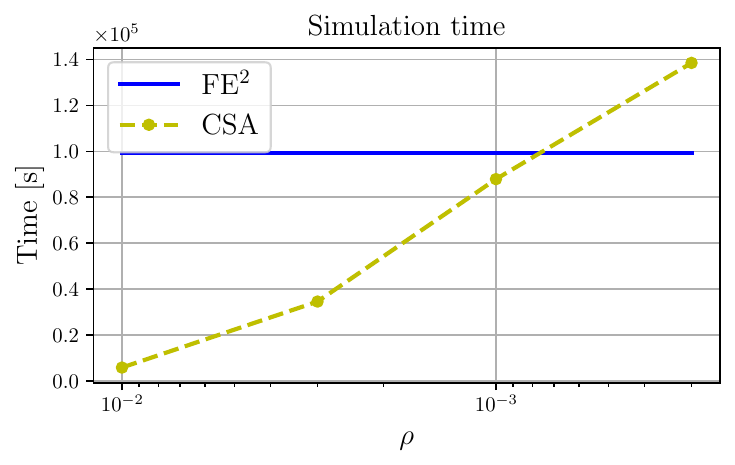}\hfill
    \includegraphics[width=0.49\linewidth]{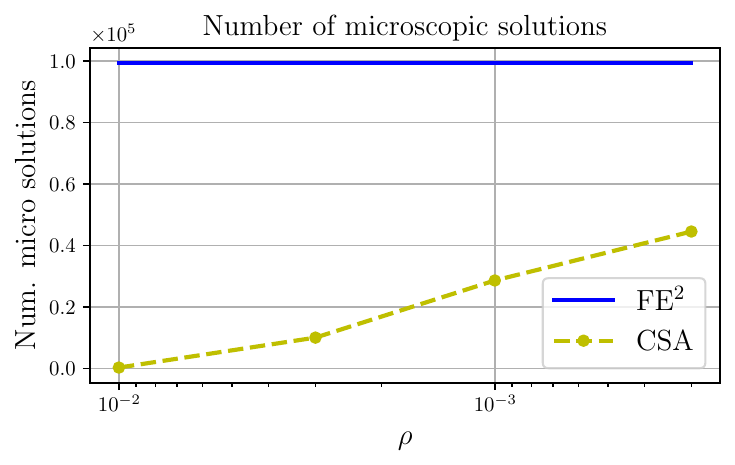}\\
    \caption{Solution time (left) and number of solved microproblems (right) as functions of $\rho$.}
    \label{fig:num-complexity-3D}
\end{figure}

Figure~\ref{fig:num-complexity-3D} compares the computational times (left) and the
numbers of solved microscopic problems (right) for the \acron simulations and the direct
two-scale (\dtss) computations. The results confirm that, even for the
three-dimensional case, an appropriate choice of the parameter $\rho$ allows
a significant reduction of computational cost while maintaining an acceptable
level of approximation accuracy.



%% file: part_conclusion.tex
\section{Conclusion}\label{sec:conclusion}

The paper contributes to the field of Model Order Reduction for multiscale
computational nonlinear analysis of deforming solids. We have introduced the CSA
method for the efficient
approximation of homogenized coefficients in hyperelastic media undergoing large
deformations. The proposed algorithm relies two ingredients: a) a linear approximation of deformation-dependent homogenized coefficients constructed using representative defortmed microstructures and their sensitivities \wrt macroscopic strain variations; b) the representative microstructures created ``on-flight'' during the simulation; for this the clustering of macroscopic
deformations based on their distances in metrics of the deformation space is employed.
This approach significantly decreases the number of
microscopic problems that need to be solved and, hence, the computational cost.
Each cluster of macroscopic deformations is represented by a single mean
deformation, for which the microscopic problem is solved and the corresponding
coefficients, together with their sensitivities with respect to the macroscopic
deformation, are computed. These values are then used to approximate the
homogenized coefficients for all required macroscopic deformation states. The
clustering procedure is dynamic and is repeatedly applied during the iterative
solution of the nonlinear problem.

The computational savings of the CSA algorithm depend on the required
approximation accuracy. With increasing accuracy the algorithm efficiency diminishes, as compared 
the full numerical simulation (FE$^2$ approach), primarily because of the
overhead associated with evaluating coefficient sensitivities. There is a
specific accuracy limit beyond which the time complexity of the CSA method
exceeds that of direct two-scale simulation. However, when a fast solution is
required and a certain approximation error is acceptable, the CSA method can
speed up simulations by several factors. For the numerical example presented in
Section~\ref{sec:num_comp_time}, allowing a relative macroscopic error of
approximately $10^{-4}$ results in a speedup factor of 2.5 compared to the
FE$^2$ method, which is comparable to the on-line phase of the POD-based
simulation.

The efficiency of the CSA approach was also studied with respect to the number
of degrees of freedom of the discretized microscopic problem. The computational
complexity of the CSA method scales nearly linearly with increasing microscopic
resolution, which is a significant advantage over classical proper orthogonal
decomposition methods, for which the cost of basis computation grows rapidly.
Consequently, the proposed method is suitable for fast approximate simulations
of problems with complex microscopic geometries.

The presented CSA method has been tested in the two-scale computational analysis
of hyperelastic periodically heterogeneous solids, however, it can be extended
with some limitations to account for: a) rate-independent inelastic solids with
the constitutive laws involving macroscopic internal variables, b) functionally
graded media with ``gently varying'' microstructures with macroscopic spatial
position, hence also to solve two-scale optimization problems, where the
required sensitivity analysis of the homogenized coefficients (properties) is
naturally treated along with the sensitivity \wrt the macroscopic deformation.
These issues will be addressed in our further research. Extensions of the CSA
method for media characterized by dissipation can be considered on a ``heuristic
basis'', in analogy with other works that employ macroscopic clustering.
However, in our view, the memory effects featuring the time evolution of each
local microstructure require tracking all the history, which, in effect,
enormously expands the micro-configuration parameterization.

The final remark concerns the potential improvements of the sensitivity-based
approximation between neighboring centroids. In this study, we reported a simple
method that uses the linear approximation within each centroid for all
macroscopic deformation states (MDS) belonging exclusively to one centroid; for
MDS points lying in the overlap of several centroids, a blending interpolation
is constructed using the weighted linear approximation provided by each of the
involved centroids. Further promising improvements have been explored and will
be reported in a forthcoming publication devoted also to the above mentioned
extensions.

%% file: statements.tex
\section*{CRediT authorship contribution statement}
\noindent
{\bf Vladim\'\i{}r Luke\v{s}}: Writing, Visualization, Software, Methodology, Investigation, Conceptualization;\\
{\bf Eduard Rohan}: Writing, Formal analysis, Methodology, Investigation, Conceptualization, Funding acquisition, Project administration.

\section*{Data Availability}
Code and data available at \href{https://data.mendeley.com/datasets/gvp9kj2nvy/1}{https://data.mendeley.com/datasets/gvp9kj2nvy/1}.

\section*{Declaration of competing interest}
The authors declare that they have no known competing financial interests or
personal relationships that could have appeared to influence the work reported
in this paper.

\section*{Acknowledgements}


The research has been supported in a part by the grant project GA~2306220S,
but mainly by project GA~2412291S of the Czech Science Foundation.

%% file: appendix_sa.tex
We use the shape sensitivity technique and the material derivative approach (see
e.g. \cite{Haslinger_1996}) to obtain the sensitivity of homogenized coefficient
\eq{eq:homog_coefs} to the configuration transformation corresponding to the
microscopic deformation.

Consider a perturbation of material points which is given in terms of a
differentiable and $Y$-periodic vector field $\Vcalbf(y): y\mapsto y' \in
\RR^d$, $y \in Y \subset \RR^d$. As explained below, construction of the
convection velocity $\Vcalbf$ is based on the microscopic displacement and
deformation fields recovery related to the corrector result of the
homogenization. Using $\Vcalbf$ defined in $Y$, perturbed configuration $Z(\tau)
= \vphibf(Y,\tau)$ can be established through mapping $\vphibf(y,\tau)$ defined
for $y \in Y$, where $\tau\in \RR$ is a ``time-like'' variable. The ``material
point'' position in $Y$ is parameterized by $\tau$, yielding the  perturbed
coordinates $z(y,\tau)$ defined, as follows
\begin{equation}\label{eq-shsa1}
  \begin{split}
    z_i(y,\tau) & = \vphi_i(y,\tau) := y_i + \tau\Vcal_i(y)\;, \;i = 1,\dots,d\;,\\
    y_i & = \vphi_i^{-1}(y,\tau) := z_i - \tau \tilde\Vcal_i(z)\;, \;i = 1,\dots,d\;,\\
    Z(\tau) & = \{z \subset \RR^d|\; z_i(y,\tau) = \vphi_i(y,\tau)\;, \;i = 1,\dots,d\;,\; y \in Y\}\;,\\
    \Fcal_{ij} & = \pd_j^y \vphi_i = \pd_j^y z_i(y,\cdot).
\end{split}
\end{equation}
The inverse mapping $\vphibf^{-1}(z,\tau):z\mapsto y$ is defined in terms of
$\tilde\Vcalbf(z(y,\tau))$, such that $\tilde\Vcalbf(z(y,0))=\Vcalbf(y)$. The
deformation gradient $\Fcal_{ij}$ is employed in the differentiation of
integrals involving gradients $\nabla_y = (\pd_i^y)$ of scalar, or vectorial
functions. Using the perturbed configuration and the associated derivatives
$\pd_i^z$, the following holds, where $\dot{\what{(~)}} = \pd_\tau(~)|_{\tau =
0}$,
\begin{equation}\label{eq-shsa2}
  \begin{split}
    \dot\Fcal_{ik} &= \pdiff{\Vcal_i}{y_k}\;,\\
    \dot{\what{(\Fcal)_{ik}^{-1}}} &= -\pdiff{\Vcal_i}{y_k}\;,\\
    \dot{\what{\det(\Fcalbf)}} &= \nabla_y\cdot\Vcalbf\;.
\end{split}
\end{equation}

Recalling that the perturbation defined using $\Vcalbf(y)$ is parametrized by
$\tau$, throughout the text below we shall use the notion of the following
derivatives: $\delta(\cdot)$ is the total ({\it material}) derivative,
$\delta_\tau(\cdot)$ is the partial ({\it shape}) derivative with respect to
$\tau$.  These derivatives are computed as the directional derivatives in the
direction of $\Vcalbf(y)$, $y \in Y$.

Let $u(y)$ be differentiable in $Y$ (whereby $u$ can represent a component of
vector field $\ub$), by $\tilde{u}(z,\tau)$ we denote the perturbed function,
an image of $u$ thus, defined in $Z$ given by \eq{eq-shsa1}$_3$. To compute the
sensitivity $\dlt \Acal_{klij}$,
the partial shape derivative of the strains is required, hence we need 
\begin{equation}\label{eq-shsa3}
  \begin{split}
\dlt_\tau (\nabla_y u)_i & = \ddt{}{\tau}\pd_i^z\tilde{u}(z,\tau)|_{\tau = 0} = \ddt{}{\tau}\left( \pdiff{u}{y_k}\pdiff{y_k}{z_i}\right)_{\tau = 0}\\ & = \ddt{}{\tau}\left( \pdiff{u}{y_k}\left(\delta_{ki} - \tau\pdiff{\tilde\Vcal_k(z)}{z_i}\right)\right)_{\tau = 0} = -\pdiff{\Vcal_k}{y_i}\pdiff{u}{y_k}\;,
\end{split}
\end{equation}
where \eq{eq-shsa2} was employed. Application of \eq{eq-shsa2} and \eq{eq-shsa3}
in the differentiation of the bilinear form \eq{eq:alin}, leads to
\eq{eq:sa_alin}, where $\dltsh\tilde\Aop$ is given by~\eq{eq:app2-daop}.

\paragraph{Sensitivity of the tangential elasticity tensor $\Acal_{klij}$} The
sensitivity formula \eq{eq:sa_elastic} is derived using the partial shape
derivative \eq{eq:sa_alin} and using the micro-problem \eq{eq:micro}, where
$\vb:=\omegabf^{ij} \in \HpbY$ can be substituted, so that
$\alin{\Xibf^{ij}}{\dlt\omegabf^{kl}} = 0$. This yields
\begin{equation}\label{eq-shsa4}
  \begin{split}
    \dlt \Acal_{klij} & = \dltsh \alin{\Xibf^{ij}}{\Xibf^{kl}} +  \alin{\Xibf^{ij}}{\dlt\omegabf^{kl} + \dltsh\Pibf^{kl}} + \alin{\dlt\omegabf^{ij} + \dltsh\Pibf^{ij}}{\Xibf^{kl}}\\
 & = \dltsh \alin{\Xibf^{ij}}{\Xibf^{kl}} +  \alin{\dltsh\Pibf^{ij}}{\Xibf^{kl}}
            + \alin{\Xibf^{kl}}{\dltsh\Pibf^{kl}}\;.
\end{split}
\end{equation}

\paragraph{Sensitivity of the effective stress $\Scal_{ij}$}

The sensitivity of $\Scal_{ij}$ is obtained using the partial ({\it shape})
derivative and, namely, making use of \eq{eq-shsa2}$_3$, so that 
\begin{equation}\label{eq-S19}
\begin{split}
  \delta \Scal_{ij} & = \delta_\tau \intY \tilde\sigmabf \dVy
= -\Shom \intY \nabla_y\cdot \Vcalbf\dVy + \intY\tilde\sigmabf \nabla_y\cdot \Vcalbf\dVy
            + \intY\dlt_\tau\tilde\sigmabf \dVy\;,
\end{split}
\end{equation}
which yields \eq{eq:sa_stress}, where $\dlt_\tau\tilde\sigmabf$ is given by~\eq{eq:app2-cs}.


The partial derivatives of the Cauchy stress tensor $\dtau\tilde\sigmabf$ and
the associated tangent moduli $\dtau\tilde{\Aop}$ are derived on the basis of
the expressions for the Kirchhoff stress $\taubf = J \tilde\sigmabf$ reported
in~\cite{Crisfield_1991}. In particular, the derivative of the Kirchhoff stress
$\dtau\taubf$ is first evaluated,
\begin{equation}\label{eq:app2-ks-diff}
    \begin{aligned}
        \dtau\tau\ij = & 
        K(2J - 1)J\delta\ij \cdivY\Vcalbf
        -\frac{2}{3} \mu J^{-2/3}\left(b\ij - \frac{1}{3}b_{kk}\delta\ij\right) \cdivY\Vcalbf\\
        &+ \mu J^{-2/3}\left(\partial_k \Vcal_i b_{kj} + b_{il}\partial_l \Vcal_j - \frac{1}{3}(\partial_k \Vcal_l b_{kl} + b_{kl}\partial_l \Vcal_k)\delta\ij\right),
    \end{aligned}
\end{equation}
from which the derivative of the Cauchy stress tensor follows as
\begin{equation}\label{eq:app2-cs}
    \begin{aligned}
        \dtau\tilde\sigmabf = \frac{1}{J}\left(\dtau\taubf - \taubf\cdivY\Vcalbf\right).
    \end{aligned}
\end{equation}
The Truesdell rate of the Kirchhoff stress tensor, denoted by $\Dop^{\mathrm{TK}}$, is
related to the Jaumann rate $\Dop^{\mathrm{JK}}$ through
\begin{equation}\label{eq:app2-tk}
    D^{\mathrm{TK}}_{ijkl} = D^{\mathrm{JK}}_{ijkl}
     - \frac{1}{2}\left(
            \delta_{ik}\tau_{lj} +
            \tau_{il}\delta_{jk} +
            \delta_{il}\tau_{kj} +
            \tau_{ik}\delta_{jl}
        \right),
\end{equation}
%
where the expression for $\Dop^{\mathrm{JK}}$ is derived in~\cite{Crisfield_1991},
\begin{equation}\label{eq:app2-jk}
    \begin{aligned}      
        D^{\mathrm{JK}}_{ijkl} = &K(2J - 1)J \delta_{ij}\delta_{kl} 
            + \mu J^{-2/3}\frac{2}{9} b_{mm}\delta_{ij}\delta_{kl} \\
            &+ \mu J^{-2/3}\left(-\frac{2}{3}\left(b_{ij}\delta_{kl} + \delta_{ij}b_{kl}\right)
        + \frac{1}{2}\left(\delta_{ik}b_{lj} + b_{il}\delta_{jk} + \delta_{il}b_{kj} + b_{ik}\delta_{jl}\right)\right).
    \end{aligned}
\end{equation}
Differentiation of the Jaumann rate yields the corresponding derivative
$\dtau \Dop^{\mathrm{JK}}$, 
\begin{equation}\label{eq:app2-djk}
    \begin{aligned}      
        \dtau D^{\mathrm{JK}}_{ijkl} = & K(4J - 1)J \cdivY\Vcalbf\delta\ij\delta\kl
             -\frac{2}{3} \mu J^{-2/3} \cdivY\Vcalbf\bigg[
                \frac{2}{9}b_{mm} \delta\ij\delta\kl
               - \frac{2}{3}\left(b_{ij}\delta_{kl} + \delta_{ij}b_{kl}\right)\\
            &\quad  + \frac{1}{2}\left(\delta_{ik}b_{lj} + b_{il}\delta_{jk} + \delta_{li}b_{kj} + b_{ik}\delta_{jl}\right)
            \bigg]\\
            &+ \mu J^{-2/3} \left(\partial^y_m \Vcal_r b_{ms} + b_{rm} \partial^y_m \Vcal_s\right)\bigg[
                \frac{2}{9}\delta_{rs}\delta\ij\delta\kl\\
            &\quad  - \frac{2}{3}\left(\delta_{ir}\delta_{js}\delta_{kl} + \delta_{ij}\delta_{kr}\delta_{ls}\right)
                + \frac{1}{2}\left(\delta_{ik}\delta_{lr}\delta_{js}  + \delta_{ir}\delta_{ls}\delta_{jk}
                + \delta_{il}\delta_{kr}\delta_{js}  + \delta_{ir}\delta_{ks}\delta_{jl}\right)
                \bigg],\\
    \end{aligned}
\end{equation}
which subsequently allows the computation of the derivative of the Truesdell
rate as
\begin{equation}\label{eq:app2-dtk}
        \dtau D^{\mathrm{TK}}_{ijkl} = \dtau D^{\mathrm{JK}}_{ijkl}
        - \frac{1}{2}\left(
            \delta_{ik}\dtau\tau_{lj} +
            \dtau\tau_{il}\delta_{jk} +
            \delta_{il}\dtau\tau_{kj} + 
            \dtau\tau_{ik}\delta_{jl}
        \right).
\end{equation}
Since $\Dop = \Dop^{\mathrm{TC}} = \frac{1}{J} \Dop^{\mathrm{TK}}$ ($\Dop^{\mathrm{TC}}$ is the Truesdell rate of the Cauchy stress tensor) in
\eq{eq:weak_elastic_operator}, the shape sensitivity of $\tilde A_{ijkl} = D_{ijkl}^{\mathrm{TC}} + \tilde\sigma_{jl}\delta_{ik}$ is
finally obtained as
\begin{equation}\label{eq:app2-daop}
    \dtau \tilde A_{ijkl} = \frac{1}{J}\left(\dtau D^{\mathrm{TK}}_{ijkl} - D^{\mathrm{TK}}_{ijkl} \cdivY\Vcalbf\right) + \dtau\tilde\sigma_{jl}\delta_{ik}.
\end{equation}